\documentclass[11pt,a4paper]{article}
\usepackage{ifthen}
\usepackage{amsmath}
\usepackage{amssymb}
\usepackage{latexsym}
\usepackage{bbm}
\usepackage{graphicx}
\usepackage{pdfsync}
\usepackage{verbatim}
\usepackage{fancyhdr}
\usepackage{url}               
\usepackage[english]{babel}
\usepackage[applemac]{inputenc}

\usepackage{dsfont}

\usepackage[pdftex,%
a4paper=true,%
bookmarks=true,%
bookmarksnumbered=true,%
colorlinks=true,%
urlcolor=blue,%
pdfstartview=FitH%
]{hyperref}

\usepackage{tikz}

\newtheorem{theorem}{Theorem}[section]
\newtheorem{lemma}[theorem]{Lemma}

\newtheorem{proposition}[theorem]{Proposition}
\newtheorem{corollary}[theorem]{Corollary}
\newtheorem{definition}[theorem]{Definition}
\newtheorem{remark}[theorem]{\it \bf{Remark}\/}

\numberwithin{equation}{section}
\numberwithin{figure}{section}

\newtheorem{Complementary lemma}[theorem]{Complementary lemma}

\newcommand{\bp}{{\it Proof. }}
\newcommand{\ep}{\hfill $\square$\\}

\newcommand{\bpl}{{\it Proof of Lemma }} 
\newcommand{\epl}{\hfill $\square$\\}

\newcommand{\bpp}{{\it Proof of Proposition }} 
\newcommand{\epp}{\hfill $\square$\\}

\newcommand{\bpt}{{\it Proof of Theorem }}
\newcommand{\ept}{\hfill $\square$\\}

\renewcommand{\Re}{\mbox{Re}}
\newcommand{\be}{\begin{equation}}
\newcommand{\ee}{\end{equation}}
\newcommand{\bea}{\begin{eqnarray}}
\newcommand{\eea}{\end{eqnarray}}
\newcommand{\bee}{\begin{eqnarray*}}
\newcommand{\eee}{\end{eqnarray*}}

\def\pa{\partial}

\def\CC{\mathbb{C}}
\def\NN{\mathbb{N}}
\def\RR{\mathbb{R}}

\def\dps{\displaystyle}
\def\ni{\noindent}

\def\eps{\vare}

\def\eps{\varepsilon}

\catcode`@=11
\def\supess{\mathop{\operator@font Sup\,ess}}
\catcode`@=12
\def\CC{\mathbb{C}}

\def\NN{\mathbb{N}}

\def\RR{\mathbb{R}}
\def\CC{\mathbb{C}}

\def\ds{\displaystyle}

\def\ni{\noindent}

\def\bar#1{{\overline #1}}

\def\R2+{\RR ^2_+}

\def\pa{\partial}

\def\lim{\mathop{\rm lim}}

\def\sup{\mathop{\rm sup}}

\def\pa{\partial}

\def\pa{\partial}

\newcommand{\ud}{\mathrm{d}}

\newcommand{\ui}{\mathrm{i}}

\begin{document}

\renewcommand{\refname}{References}
\bibliographystyle{abbrv}

\pagestyle{fancy}
\fancyhead[L]{ }
\fancyhead[R]{}
\fancyfoot[C]{}
\fancyfoot[L]{ }
\fancyfoot[R]{}
\renewcommand{\headrulewidth}{0pt} 
\renewcommand{\footrulewidth}{0pt}

\newcommand{\montitre}{Fractional diffusion for Fokker-Planck equation with heavy tail equilibrium: an \`a la Koch spectral method in any dimension}

\newcommand{\auteur}{\textsc{Dahmane Dechicha$^\ast$,  Marjolaine Puel$^\ddagger$ }}
\newcommand{\affiliation}{$^\ast$Laboratoire J.-A. Dieudonn\'e.  Universit\'e C\^{o}te d'Azur. UMR 7351 \\ Parc Valrose, 06108 Nice Cedex 02,   France   \\

$^\ddagger$Laboratoire de recherche AGM. CY Cergy Paris Universit\'e. UMR CNRS 8088 \\

2 Avenue Adolphe Chauvin 95302 Cergy-Pontoise Cedex, France \\
\url{dechicha@unice.fr}; \url{ mpuel@cyu.fr
} }

 \begin{center}
{\bf  {\LARGE \montitre}}\\ \bigskip \bigskip
 {\large\auteur}\\ \bigskip \smallskip
 \affiliation \\ \bigskip
\today
 \end{center}
 \begin{abstract} 
 In this paper, we extend the spectral method developed \cite{DP} to any dimension $d\geqslant 1$, in order to construct an eigen-solution for the Fokker-Planck operator with heavy tail equilibria, of the form $(1+|v|^2)^{-\frac{\beta}{2}}$, in the range $\beta \in ]d,d+4[$. The method developed in dimension 1 was inspired by the work of H. Koch on nonlinear KdV equation \cite{Koch}. The strategy in this paper is the same as in dimension 1 but the tools are different, since dimension 1  was based on ODE methods. As a direct consequence of our construction, we obtain the fractional diffusion limit for the kinetic Fokker-Planck equation, for the correct density $\rho := \int_{\RR^d} f \mathrm{d}v$, with a fractional Laplacian $\kappa(-\Delta)^{\frac{\beta-d+2}{6}}$ and a positive diffusion coefficient $\kappa$.

\end{abstract}

\tableofcontents
\newpage

 \pagestyle{fancy}
\fancyhead[R]{\thepage}
\fancyfoot[C]{}
\fancyfoot[L]{}
\fancyfoot[R]{}
\renewcommand{\headrulewidth}{0.2pt} 
\renewcommand{\footrulewidth}{0pt} 

\section{Introduction}
\subsection{Setting of the problem}
In this present paper, we deal with the kinetic Fokker-Planck (FP) equation, which describes in a deterministic way the Brownian motion of a set of particles. It is given by the following form
\begin{equation}\label{FP}
\left\{\begin{array}{l}
 \partial_t f + v\cdot \nabla_x f = Q(f), \; \; \; t \geqslant 0 , \ x \in \RR^d , \ v \in \RR^d ,\\
\\
    f(0,x,v) = f_0(x,v), \; \; \quad x \in \RR^d , \ v \in \RR^d ,
\end{array}\right.
\end{equation}
where the collisional Fokker-Planck operator $Q$ is given by
\begin{equation}\label{defQ}
Q(f)=\nabla_v\cdot\bigg(F\nabla_v\bigg(\frac{f}{F}\bigg)\bigg) ,
\end{equation}
and $F$ is the equilibrium of $Q$, a fixed function which depends only on $v$ and satisfying  
$$ Q(F)=0 \quad \mbox{ and } \quad \int_{\RR^d} F(v)\ \mathrm{d}v = 1 . $$ 
Provided $f_0 \geqslant 0$, the unknown $f(t,x,v) \geqslant 0$ can be interpreted as the density of particles occupying at time $t \geqslant 0$, the position $x \in \RR^d$ with velocity $v \in \RR^d$. \\

Recall that one of the motivations for studying the \emph{classical} or \emph{fractional diffusion limit} is to simplify the equations for some collisional kinetic models when the interaction between particles are the dominant phenomena and when the observation time is very large. For that purpose, we introduce a small parameter, $\eps \ll 1$, the mean free path and we proceed to rescaling the distribution function $f(t,x,v)$ in time and space
$$ t=\frac{t'}{\theta(\eps)} \quad \mbox{ and } \quad x=\frac{x'}{\eps} \quad \mbox{ with } \quad \theta(\eps) \underset{\eps \to 0}{\longrightarrow } 0 , $$
which leads to the following rescaled equation (without primes)
\begin{equation}\label{fp-theta}
\left\{\begin{array}{l}
\theta(\varepsilon) \partial_t f^\varepsilon + \varepsilon v\cdot \nabla_x f^\varepsilon = Q(f^\varepsilon), \; \; \;  t \geqslant 0 , x \in \RR^d , v \in \RR^d , \\
\\
    f^\varepsilon(0,x,v) = f_0(x,v), \; \; \; \; x \in \RR^d , v \in \RR^d .
\end{array}\right.
\end{equation}
Note that initial condition written in non rescaled variable are well prepared conditions. \\

The goal is then to study the behavior of the solution $f^\eps$ as $\eps \to 0$. Formally, passing to the limit when $\eps \to 0$ in equation \eqref{fp-theta}, we obtain that the limit $f^0$ is in the kernel of $ Q$ which is spanned by the equilibrium $ F $, which means that $f^0=\rho(t,x)F(v)$. Thus, it amounts to  find the equation satisfied by the density $\rho$. Note that this limit depends on the nature of the equilibrium $F$ considered as well as on the chosen change of time scale $\theta(\eps)$. \\

For Gaussian equilibria, it is classical (see \cite{BaSaSe},\cite{BeLi},\cite{DeGoPo},\cite{LaKe},\cite{De2} for Boltzmann and \cite{DeMa-Ga} for Fokker Planck) that by taking the classical time scaling $\theta(\eps)=\eps^2$, we obtain  a \emph{diffusion equation}
\begin{equation}\label{d}
\partial_t\rho-\nabla_x\cdot(D\nabla_x\rho)=0 ,
\end{equation}
where 
\begin{equation}\label{dc}
D=\int vQ^{-1}(-vF) \mathrm{d}v.
\end{equation}

 For slowly decreasing equilibria, or so-called \emph{heavy tail} equilibria of the form $F(v) \sim \langle v \rangle^{-\beta}$, it is more complicated, and this study has been the interest of many papers in the last few years, with different methods and for different collision operators. Fractional diffusion limit has been obtained in the case of the linear Boltzmann equation when the cross section is such that the operator has a spectral gap, see \cite{MMM} for the pioneer paper in the case of space independent cross section, where the authors used a method based on Fourier-Laplace transformation, and see \cite{M} for a weak convergence result obtained by the Moment method, which also applies to cross sections that depend on the position variable. See also \cite{MKO} for a probabilistic approach.  \\
 
 In the present work, we consider for any $\beta > d$, heavy tail equilibria 
$$ F(v) = \frac{C_\beta^2}{(1+|v|^2)^\frac{\beta}{2}} ,$$ 
where $C_\beta$ is a normalization constant. \\

The diffusion limit for the FP equation seems more complicated then the linear Boltzmann one, and the main difficulty is due to the fact that the Fokker-Planck operator $Q$ has no spectral gap. In addition, for this equation, all the terms of the operator participate in the limit, i.e. the collision and advection parts. In \cite{NaPu}, the classical scaling is studied and it is proved in any dimension $d$ that we obtain a diffusion equation \eqref{d}, with diffusion coefficient \eqref{dc} as soon as $\beta>d+4$. The critical case where $\beta=d+4$ is studied in \cite{CNP}, where the expected result of classical diffusion with an anomalous time scaling is proved, $\theta(\eps)=\eps^2|\ln \eps|$. A unified presentation of the result for even more general cases of $\beta$ can be found in recent papers where the result has been obtained, by probabilistic method in  \cite{FT-d1} and  \cite{FT}, and using a \emph{quasi-spectral} problem in \cite{BM}. In this last paper, in addition to the diffusion limit results, an estimates on the fluid approximation error have been obtained. We refer also to \cite{BDL} for this last point, where the authors have developed an $L^2$-hypocoercivity approach and established an optimal decay rate, determined by a fractional Nash type inequality, compatible with the fractional diffusion limit.\\

 In this paper we focus on the case $d<\beta<d+4$. By taking as test function the eigenvector of the whole Fokker-Planck operator (advection $+$ collisions), which converges towards equilibrium $F$, we capture at the limit the ``diffusion'' equation for any $\beta>d$. The computation of the eigenvalue  gives us the right scaling in time, $\theta(\eps)$, and the diffusion coefficient $\kappa$ at the same time. We are therefore interested in a new problem: \emph{the construction of an eigen-solution for the whole Fokker-Planck operator}, which is the main subject of this paper. \\
 
This spectral problem for the FP operator has already been obtained recently in dimension 1 \cite{LebPu} with a method based on the reconnection of two branches on $\RR_+$ and $\RR_-$, but this  method of reconnection is difficult to adapt in dimension $d$. This led us to look for another strategy, which was the subject of \cite{DP}, a method inspired by the work of H. Koch on nonlinear KdV equation \cite{Koch}, which allowed us to construct an eigen-solution for the spectral problem associated to the whole Fokker-Planck operator with ODE methods in dimension 1. The aim of this paper is to develop PDE methods in order to obtain the result in any dimension. This method is interesting since it can be used for different potentials like convolution, or  for nonlinear equations as well. Moreover, as in dimension 1, a splitting of the Fokker-Planck operator is involved, which recalls the \emph{enlargement theory} for nonlinear Boltzmann operator when there are spectral gap issues. This theory was developed by Gualdani, Mischler and Mouhot in \cite{GMM} whose key idea was based on the decomposition of the operator into two parts, a dissipative part plus a regularizing part. See also \cite{Gervais} and references therein.  \\

Note that we don't look at the same spectral problem as in the paper by E. Bouin and C. Mouhot \cite{BM}.  Indeed,  in this paper we were interested in the improvement and generalization of the construction given in \cite{LebPu} to solve the problem
$$ \big[ Q + \ui \eps \xi \cdot v \big] M_{\mu,\eps} = \mu \ M_{\mu,\eps} ,$$
with $\xi$ being the Fourier variable of $x$. While in \cite{BM} the authors considered the following quasi-spectral problem:
\begin{equation}\label{Bouin-Mouhot problem}
  \big[ Q + \ui \eps \xi \cdot v \big] \phi_{\mu,\eps} = \mu \ \frac{\phi_{\mu,\eps}}{\langle v \rangle^{2}} ,
\end{equation}
with $\phi_{\mu,\eps} \in L^2 \big(\RR^d ; \frac{\ud v}{\langle v \rangle^2}\big)$ satisfying $\ds \int_{\RR^d} \phi_{\mu,\eps}(v)M(v) \ \frac{\ud v}{\langle v \rangle^2} =1 $.
The key idea in \eqref{Bouin-Mouhot problem} is the introduction of a weight that allowed to recover the spectral gap inequality for the latter operator thanks to the Hardy-Poincar\'e inequality
$$ \int_{\RR^d} f Q(f) \ \ud v \geqslant C \int_{\RR^d} |f-r M |^2 \ \frac{\ud v}{\langle v \rangle^2} ,$$
where $r$ is a weighted density defined by
\begin{equation}\label{densite r}
  r(t,x) := \int_{\RR^d} f \ \frac{\ud v}{\langle v \rangle^2} .
\end{equation}
Thus,  by totally different techniques based on energy estimates and the study of the resolvent,  E. Bouin and C. Mouhot showed the existence of a ``fluid mode'',  a couple $ \big(\mu(\eps),\phi_{\mu,\eps}\big)$ solution of  problem \eqref{Bouin-Mouhot problem}. Thanks to this construction, they obtain the convergence of $f^\eps/F$ towards $(\int_{\RR^d}\frac{F}{\langle v \rangle^2}\ud v)^{-1}r(t,x)$ in $L^2_t\big([0,T];H^{-\frac{\beta-d+2}{3}}_xL^2_v(\frac{F}{\langle v\rangle^2})\big)$, when $\eps$ goes to $0$,  with $r$ solution to a fractional diffusion equation.  Finally,  the diffusion limit with the \emph{classical} density $\rho := \int f  \ud v$ is recovered, in a weak sense.

\subsection{Setting of the result}
Before stating our main result, let us give some notations that we will use along this paper.\\

\noindent \textbf{Notations.} 
As in \cite{LebPu}, in order to simplify the computation and work with a self-adjoint operator in $L^2$, we proceed to a change of unknown by writing
$$ f = F^\frac{1}{2}g = C_\beta M g $$
with  $$ M :=C_\beta^{-1} F^\frac{1}{2}=\frac{1}{(1+|v|^2)^\frac{\gamma}{2}} , $$ 
since we impose $\gamma :=\frac{\beta}{2}>\frac{d}{2}$, $F \in L^1(\RR^d)$ then, $M \in L^2(\RR^d)$ and $C_\beta$ is chosen such that $$ \int_{\RR^d} F \mathrm{d}v = 1 .$$
The equation \eqref{fp-theta} becomes
$$ \theta(\eps) \partial_t g^\eps + \eps v\cdot \nabla_x g^\eps =\frac{1}{M}\nabla_v \cdot  \bigg(M^2\nabla_v\bigg(\frac{g^\eps}{M}\bigg)\bigg) = \Delta_v g^\eps - W(v)g^\eps ,
$$
with 
$$ W(v)= \frac{\Delta_vM}{M} = \frac{\gamma(\gamma-d+2)|v|^2 -\gamma d}{(1+|v|^2)^2}.
$$
We see the equation as 
$$
\theta(\eps) \partial_t g^\eps = - \mathcal{L}_\eps g^\eps ,
$$
where
$$ \mathcal{L}_\eps := -\Delta_v + W(v) + \eps v \cdot \nabla_x = -(Q-\eps v \cdot \nabla_x) $$
and $$ Q  :=  -\Delta_v + W(v) .$$

\noindent We operate a Fourier transform in $x$ and since the operator $Q$ has coefficient that do not depend on $x$, we get:
\begin{equation}\label{rescaled}
\theta(\eps) \partial_t \hat g^\eps = - \mathcal{L}_\eta \hat g^\eps ,
\end{equation}
where
$$ \mathcal{L}_\eta := -\Delta_v + W(v) + \mathrm{i} \eta v_1 $$
with $$ \eta := \eps |\xi| \quad \mbox{ and } \quad v_1 := v \cdot \frac{\xi}{|\xi|} ,$$ 
where $\xi$ being the space Fourier variable. \\

\noindent The operator $\mathcal{L}_\eta$ is an unbounded self-adjoint operator acting on $L^2$. Its domain is given by
$$D(\mathcal{L}_\eta) = \big\{  g \in L^2(\RR^d) \ ; \ \Delta_v g \in L^2(\RR^d), v_1 g \in L^2(\RR^d)\big\} . $$

\subsubsection*{Main results.}
\begin{theorem}[Eigen-solution for the Fokker-Planck operator]\label{main}
\item Assume that $d<\beta<d+4$ with $\beta\neq d+1$. Let $\eta_0>0$ and $\lambda_0>0$ small enough.  Then, for all $\eta\in [0,\eta_0]$,  there exists a unique eigen-couple $\big(\mu(\eta),M_\eta\big)$ in $\{\mu\in \mathbb C, |\mu|\leqslant \eta^{\frac{2}{3}}\lambda_0\}\times L^2(\RR^d,\CC)$, solution to the spectral problem 
\begin{equation}\label{M_mu,eta}
\mathcal{L}_\eta(M_{\mu,\eta})(v)=\big[-\Delta_v +W(v) + \mathrm{i} \eta v_1 \big] M_{\mu,\eta} (v)= \mu M_{\mu,\eta}(v) , \  v\in \RR^d .
\end{equation}
Moreover,\\
1. The following convergence in the Sobolev space $H^1(\RR^d)$ holds:  
\begin{equation}\label{M_eta--->M dans H^1}
 \big\|M_\eta - M \big\|_{H^1(\RR^d)} \underset{\eta \to 0}{\longrightarrow} 0 .
\end{equation}
2. The  eigenvalue $\mu(\eta)$ is given by 
\begin{equation}\label{glres}
\mu(\eta)=\bar\mu(-\eta)= \kappa |\eta|^{\frac{\beta-d+2}{3}}\big(1+O(|\eta|^{\frac{\beta-d+2}{3}})\big) ,
\end{equation}
where $\kappa$ is a positive constant given by 
\begin{equation}\label{kappa thm}
\kappa= -2C_\beta^2 \int_{\{s_1>0\}} s_1|s|^{-\gamma} \mathrm{Im} H_0(s) \mathrm{d}s ,
\end{equation}
and where $H_0$ is the unique solution to the equation 
\begin{equation}
\big[-\Delta_s+\frac{\gamma(\gamma-d+2)}{|s|^2} + \mathrm{i} s_1 \big]H_0(s) = 0 , \quad \forall s \in \RR^d \setminus \{0\} ,
\end{equation}
satisfying 
\begin{equation}
\int_{\{|s_1|\geqslant 1 \}} |H_0(s)|^2 \mathrm{d}s < \infty \quad \mbox{ and } \quad H_0(s) \underset{0}{\sim} |s|^{-\gamma} .
\end{equation}
\end{theorem}

Introduce  $V$,  the space defined by
$$ V := \left\{f: \mathbb{R}^d\rightarrow \mathbb{R},\int_{\mathbb{R}^{d}} \frac{|f|^2}{F} \ \mathrm{d}v< \infty \mbox{ and } \int_{\mathbb{R}^{d}}  \bigg|\nabla_v\bigg(\frac{f}{F}\bigg)\bigg|^2  F \mathrm{d}v < \infty \right\} , $$ 
$V'$ being its dual, and 
$$ Y := \left\{ f \in L^2\big([0,T]\times\RR^d; V\big); \ \theta(\eps) \pa_t f + \eps v \cdot \nabla_x f \in L^2\big([0,T]\times\RR^d; V'\big) \right\} . $$
\begin{theorem}[Fractional diffusion limit for the Fokker-Planck equation]\label{main2}
\item Assume that $d<\beta<d+4$ with $\beta\neq d+1$. Assume that  $f_0\in L^1(\RR^{2d})$ is a non-negative function in  $ L^2_{F^{-1}}(\RR^{2d})\cap L^\infty_{F^{-1}}(\RR^{2d})$. Let $f^\varepsilon$ be the solution of \eqref{fp-theta} in $Y$ with initial data $f_0$, with $\theta(\varepsilon)=\varepsilon^{\frac{\beta-d+2}{3}} $. Let $ \kappa$ be the constant given by \eqref{kappa thm}.\\
Then $f^\varepsilon$ converges weakly star in $L^\infty\big([0,T],\  L^2_{F^{-1}}(\RR^{2d})\big)$ towards $\rho(t,x) F(v)$ where $\rho(t,x)$ is the solution to 
\begin{equation}
\partial_t\rho +\kappa (-\Delta)^{\frac{\beta-d+2}{6}}\rho =0,\quad \rho(0,x)=\int_{\RR^d} f_0 \mathrm{d}v .
 \end{equation}
\end{theorem}
\begin{remark}The hypothesis $\beta\neq d+1$ is  technical. It avoids to introduce logarithmic terms in the expression of $\mu(\eta)$. 
\end{remark}

\subsubsection*{Ideas of the proof and outline of the paper.}
The proof of Theorem \ref{main} is done in two main steps, both based on the implicit function theorem (IFT).
First, we consider what we call a \emph{penalized equation}, given by
\begin{equation}\label{eq penalisee1}
\left\{ \begin{array}{l}
\big[-\Delta_v + W(v) + \mathrm{i} \eta v_1 \big] M_{\mu,\eta} (v)= \mu M_{\mu,\eta}(v) - \langle M_{\mu,\eta}-M,\Phi\rangle \Phi(v) , \  v\in \mathbb{R}^d ,\\
\\
M_{\mu,\eta} \in L^2(\mathbb{R}^d).
\end{array}\right.
\end{equation}
where $\Phi$ is a function, that satisfies some assumptions, that we will determine later.  The additional term allows us to avoid the problem of reconnection by ensuring existence of a solution to equation \eqref{eq penalisee1} on the whole space $\RR^d$ for any $\eta$ and $\mu$. This is one of the key points of this method. Also, note that the sign before the scalar product $\langle M_{\mu,\eta}-M,\Phi\rangle$ is important.\\

The aim of the first step is to show the existence of a unique solution for equation \eqref{eq penalisee1} for $\eta$ and $\mu$ fixed, which is the purpose of Section 2. As we said above, we will decompose the operator $``-\Delta_v + W(v) + \mathrm{i} \eta v_1 - \mu"$ in two parts. The first one is chosen such that it admits an inverse that is continuous as a linear operator between two suitable functional spaces, continuous with respect to the parameters $\eta$ and $\mu$ and compact at $\eta=\mu=0$. The second part of the operator is left in the right-hand side of the equation, i.e. is considered as a source term. The invertibility of the first part is the subject of the first subsection, and it is based on an elaborated version of the Lax-Milgram theorem. While the study of the inverse operator and its properties is the subject of the second subsection whose main result is the existence of solutions for equation \eqref{eq penalisee1}.  \\

In the second step, to ensure that the additional term vanishes, we have to chose $\mu(\eta)$ obtained via the implicit function theorem around the point $(\mu,\eta)=(0,0)$. The study of this constraint is the subject of a large part of section 3 which is composed of three subsections. The first one is dedicated to the $L^2$ estimates for the solution of the penalized equation \eqref{eq penalisee1}. It consists in improving the space to which the solution found by Lax-Milgram belongs.  It is the objective of the second subsection. The last subsection is dedicated to the approximation of the eigenvalue and the computation of the diffusion coefficient. \\

The last section is devoted to the proof of Theorem \ref{main2}. It consists of two subsections, a priori estimates and limiting process in the weak formulation of equation \eqref{rescaled}.
\section{Existence of solutions for the penalized equation}
We start this section by some notations and definition of the considered operators. Let $\mu =\lambda\eta^\frac{2}{3}$ with $\lambda \in \mathbb{C}$ and let denote by $L_{\lambda,\eta}$ the operator
$$ L_{\lambda,\eta} := -\Delta_v + \tilde W(v) + \mathrm{i} \eta v_1 -  \lambda \eta^\frac{2}{3} ,$$
where $$\tilde W(v) :=\frac{\gamma(\gamma-d+2)}{1+|v|^2}.$$
Let denote by $V := \tilde W - W$.  We have
$$ V(v)=\frac{\gamma(\gamma+2)}{(1+|v|^2)^2}.$$
We will rewrite  equation \eqref{eq penalisee1} as follows
\begin{equation}\label{eq penalisee2}
\left\{ \begin{array}{l}
L_{\lambda,\eta} (M_{\lambda,\eta})  =  V(v) M_{\lambda,\eta} - \langle M_{\lambda,\eta}-M,\Phi\rangle \Phi , \ v\in \RR^d ,\\
\\
M_{\lambda,\eta} \in L^2(\RR^d).
\end{array}\right.
\end{equation}
The two equations \eqref{eq penalisee1} and \eqref{eq penalisee2} are equivalent.
\begin{remark}\label{symetrie de la sol}
\item \begin{enumerate}
\item Since $L_{\lambda,0}$ does not depend on $\lambda$, let's denote it by $L_0$, $L_0:=L_{\lambda,0}$.
\item If $ \ \bar \Phi(-v)=\Phi(v)$ and $M_{\lambda,\eta}(v_1,v')$ satisfies the equation \eqref{eq penalisee2}, then $\overline{M}_{\bar \lambda,\eta}(-v_1,v')$ satisfies also \eqref{eq penalisee2}, since the potential $ W$ is symmetric for a symmetric equilibrium $M$. Note that this is where the symmetry of the equilibrium $M$ is used and therefore this is a ``non-drift condition". 
\item Note that the splitting of the potential $W$ into $\tilde W$ and $V$ is crucial in our study. It plays a very important role whether in the invertibility of the operator $L_{\lambda,\eta}$ or in the compactness of its inverse at the point $(\lambda,\eta)=(0,0) $. 
\end{enumerate}
\end{remark}

\subsection{Coercivity and Lax-Milgram theorem}
The purpose of this subsection is to show that the operator $L_{\lambda,\eta}$ defined above is invertible. For this, we are going to define a Hilbert space $\mathcal{H}_\eta$ as well as a scalar product $\langle \cdot,\cdot\rangle_{\mathcal{H}_\eta}$ on which we apply a Lax-Milgram theorem. 
\begin{definition} 
\item $\bullet$ We define the Hilbert space $\mathcal{H}_\eta$ as being the completion of the space $C_c^\infty(\RR^d,\CC)$ for the norm $\|\cdot\| _{\mathcal{H}_\eta}$ induced from the scalar product $\langle \cdot,\cdot\rangle_{\mathcal{H}_\eta}$
$$
\mathcal{H}_\eta := \overline{\big\{ \psi \in C_c^\infty(\RR^d,\CC) ; \ \| \psi \|_{\tilde{\mathcal{H}}_\eta}^2:= \langle \psi,\psi\rangle_{\tilde{\mathcal{H}}_\eta} < +\infty \big\} } ,
$$
where
$$
\langle \psi,\phi\rangle_{\tilde{\mathcal{H}}_\eta}:=  \int_{\RR^d} \nabla_v \bigg(\frac{\psi}{M}\bigg) \cdot \nabla_v \bigg(\frac{\bar\phi}{M}\bigg) \ M^2 \mathrm{d}v + \int_{\RR^d}  V \psi \bar \phi  \ \mathrm{d}v +  \eta \int_{\RR^d} |v_1| \psi \bar \phi \ \mathrm{d}v ,
$$
and where $\ V(v):= \tilde W(v)-W(v) = \frac{\gamma(\gamma+2)}{(1+|v|^2)^2} > 0 \ $ for all $\ v \in \RR^d$.  \\

\noindent We have the embeddings
$$ \mathcal{H}_\eta \subseteq \mathcal{H}_{\eta^*} \subseteq \mathcal{H}_{0}  ,   \quad \forall \ 0 \leqslant \eta^* \leqslant \eta $$
since $\ \| \cdot \|_{\tilde{\mathcal{H}}_0} \leqslant \| \cdot \|_{\tilde{\mathcal{H}}_{\eta^*}} \leqslant \| \cdot \|_{\tilde{\mathcal{H}}_\eta} \ $  for all $\ 0\leqslant \eta^* \leqslant \eta$.
\\

\noindent $\bullet$ We define the sesquilinear form $a$ on $\mathcal{H}_\eta\times\mathcal{H}_\eta$ by
$$
a(\psi,\phi):= \int_{\RR^d} \nabla_v \bigg(\frac{\psi}{M}\bigg) \cdot \nabla_v \bigg(\frac{\bar\phi}{M}\bigg) \ M^2 \mathrm{d}v + \int_{\RR^d}  V \psi \bar \phi  \ \mathrm{d}v +  \mathrm{i} \eta \int_{\RR^d} v_1 \psi \bar \phi \ \mathrm{d}v  - \lambda \eta^{\frac{2}{3}} \int_{\RR^d} \psi \bar \phi \ \mathrm{d}v  .
$$
\end{definition}
\begin{remark}\label{tilde Q = Q + V}
\item \begin{enumerate}
\item Note that $a(\psi,\psi)\neq \|\psi\|_{\tilde{\mathcal{H}}_\eta}^2$.
\item Note that the sesquilinear form $a$ depends on $\lambda$ and $\eta$ and in order to simplify the notation, we omit the subscript when no confusion is possible.
\item Let us denote by $\tilde Q$ the operator $\tilde Q := -\Delta_v + \tilde W(v)$. We have $ \tilde Q = Q + V $. Thus, the operator $\tilde Q$ is dissipative since
$$  \int_{\RR^d} \tilde Q(\psi)\psi \ \mathrm{d}v = \int_{\RR^d} Q(\psi)\psi \ \mathrm{d}v + \int_{\RR^d} V |\psi|^2 \ \mathrm{d}v =  \int_{\RR^d} \bigg| \nabla_v \bigg(\frac{\psi}{M}\bigg)\bigg|^2 M^2 + V |\psi|^2 \ \mathrm{d}v\geqslant 0 . $$
Note that we have also the equality
$$ \int_{\RR^d} \tilde Q(\psi)\psi \ \mathrm{d}v  = \int_{\RR^d} \big| \nabla_v \psi \big|^2 \ \mathrm{d}v + c_{\gamma,d} \int_{\RR^d} \frac{|\psi |^2}{\langle v\rangle^2} \ \mathrm{d}v , $$
with $\ c_{\gamma,d}:=\gamma(\gamma-d+2)$. Observe that $\ c_{\gamma,d} < 0 \ $ for $\ \gamma \in (\frac{d}{2},\frac{d+4}{2}) \ $ with $\ d> 4$.  
\item Since $\ \tilde Q = Q + V \ $ then, the sesquilinear form $a$ can be written as follows:
$$ a(\psi,\phi) = \int_{\RR^d} \nabla_v \psi \cdot \nabla_v \bar\phi \ \mathrm{d}v + c_{\gamma,d} \int_{\RR^d} \frac{\psi \bar \phi}{\langle v\rangle^2} \ \mathrm{d}v +  \mathrm{i} \eta \int_{\RR^d} v_1 \psi \bar \phi \ \mathrm{d}v - \lambda \eta^{\frac{2}{3}} \int_{\RR^d} \psi \bar \phi \ \mathrm{d}v . $$
\end{enumerate}
\end{remark}
\begin{lemma}\label{equivalence des normes}
The norm defined by 
$$ \| \psi \|_{{\mathcal{H}}_\eta}^2 := \int_{\RR^d} \big |\nabla_v \psi \big|^2 \ \mathrm{d}v +  \int_{\RR^d} \frac{|\psi|^2}{\langle v\rangle^2} \ \mathrm{d}v +  \eta \int_{\RR^d} |v_1| |\psi|^2 \ \mathrm{d}v $$
is induced from the scalar product
 $$
\langle \psi,\phi\rangle_{{\mathcal{H}}_\eta}:=  \int_{\RR^d} \nabla_v \psi \cdot \nabla_v \bar\phi \ \mathrm{d}v +  \int_{\RR^d} \frac{\psi \bar \phi}{\langle v\rangle^2} \ \mathrm{d}v +  \eta \int_{\RR^d} |v_1| \psi \bar \phi \ \mathrm{d}v ,
$$
and the two norms $\| \cdot \|_{\mathcal{H}_{\eta}} $ and $\| \cdot \|_{\tilde{\mathcal{H}}_\eta} $ are equivalent, i.e., there are two positive constants $C_1$ and $C_2$ such that
$$ C_1 \| \psi \|_{\tilde{\mathcal{H}}_\eta} \leqslant \| \psi \|_{{\mathcal{H}}_\eta} \leqslant C_2 \| \psi \|_{\tilde{\mathcal{H}}_\eta} , \quad \forall \psi \in \mathcal{H}_\eta . $$
\end{lemma}
To prove this Lemma, we need the Hardy-Poincar\'e inequality that we recall in the following 
\begin{lemma}[Hardy-Poincar\'e inequality]\label{Hardy-Poincare} $\cite{BDGV}$
Let $d\geqslant 1$ and $\alpha_*=\frac{2-d}{2}$. For any $\alpha<0$,  and $\alpha \in (-\infty,0)\setminus \{ \alpha_*\}$ for $d\geqslant 3$, there is a positive constant $\Lambda_{\alpha,d}$ such that
\begin{equation}
\Lambda_{\alpha,d} \int_{\RR^d} |f|^2(D+|x|^2)^{\alpha-1}\mathrm{d}x \leqslant \int_{\RR^d} |\nabla f|^2 (D+|x|^2)^{\alpha}\mathrm{d}x
\end{equation}
holds for any function $f \in H^1\big((D+|x|^2)^{\alpha}\mathrm{d}x\big)$ and any $D \geqslant 0$, under the additional condition $\int_{\RR^d} f (D+|x|^2)^{\alpha-1}\mathrm{d}x = 0$ and $D>0$ if $\alpha<\alpha_*$.
\end{lemma}
\begin{remark}
For $\ds f=\frac{g}{M}$, $D=1$ and $\alpha = -\gamma$ in the previous lemma, the inequality becomes
\begin{equation}\label{ineg.  Hardy-Poincare}
\Lambda_{\alpha,d} \int_{\RR^d} \frac{|g|^2 }{\langle v\rangle^2} \mathrm{d}v \leqslant \int_{\RR^d} \bigg|\nabla_v \bigg(\frac{g}{M}\bigg)\bigg|^2 M^2 \ \mathrm{d}v ,
\end{equation}
and the orthogonality condition becomes
\begin{equation}\label{condition d'orthogonalite}
 \int_{\RR^d}  \frac{g M}{\langle v \rangle^2} \mathrm{d}v = 0 
\end{equation}
since $\  -\gamma < \frac{2-d}{2} =: \alpha_* \ $ for $\ \gamma \in (\frac{d}{2},\frac{d+4}{2})$.  \\
If we denote by
$$ \mathcal{P}(g) := \left(\int_{\RR^d} \frac{M^2}{\langle v \rangle^2} \mathrm{d}v\right)^{-1} \int_{\RR^d}  \frac{g M}{\langle v \rangle^2} \mathrm{d}v .$$
Then,  inequality \eqref{ineg. Hardy-Poincare} can be written for all $g \in \mathcal{H}_0$
\begin{equation}\label{Hardy-Poincare avec orthogonalite}
\Lambda_{\alpha,d} \int_{\RR^d}  \frac{\big| g -\mathcal{P}(g)M \big|^2}{\langle v\rangle^2} \mathrm{d}v \leqslant \int_{\RR^d} \bigg|\nabla_v \bigg(\frac{g}{M}\bigg)\bigg|^2 M^2 \ \mathrm{d}v .
\end{equation}
\end{remark}
\bpl \ref{equivalence des normes}.  Let's start with the right inequality: $ \| \psi \|_{{\mathcal{H}}_\eta} \leqslant C_2 \| \psi \|_{\tilde{\mathcal{H}}_\eta} $.   Let $\psi \in \mathcal{H}_\eta$.  
Then, since $M \in L^2(\RR^d)$, by Cauchy-Schwarz inequality we get
$$\bigg| \int_{\RR^d} \frac{\psi M}{\langle v \rangle^2} \mathrm{d}v \bigg| \leqslant \left(\int_{\RR^d} \frac{|\psi|^2}{\langle v \rangle^4} \mathrm{d}v\right)^{\frac{1}{2}}  \left(\int_{\RR^d} M^2 \mathrm{d}v\right)^{\frac{1}{2}} \leqslant \frac{1}{\gamma(\gamma+2)} \| \psi \|_{\tilde{\mathcal{H}}_\eta} . $$
Now, since the function $\psi - \mathcal{P}(\psi)M$ satisfies  condition \eqref{condition d'orthogonalite}, $ \mathcal{P}\big(\psi-\mathcal{P}(\psi)M\big)=0$,  then  inequality \eqref{ineg.  Hardy-Poincare} can be used and therefore
\begin{align*}
\int_{\RR^d} \frac{|\psi|^2}{\langle v \rangle^2} \mathrm{d}v &= \int_{\RR^d} \frac{|\psi- \mathcal{P}(\psi)M + \mathcal{P}(\psi)M|^2}{\langle v \rangle^2} \mathrm{d}v \\ 
&\leqslant 2 \left( \Lambda_{\alpha,d}^{-1} \int_{\RR^d} \bigg|\nabla_v \bigg(\frac{\psi}{M}\bigg)\bigg|^2 M^2 \ \mathrm{d}v+ |\mathcal{P}(\psi)|^2 \int_{\RR^d} \frac{M^2}{\langle v\rangle^2}\mathrm{d}v \right)  \\
&\leqslant 2 \left( \Lambda_{\alpha,d}^{-1}  \|\psi\|_{\tilde{\mathcal{H}}_\eta}^2 + \frac{1}{\gamma^2(\gamma+2)^2}\left(\int_{\RR^d} \frac{M^2}{\langle v \rangle^2} \mathrm{d}v\right)^{-1}  \| \psi \|_{\tilde{\mathcal{H}}_\eta}^2 \right)   \\
&\leqslant C_{\gamma,d}  \|\psi\|_{\tilde{\mathcal{H}}_\eta}^2 .
\end{align*}
We have by the first point of Remark \ref{tilde Q = Q + V}
$$  \int_{\RR^d} \tilde Q(\psi)\psi \ \mathrm{d}v =  \int_{\RR^d} \bigg| \nabla_v \bigg(\frac{\psi}{M}\bigg)\bigg|^2 M^2 + V |\psi|^2 \ \mathrm{d}v = \int_{\RR^d} \big| \nabla_v \psi \big|^2 \ \mathrm{d}v + c_{\gamma,d} \int_{\RR^d} \frac{|\psi |^2}{\langle v\rangle^2} \ \mathrm{d}v . $$
From where we get 
$$  \int_{\RR^d} \big| \nabla_v \psi \big|^2 \ \mathrm{d}v \leqslant \int_{\RR^d} \bigg| \nabla_v \bigg(\frac{\psi}{M}\bigg)\bigg|^2 M^2 + V |\psi|^2 \ \mathrm{d}v  +  |c_{\gamma,d}| \int_{\RR^d} \frac{|\psi |^2}{\langle v\rangle^2} \ \mathrm{d}v  \leqslant \big(1+\tilde C_{\gamma,d}\big)  \|\psi\|_{\tilde{\mathcal{H}}_\eta}^2 . $$ 
Hence, $$\|\psi\|_{{\mathcal{H}}_\eta} \leqslant C_2 \|\psi\|_{\tilde{\mathcal{H}}_\eta} , $$
with $C_2 := \sqrt{2(1+\tilde C_{\gamma,d})}$ a positive constant which depends only on $\gamma$ and $d$. \\
To get inequality $ C_1 \| \psi \|_{\tilde{\mathcal{H}}_\eta} \leqslant  \| \psi \|_{{\mathcal{H}}_\eta} $, it is enough just to write
 $$\int_{\RR^d} \bigg| \nabla_v \bigg(\frac{\psi}{M}\bigg)\bigg|^2 M^2 + V |\psi|^2 \ \mathrm{d}v = \int_{\RR^d} \big| \nabla_v \psi \big|^2 \ \mathrm{d}v + c_{\gamma,d} \int_{\RR^d} \frac{|\psi |^2}{\langle v\rangle^2} \ \mathrm{d}v \leqslant (1+|c_{\gamma,d}|) \| \psi \|_{{\mathcal{H}}_\eta}^2 . $$
 Hence, $$C_1 \|\psi\|_{\tilde{\mathcal{H}}_\eta} \leqslant \|\psi\|_{{\mathcal{H}}_\eta} , $$
with $C_1:= (2+|c_{\gamma,d}|)^{-\frac{1}{2}}$ a positive constant which depends only on $\gamma$ and $d$. 
\epl

\noindent In the remainder of this section, we will work with the norm $\| \cdot \|_{{\mathcal{H}}_\eta}$. \\

\noindent Before moving on to the continuity of $a$, we will prove a Poincar\'e type inequality which we give in the following lemma:
\begin{lemma}\label{Poincare1} Let $\eta>0$ be fixed. Then, there exists a constant $C_0>0$, independent of $\eta$ such that the following inequality holds true
$$
\| \psi \|_{L^2(\RR^d)} \leqslant C_0 \eta^{-\frac{1}{3}} \| \psi \|_{\mathcal{H}_\eta}  , \quad \forall \psi \in \mathcal{H}_\eta .
$$
\end{lemma}
\bp We will split the integral of $\| \psi \|_{L^2(\RR^2)}^2$ into two parts $\{ |v_1| \leqslant \eta^{-\frac{1}{3}}\}$ and $\{ |v_1| \geqslant \eta^{-\frac{1}{3}}\}$.\\
$\bullet$ On $\{ |v_1| \geqslant \eta^{-\frac{1}{3}}\}$, we simply have 
$$ \eta^{\frac{2}{3}} \int_{\{ |v_1| \geqslant \eta^{-\frac{1}{3}}\}} |\psi|^2\mathrm{d}v \leqslant  \int_{\{ |v_1| \geqslant \eta^{-\frac{1}{3}}\}} \eta |v_1||\psi|^2\mathrm{d}v \leqslant  \| \psi \|^2_{\mathcal{H}_\eta} . $$
$\bullet$ While on $\{ |v_1| \leqslant \eta^{-\frac{1}{3}}\}$, we introduce the function $\zeta_\eta$ defined by: $\zeta_{\eta}(v_1) := \zeta\big(\eta^{\frac{1}{3}} v_1)$, where $\zeta \in C^\infty(\RR)$  such that $0\leqslant \zeta \leqslant 1$, $\zeta \equiv 1$ on $B(0,1)$ and $\zeta \equiv 0$ outside of $B(0,2)$. Then, one has 
\begin{align*}
\eta^{\frac{2}{3}} \int_{\{ |v_1| \leqslant \eta^{-\frac{1}{3}}\}} |\psi|^2\mathrm{d}v &\leqslant \eta^{\frac{2}{3}} \int_{\{ |v_1| \leqslant 2\eta^{-\frac{1}{3}}\}} |\zeta_\eta \psi|^2 \mathrm{d}v \\
&= \eta^{\frac{2}{3}} \int_{\{ |v_1| \leqslant 2\eta^{-\frac{1}{3}}\}} \bigg|\int\int_{-2\eta^{-\frac{1}{3}}}^{v_1} \pa_{w_1}(\zeta_\eta \psi)\mathrm{d}w_1\bigg|^2\mathrm{d}v' \mathrm{d}v_1\\
&\leqslant \eta^{\frac{2}{3}} \int_{\{ |v_1| \leqslant 2\eta^{-\frac{1}{3}}\}}  \bigg(\int_{-2\eta^{-\frac{1}{3}}}^{v_1}\mathrm{d}w_1\bigg)\bigg(\int_{-2\eta^{-\frac{1}{3}}}^{v_1} |\pa_{w_1}(\zeta_\eta \psi)|^2\mathrm{d}w_1\bigg) \mathrm{d}v \\
&\leqslant 16 \big\|\pa_{v_1}(\zeta_\eta \psi)\big\|^2_{L^2(\{ |v_1| \leqslant 2\eta^{-\frac{1}{3}}\})} .
\end{align*} 
On the other hand, one has
\begin{align*}
 \big|\pa_{v_1}(\zeta_\eta \psi)\big|^2 &= |\zeta_\eta' \psi|^2+ |\zeta_\eta \pa_{v_1}\psi|^2 + \zeta_\eta\zeta_\eta'\big(\bar \psi \pa_{v_1}\psi + \psi \overline{ \pa_{v_1} \psi } \big) \\
 &\leqslant \big(\eta^{-1}|v_1|^{-1} |\zeta_\eta'|^2\big) \eta |v_1||\psi|^2+ |\pa_{v_1}\psi|^2 + 2\eta^{-\frac{1}{2}}|v_1|^{-\frac{1}{2}}|\zeta_\eta'| \big(\eta^{\frac{1}{2}}|v_1|^{\frac{1}{2}}|\psi|\big) |\pa_{v_1}\psi| ,
\end{align*}
and since $\zeta_\eta'=0$ except on $\{\eta^{-\frac{1}{3}}\leqslant |v_1|\leqslant 2\eta^{-\frac{1}{3}}\}$, where $|\eta^{-\frac{1}{2}}|v_1|^{-\frac{1}{2}}\zeta_\eta'(v_1)|\leqslant C$. Then, by integrating the last inequality and using Cauchy-Schwarz for the last term we get
\begin{align*}
\big\|\pa_{v_1}(\zeta_\eta \psi)\big\|^2_{L^2(\{ |v_1| \leqslant 2\eta^{-\frac{1}{3}}\})} &\lesssim \int_{\{ |v_1| \geqslant \eta^{-\frac{1}{3}}\}} \eta |v_1||\psi|^2 \mathrm{d}v + \|\pa_{v_1}\psi\|^2_{L^2(\{ |v_1| \leqslant 2\eta^{-\frac{1}{3}}\})} \lesssim \| \psi \|^2_{\mathcal{H}_\eta} .
\end{align*}
Note that we used the inclusion $\{\eta^{-\frac{1}{3}}\leqslant |v_1|\leqslant 2\eta^{-\frac{1}{3}}\} \subset \{ |v_1| \geqslant \eta^{-\frac{1}{3}} \}$. Hence, the inequality of Lemma \ref{Poincare1} holds.
\ep
\begin{lemma} The sesquilinear form $a$ is continuous on $\mathcal{H}_\eta\times\mathcal{H}_\eta$. Moreover, there exists a constant $C>0$, independent of $\lambda$ and $\eta$ such that, for all $\psi,\phi \in \mathcal{H}_\eta$
$$
|a(\psi,\phi)| \leqslant C \| \psi \|_{\mathcal{H}_\eta}\| \phi \|_{\mathcal{H}_\eta} .
$$
\end{lemma}
\bp It follows from the previous lemma that allows to handle the term $\ds \eta^{\frac{2}{3}}\lambda \int \psi \bar \phi \ \ud v$.
\ep
\begin{remark}\label{A representant de a}
By application of Riesz's theorem to continuous sesquilinear forms, there exists a continuous linear map $ A_{\lambda,\eta}\in \mathcal{L}({\mathcal {H}_\eta})$ such that $\ a(\psi,\phi) = \langle A_{\lambda,\eta}\psi, \phi\rangle_{\mathcal{H}_\eta}$ for all $\psi, \phi \in \mathcal{H}_\eta$.  \\
Note that $A_{\lambda,\eta}$ depends on $\lambda$ and $\eta$ since the form $a$ depends on these last parameters.
\end{remark}
\begin{lemma}\label{coercivite de a} Let $\eta>0$ and $\lambda\in\CC$ fixed, such that $|\lambda|\leqslant \lambda_0$ with $\lambda_0$ small enough. Let $A_{\lambda,\eta}$ be the linear operator representing the sesquilinear form $a$. Then, there exists a constant $C>0$, independent of $\lambda$ and $\eta$ such that
\begin{equation}\label{psi < A psi}
\| \psi \|_{\mathcal{H}_\eta} \leqslant C \| A_{\lambda,\eta} \psi \|_{\mathcal{H}_\eta} , \quad \forall \psi \in \mathcal{H}_\eta .
\end{equation}
\end{lemma}
\bp We have for all $a, b \in \RR$ and $z\in \CC$: $ |a+ \mathrm{i} b + z| \geqslant |a|-|z|$. Now, applying this inequality to $|a(\psi,\psi)|$ and using  Lemma \ref{Poincare1} for the term which contains $\lambda$, we write
\begin{align*}
|a(\psi,\psi)| &= \bigg| \int_{\RR^d} \big( |\nabla_v \psi|^2 + c_\gamma \frac{|\psi|^2}{\langle v \rangle^2} +  \mathrm{i} \eta v_1 |\psi|^2 -\lambda \eta^{\frac{2}{3}}|\psi|^2 \big) \mathrm{d}v \bigg|  \\
& \geqslant \bigg| \int_{\RR^d} \big( |\nabla_v \psi|^2 + c_\gamma \frac{|\psi|^2}{\langle v \rangle^2}\big) \mathrm{d}v \bigg| - |\lambda| \eta^{\frac{2}{3}}\| \psi \|_2^2 \\
&\geqslant \|\psi\|^2_{\mathcal{H}_0} - C_0|\lambda| \|\psi\|^2_{\mathcal{H}_\eta} .
\end{align*}
Then,  since $|a(\psi,\psi)|=|\langle A_{\lambda,\eta} \psi,\psi \rangle_{\mathcal{H}_\eta}|\leqslant \|A_{\lambda,\eta} \psi\|_{\mathcal{H}_\eta} \|\psi\|_{\mathcal{H}_\eta}$, we get
\begin{equation}\label{psi sur H_0 et nabla psi sur L^2}
\|\psi\|^2_{\mathcal{H}_0} = \|\nabla_v \psi\|^2_2 + c_\gamma \bigg\| \frac{\psi}{\langle v \rangle} \bigg\|^2_2 \leqslant \|A_{\lambda,\eta} \psi\|_{\mathcal{H}_\eta} \|\psi\|_{\mathcal{H}_\eta} + C_0 |\lambda| \|\psi\|^2_{\mathcal{H}_\eta} .
\end{equation}
Let denote 
$$ I^\eta_1 := \int_{\{ |v_1| \leqslant \eta^{-\frac{1}{3}}\}} \eta |v_1||\psi|^2\mathrm{d}v \quad \mbox{ and } \quad  I^\eta_2 := \int_{\{ |v_1| \geqslant \eta^{-\frac{1}{3}}\}} \eta |v_1||\psi|^2\mathrm{d}v . $$
Note that $\|\psi\|^2_{\mathcal{H}_\eta} = \|\psi\|^2_{\mathcal{H}_0} + I^\eta_1 + I^\eta_2$. To estimate $I^\eta_1$ and $I^\eta_2$, we need the following two steps.  \\

\noindent \textbf{Step 1: Estimation of $I^\eta_1$.} Let $\zeta_\eta$ be the function defined in the proof of Lemma \ref{Poincare1}. Then, 
$$ I^\eta_1 \leqslant  \int_{\{ |v_1| \leqslant 2\eta^{-\frac{1}{3}}\}} \eta |v_1||\zeta_\eta \psi|^2 \mathrm{d}v \leqslant \eta^{\frac{2}{3}} \int_{\{ |v_1| \leqslant 2\eta^{-\frac{1}{3}}\}} |\zeta_\eta \psi|^2 \mathrm{d}v  \leqslant 16 \big\|\pa_{v_1}(\zeta_\eta \psi)\big\|^2_{L^2(\{ |v_1| \leqslant 2\eta^{-\frac{1}{3}}\})} .$$
By the same calculations as in the proof of Lemma \ref{Poincare1} for $\big\|\pa_{v_1}(\zeta_\eta \psi)\big\|^2_{L^2(\{ |v_1| \leqslant 2 \eta^{-\frac{1}{3}}\})}$, we get
\begin{equation}\label{I^eta_1 < I^eta_2 + |grad psi|^2}
I^\eta_1 \leqslant C_1\bigg(I^\eta_2 + \|\nabla_v \psi\|^2_2\bigg) .
\end{equation}
\textbf{Step 2: Estimation of $I^\eta_2$.} Let $\chi_\eta$ this time be the function defined by $\chi_\eta(v_1):=\chi(\eta^{\frac{1}{3}}v_1)$ with $\chi \in C^\infty(\RR)$ such that: $-1 \leqslant \chi \leqslant 1$, $\chi \equiv -1$ on $]-\infty,-1]$,  $\chi \equiv 1$ on $[1,+\infty[$ and $\chi \equiv 0$ on $B(0,\frac{1}{2})$. Then,
$$ I^\eta_2 := \int_{\{ |v_1| \geqslant \eta^{-\frac{1}{3}}\}} \eta |v_1||\psi|^2\mathrm{d}v \leqslant  \int_{\{ |v_1| \geqslant \frac{1}{2}\eta^{-\frac{1}{3}}\}} \eta v_1 \chi_\eta \psi \bar \psi \mathrm{d}v .$$
By integrating the equation of $\psi$ multiplied by $\chi_\eta \bar \psi$ over $\{ |v_1| \geqslant \frac{1}{2}\eta^{-\frac{1}{3}}\}$ and taking the imaginary part, we obtain
$$ \int_{\{ |v_1| \geqslant \frac{1}{2}\eta^{-\frac{1}{3}}\}} \eta v_1 \chi_\eta \psi \bar \psi \mathrm{d}v =\mathrm{Im}\bigg( a(\psi,\chi_\eta \psi) - \int_{\{ |v_1| \geqslant \frac{1}{2}\eta^{-\frac{1}{3}}\}} \big[ \nabla_v \psi \cdot \nabla_v(\chi_\eta \bar \psi) - \lambda \eta^{\frac{2}{3}}  \chi_\eta \psi \bar \psi  \big] \mathrm{d}v \bigg) . $$  
For the first term, by Cauchy-Schwarz: $|\mathrm{Im}\ a(\psi,\chi_\eta \psi)|\leqslant \|A_{\lambda,\eta} \psi\|_{ \mathcal{H}_\eta}\|\chi_\eta\psi\|_{\mathcal{H}_\eta}$, and for the last term, by Lemma \ref{Poincare1}: $$
\bigg|\mathrm{Im}\lambda \eta^{\frac{2}{3}} \int_{\{ |v_1| \geqslant \frac{1}{2}\eta^{-\frac{1}{3}}\}} \chi_\eta \psi \bar \psi \mathrm{d}v \bigg| \leqslant C_0 |\lambda| \|\psi\|^2_{\mathcal{H}_\eta} .
$$ 
Finally, for the second term, we write 
\begin{align*}
\bigg| \mathrm{Im} \int_{\{ |v_1| \geqslant \frac{1}{2}\eta^{-\frac{1}{3}}\}} \nabla_v \psi \cdot \nabla_v(\chi_\eta \bar \psi) \mathrm{d}v \bigg| &= \bigg| \mathrm{Im} \int_{\{ |v_1| \geqslant \frac{1}{2}\eta^{-\frac{1}{3}}\}}  \chi_\eta'\bar\psi \pa_{v_1}\psi \mathrm{d}v \bigg| \\
&= \bigg| \mathrm{Im} \int_{\{ \frac{1}{2}\eta^{-\frac{1}{3}} \leqslant |v_1| \leqslant \eta^{-\frac{1}{3}}\}}  \chi_\eta'\bar\psi \pa_{v_1}\psi \mathrm{d}v\bigg| \\
&\leqslant 2C_2 \bigg|\int_{\{ \frac{1}{2}\eta^{-\frac{1}{3}} \leqslant |v_1| \leqslant \eta^{-\frac{1}{3}}\}}  \eta^{\frac{1}{2}}|v_1|^{\frac{1}{2}}|\psi||\pa_{v_1}\psi| \mathrm{d}v\bigg|  \\
&\leqslant 2C_2 \big(I^\eta_1\big)^{\frac{1}{2}} \big\| \nabla_v\psi\big\|_2 \hspace{2.35cm} (\mbox{by Cauchy-Schwarz})\\
&\leqslant C_3\bigg(I^\eta_2 + \|\nabla_v \psi\|^2_2\bigg)^{\frac{1}{2}} \big\| \nabla_v\psi\big\|_2 , \quad (\mbox{by  inequality } \eqref{I^eta_1 < I^eta_2 + |grad psi|^2}) \\
&\leqslant \frac{1}{4} I^\eta_2 + C\|\nabla_v \psi\|^2_2 , 
\end{align*}
where we used the inequality: $ab \leqslant C_3a^2+\frac{b^2}{4C_3} $ in the last line and where $C_2=\underset{\frac{1}{2 }\leqslant |t| \leqslant 1}{\sup}|t^{-\frac{1}{2}} \chi'(t)|=\|\eta^{-\frac{1}{2}}|v_1|^{-\frac{1}{2}}\chi_\eta'\|_{L^\infty(\{ \frac{1}{2}\eta^{-\frac{1}{3}} \leqslant |v_1| \leqslant \eta^{-\frac{1}{3}}\})}$, $C_3=2\sqrt{C_1}C_2$ and $C= C_3+\frac{1}{4}$. Therefore,
$$I^\eta_2 \leqslant \int_{\{ |v_1| \geqslant \frac{1}{2}\eta^{-\frac{1}{3}}\}} \eta v_1 \chi_\eta \psi \bar \psi \mathrm{d}v \leqslant \|A_{\lambda,\eta} \psi\|_{ \mathcal{H}_\eta}\|\chi_\eta\psi\|_{\mathcal{H}_\eta} + \frac{1}{4} I^\eta_2 + C \big\|\nabla_v \psi \big\|^2_2 +C_0 |\lambda| \|\psi\|^2_{\mathcal{H}_\eta} . $$
Thus,
$$ I^\eta_2 \leqslant C \bigg( \|A_{\lambda,\eta} \psi\|_{ \mathcal{H}_\eta} \|\chi_\eta\psi\|_{\mathcal{H}_\eta} +  \big\|\nabla_v \psi \big\|^2_2 + |\lambda| \|\psi\|^2_{\mathcal{H}_\eta} \bigg).$$
Recall that we have $\|\nabla_v \psi \|^2_2\leqslant \|A_{\lambda,\eta} \psi\|_{ \mathcal{H}_\eta} \|\psi\|_{\mathcal{H}_\eta}$ thanks to \eqref{psi sur H_0 et nabla psi sur L^2}. Hence,
\begin{equation}\label{I^eta_2 (1)}
I^\eta_2 \leqslant C \bigg( \|A_{\lambda,\eta} \psi\|_{ \mathcal{H}_\eta} \|\chi_\eta\psi\|_{\mathcal{H}_\eta} + |\lambda| \|\psi\|^2_{\mathcal{H}_\eta} \bigg).
\end{equation}
It only remains to handle the term $\|\chi_\eta\psi\|_{\mathcal{H}_\eta}$. We have, as in the proof of Lemma \ref{Poincare1},
\begin{align*}
 \big|\nabla_v(\chi_\eta \psi)\big|^2 &= |\chi_\eta' \psi|^2+ \big|\chi_\eta \nabla_v\psi \big|^2 + \chi_\eta\chi_\eta'\big(\bar \psi \pa_{v_1}\psi + \psi \overline{ \pa_{v_1} \psi } \big) \\
 &\leqslant \big(\eta^{-1}|v_1|^{-1} |\chi_\eta'|^2\big) \eta |v_1||\psi|^2+ \big|\nabla_v\psi \big|^2 + 2\eta^{-\frac{1}{2}}|v_1|^{-\frac{1}{2}}|\chi_\eta'| \big(\eta^{\frac{1}{2}}|v_1|^{\frac{1}{2}}|\psi|\big) |\pa_{v_1}\psi| .
\end{align*}
Then, $ \|\nabla_v(\chi_\eta \psi)\|^2_2 \leqslant C\big(I^\eta_1+\|\nabla_v \psi\|^2_2\big) \leqslant C \|\psi\|^2_{\mathcal{H}_\eta}$ and therefore, $\|\chi_\eta \psi\|^2_{\mathcal{H}_\eta} \leqslant C \|\psi\|^2_{\mathcal{H}_\eta}$. 
By injecting this last inequality into \eqref{I^eta_2 (1)}, we get 
\begin{equation}\label{I^eta_2}
I^\eta_2 \leqslant C \bigg( \|A_{\lambda,\eta} \psi\|_{ \mathcal{H}_\eta} \|\psi\|_{\mathcal{H}_\eta} + |\lambda| \|\psi\|^2_{\mathcal{H}_\eta} \bigg) .
\end{equation}
Thus, by summing \eqref{psi sur H_0 et nabla psi sur L^2}, \eqref{I^eta_1 < I^eta_2 + |grad psi|^2} and \eqref{I^eta_2}, we obtain
$$ \|\psi\|_{\mathcal{H}_\eta}^2 \leqslant C \bigg( \|A_{\lambda,\eta} \psi\|_{ \mathcal{H}_\eta} \|\psi\|_{\mathcal{H}_\eta} + |\lambda| \|\psi\|^2_{\mathcal{H}_\eta} \bigg) .$$
Finally, we obtain  inequality \eqref{psi < A psi} by the inequality $ab \leqslant C a^2+\frac{b^2}{4C}$ applied to the term $\|A_{\lambda,\eta} \psi\|_{ \mathcal{H}_\eta} \|\psi\|_{\mathcal{H}_\eta}$, and with $\lambda$ small enough: $|\lambda|\leqslant \frac{1}{4C}$. 
\ep

\begin{lemma}[Complementary Lemma]\label{lemme compl.}
Let $\eta>0$ fixed and let $\lambda_0>0$ small enough. Let $\lambda \in \CC$ such that $|\lambda|\leqslant \lambda_0$. Then, for all $\psi, F \in \mathcal{H}_\eta$ such that $|a(\psi,\psi)| \leqslant C \| F \|_{\mathcal{H}_\eta} \| \psi \|_{\mathcal{H}_\eta}$, the following inequality holds
\begin{equation}\label{psi < C F dans H_eta}
\| \psi \|_{\mathcal{H}_\eta} \leqslant \tilde C \| F \|_{\mathcal{H}_\eta} ,
\end{equation}
where $C$ and $\tilde C$ are two positive constants that do not depend on $\lambda$ and $\eta$.
\end{lemma}
\bp 
The proof is identical to that of the previous Lemma,  just replace the inequality $|a(\psi,\psi)|=|\langle A_{\lambda,\eta} \psi,\psi \rangle_{\mathcal {H}_\eta}|\leqslant \|A_{\lambda,\eta} \psi\|_{\mathcal{H}_\eta} \|\psi \|_{\mathcal{H}_\eta}$ by $|a(\psi,\psi)| \leqslant C \| F \|_{\mathcal{H}_\eta} \| \psi \|_{\mathcal{H}_\eta}$.
\ep

\noindent Let denote by $\mathcal{H}_\eta'$ the topological dual of $\mathcal{H}_\eta$. By the Riesz representation theorem, for all $F \in \mathcal{H}_\eta'$, there exists a unique $f \in \mathcal{H}_\eta$ such that 
$$ (F,\phi) = \langle f,\phi\rangle_{\mathcal{H}_\eta} , \quad \forall \phi \in \mathcal{H}_\eta , $$ 
where $(F,\phi)$ denotes the value taken by $F \in \mathcal{H}_\eta'$ in $\phi \in \mathcal{H}_\eta$. Then, by Remark \ref{A representant de a}, the problem  
\begin{equation}\label{formulation fable}
a(\psi,\phi) = (F,\phi) , \quad \forall \phi \in \mathcal{H}_\eta
\end{equation}
is equivalent to the problem $A_{\lambda,\eta} \psi = f$, $f \in \mathcal{H}_\eta$. Therefore, equivalent to the invertibility of $A_{\lambda,\eta}$.
\begin{proposition}[Existence of solution to the the variational problem]\label{inversibilite de A} Let $\eta_0>0$ and $\lambda_0>0$ small enough. Let $\eta \in [0,\eta_0]$ and $\lambda \in \CC$ fixed, with $|\lambda|\leqslant \lambda_0$. 
For all $F \in \mathcal{H}_\eta'$,  equation \eqref{formulation fable} admits a unique solution $\psi^{\lambda,\eta} \in \mathcal{H}_\eta\subset  \mathcal{H}_0$, satisfying the following estimate
\begin{equation}\label{psi_H-eta < C F_H-eta'}
\| \psi^{\lambda,\eta} \|_{\mathcal{H}_0} \leqslant \| \psi^{\lambda,\eta} \|_{\mathcal{H}_\eta} \leqslant C \| F \|_{\mathcal{H}_\eta'} ,
\end{equation}
where $C$ is a positive constant that does not depend on $\lambda$ and $\eta$. Moreover, for $F \in L^2_{\langle v\rangle^2} \subset \mathcal{H}_\eta'$ we have
\begin{equation}\label{psi_H-eta < C <v>F_L^2}
\| \psi^{\lambda,\eta} \|_{\mathcal{H}_0} \leqslant\| \psi^{\lambda,\eta} \|_{\mathcal{H}_\eta} \leqslant C \| F \|_{L^2_{\langle v\rangle^2}} ,
\end{equation}
where $L^2_{\langle v\rangle^2}$ denote the weighted $L^2$space: $\ds L^2_{\langle v\rangle^2}:=\left\{f: \RR^d \longrightarrow \CC ; \ \int_{\RR^d} |f|^2 \langle v\rangle^2 \mathrm{d}v < \infty \right\}$.
\end{proposition}
\begin{remark}
The sesquilinear form $a$ depends continuously on $\eta$ and holomorphically on $\lambda$. \\
The solution in the previous Proposition, is for $\lambda$ and $\eta$ fixed, and it depends on $\lambda$ and $\eta$ since $a$ depends on these last parameters.
\end{remark}
\bpp \ref{inversibilite de A}. This proof was taken from \cite{VKK} to prove the first statement of the Lax-Milgram lemma  [page 235]. We want to prove that the linear map $A_{\lambda,\eta}$ representing the sesquilinear form $a$ is invertible with continuous inverse, since it implies that for all $f \in \mathcal{H}_\eta$, the equation $A_{\lambda,\eta}\psi=f$ admits a unique solution $\psi^{\lambda,\eta} \in \mathcal{H}_\eta$. 

First,  inequality \eqref{psi < A psi} of Lemma \ref{coercivite de a}, $\|\psi\|_{\mathcal{H}_\eta}\leqslant C\|A_{\lambda,\eta}\psi\|_{\mathcal{H}_\eta}$, shows that $A_{\lambda,\eta}$ is injective with continuous inverse, so it is a topological isomorphism from $\mathcal{H}_\eta$ to $\mathrm{R}(A_{\lambda,\eta})$; in particular $\mathrm{R}(A_{\lambda,\eta})$ is complete and therefore closed in $\mathcal{H}_\eta$, where we denote by $\mathrm{R}(A_{\lambda,\eta})$ the range of the operator $A_{\lambda,\eta}$, i.e., $\mathrm{R}(A_{\lambda,\eta}):=\{ f \in \mathcal{H}_\eta ; \ f = A_{\lambda,\eta}\psi ,\ \psi \in \mathcal{H}_\eta \}$. To show that $A_{\lambda,\eta}$ is surjective, it is enough to prove that $\mathrm{R}(A_{\lambda,\eta})$ is dense; for this, let $\phi_0 \in \mathcal{H}_\eta$ such that $\langle A_{\lambda,\eta} \psi,\phi_0\rangle_{\mathcal{H}_\eta}=0$ for all $\psi \in \mathcal{H}_\eta$; taking $\psi = \phi_0$ we get $a(\phi_0,\phi_0)=0$, which gives $\phi_0 =0$. \\
Inequality \eqref{psi_H-eta < C F_H-eta'} comes from
$$
\| \psi^{\lambda,\eta} \|_{\mathcal{H}_\eta} \leqslant C\| A_{\lambda,\eta}\psi^{\lambda,\eta}\|_{\mathcal{H}_\eta}\leqslant \| f \|_{\mathcal{H}_\eta}\leqslant\| F \|_{\mathcal{H}_\eta'}.
$$
 For the second one, it comes from the fact that the weighted space $L^2_{\langle v\rangle^2}$ is continuously embedded in $\mathcal{H}_\eta'$.
\epp

\noindent We will denote by $T_{\lambda,\eta}$ the inverse operator of $L_{\lambda,\eta}$ for $\lambda$ and $\eta$ fixed, i.e., the operator which associates to $F$ the solution $\psi^{\lambda,\eta} =: T_{\lambda,\eta}(F)$.
\subsection{Implicit function theorem}
In this subsection, we use the operator $T_{\lambda,\eta}$ to rewrite equation \eqref{eq penalisee2} as a fixed point problem for the identity plus a compact map. Then, the Fredholm Alternative will allow us to apply the implicit function theorem in order to have the existence of solutions. For this purpose, let's define $F: \{\lambda\in\CC;|\lambda|\leqslant\lambda_0\}\times[0,\eta_0]\times \mathcal{H}_0\longrightarrow \mathcal{H}_0$ by 
$$F(\lambda,\eta,h):= h - \mathcal{T}_{\lambda,\eta}(h) ,$$ 
with
$$\mathcal{T}_{\lambda,\eta}(h):= T_{\lambda,\eta}\big[Vh-\langle h-M,\Phi \rangle \Phi \big] .$$
Note that finding a solution $h(\lambda,\eta)$ solution to $F\big(\lambda,\eta,h(\lambda,\eta)\big)=0$ gives a solution to the penalized equation by taking $M_{\lambda,\eta}=h(\lambda,\eta)$.\\

\ni The function $\Phi$ satisfies the following assumptions:
\begin{enumerate}
\item For all $v$ in $\RR^d$, $\Phi(v)=\Phi(-v)>0$.
\item The function $ \Phi$ belongs to the weighted Sobolev space $H^1_{\langle v\rangle^2}:=H^1(\RR^d,\langle v\rangle^2\mathrm{d}v)$,  and for all $v$ in $\RR^d$, $\Phi(v) \leqslant \frac{M(v)}{\langle v \rangle^2}$.  
\item Even if it means multiplying $\Phi$ by a constant, we can take it such that $\langle \Phi, M \rangle = 1$.
\end{enumerate}

\noindent For the following, we will take the function $\Phi:=c_{\gamma,d} \langle v\rangle^{-2-\gamma} $ which satisfies all the previous assumptions, where $c_{\gamma,d }=\big(\int_{ \RR^d}\langle v\rangle^{-2-2\gamma}\mathrm{d}v\big)^{-1}$.
\begin{remark}
Note that the operator $\mathcal{T}_{\lambda,0}$ does not depend on $\lambda$ since $T_{\lambda,0}$ does not. Let's denote it by $\mathcal{T}_0$. Also, $\mathcal{T}_{\lambda,\eta}$ is affine with respect to $h$, we denote by $\mathcal{T}^l_{\lambda,\eta}$ its linear part.
\end{remark}
\begin{lemma}[Continuity of $\mathcal{T}_{\lambda,\eta}$]\label{continuite de T_eta}
Let $\eta_0>0$ and $\lambda_0>0$ small enough. Let $\eta \in [0,\eta_0]$ and $\lambda \in \mathbb{C}$ such that $|\lambda|\leqslant\lambda_0$. Then, 
\begin{enumerate}
\item The map $\mathcal{T}_{\lambda,\eta} : \mathcal{H}_0 \longrightarrow \mathcal{H}_\eta$ is continuous.  Moreover, there exists a constant $C>0$, independent of $\lambda$ and $\eta$ such that
\begin{equation}\label{continuite de H_0 dans H_eta}
\| \mathcal{T}_{\lambda,\eta}^l(h)\|_{\mathcal{H}_\eta} \leqslant C \| h \|_{\mathcal{H}_0} , \quad \forall h \in \mathcal{H}_0 ,
\end{equation}
and the embedding $ \mathcal{T}_{\lambda,\eta}^l(\mathcal{H}_0) \subset \mathcal{H}_\eta \subset \mathcal{H}_0 $ holds for all $ \eta \in [0,\eta_0]$ and for all $\lambda \in \{|\lambda|\leqslant \lambda_0\}$.  Hence the map $\mathcal{T}_{\lambda,\eta} : \mathcal{H}_0 \longrightarrow \mathcal{H}_0$ is continuous. 
\item The map $\mathcal{T}_{\lambda,\eta}$ is continuous with respect to $\lambda$ and $\eta$.  Moreover, there exists a constant $C>0$, independent of $\lambda$ and $\eta$ such that, for all $\eta' \in [0,\eta_0]$ and for all $|\lambda'| \leqslant \lambda_0$
\begin{equation}\label{T_eta(h)-T_0(h) -->0}
\|\mathcal{T}_{\lambda,\eta}(h)-\mathcal{T}_{\lambda,\eta'}(h)\|_{\mathcal{H}_0} \leqslant C \bigg(\bigg|1-\frac{\eta'}{\eta} \bigg|+\bigg|1-\bigg|\frac{\eta'}{\eta}\bigg|^{\frac{2}{3}} \bigg| \bigg) \big( \|h \|_{\mathcal{H}_0} +\|\Phi\|_{L^2_{\langle v \rangle^2}} \big)
\end{equation}
and
\begin{equation}
\|\mathcal{T}_{\lambda,\eta}(h)-\mathcal{T}_{\lambda',\eta}(h)\|_{\mathcal{H}_0} \leqslant C  |\lambda-\lambda'| \big( \|h \|_{\mathcal{H}_0} +\|\Phi\|_{L^2_{\langle v \rangle^2}} \big), 
\end{equation}
for all $h \in \mathcal{H}_0$.
\end{enumerate}
\end{lemma}

\noindent \bp \textbf{1. }The first point follows from the second inequality of  Proposition \ref{psi_H-eta < C <v>F_L^2}.  Indeed,  we have by \eqref{psi_H-eta < C <v>F_L^2}, for all $\ F \in L^2_{\langle v\rangle^2}$
$$  \|T_{\lambda,\eta}(F) \|_{\mathcal{H}_\eta} \leqslant C \| F \|_{L^2_{\langle v\rangle^2}}  $$ 
For $h_1, h_2 \in \mathcal{H}_0$, we have $\mathcal{T}_{\lambda,\eta}(h_1)-\mathcal{T}_{\lambda,\eta}(h_2)=\mathcal{T}_{\lambda,\eta}^l(h_1-h_2)$. Let denote $h:=h_1-h_2$ and $F := V h - \langle h,\Phi\rangle \Phi \in L^2_{\langle v \rangle^2}$. We have $\mathcal{T}_{\lambda,\eta}^l(h)=T_{\lambda,\eta}(F) $. Thus, by the last inequality and by Cauchy-Schwarz for the term $|\langle h,\Phi\rangle|$, we obtain
$$ \|\mathcal{T}_{\lambda,\eta}^l(h) \|_{\mathcal{H}_\eta} \leqslant C \bigg\|\langle v\rangle^2 V \frac{h}{\langle v\rangle} -\langle h,\Phi\rangle \langle v\rangle \Phi\bigg\|_2 \leqslant C\bigg(\big\|\langle v\rangle^2 V \big\|_\infty+ \big\|\langle v\rangle \Phi \big\|_2^2\bigg) \bigg\|\frac{h}{\langle v\rangle }\bigg\|_{L^2}  \leqslant \tilde C \|h\|_{\mathcal{H}_0}. $$
The embedding $\ \mathcal{T}_{\lambda,\eta}^l(\mathcal{H}_0) \subset \mathcal{H}_\eta \subset \mathcal{H}_0 \ $ comes from the previous inequality and the fact that $\|\mathcal{T}_{\lambda,\eta}^l(h)\|_{\mathcal{H}_0} \leqslant  \|\mathcal{T}_{\lambda,\eta}^l(h)\|_{\mathcal{H}_\eta}$ for all $h \in \mathcal{H}_0$. \\

\noindent \textbf{2. } Let $\eta_0>0$ and $\lambda_0>0$ small enough.  Let $\eta \in [0,\eta_0]$ and $\lambda\in \CC$ such that $|\lambda| \leqslant \lambda_0$.  Recall that $T_{\lambda,\eta}$ is the inverse of $L_{\lambda,\eta}:= \tilde Q +  \mathrm{i} \eta v_1 - \lambda \eta^{\frac{2}{3}}$ with $\ \tilde Q := -\Delta_v + \tilde W(v)$.  \\

\noindent \textbf{Continuity of $\mathcal{T}_{\lambda,\eta}$ with respect to $\lambda$.} Let $\lambda' \in \CC$ such that $|\lambda'|\leqslant \lambda_0$.  We have  for $h \in \mathcal{H}_0$
$$ \big[\tilde Q +  \mathrm{i} \eta v_1 - \lambda\eta^{\frac{2}{3}}\big]\big(T_{\lambda,\eta}\big[Vh-\langle h-M,\Phi \rangle \Phi \big]\big) =  Vh-\langle h-M,\Phi \rangle \Phi$$
and 
$$ \big[\tilde Q +  \mathrm{i} \eta v_1 - \lambda'\eta^{\frac{2}{3}}\big]\big(T_{\lambda',\eta}\big[Vh-\langle h-M,\Phi \rangle \Phi \big]\big) =  Vh-\langle h-M,\Phi \rangle \Phi .$$
Thus,  the function $\ \mathcal{T}_{\lambda,\eta}(h)-\mathcal{T}_{\lambda',\eta}(h) = (T_{\lambda,\eta}-T_{\lambda',\eta})\big[Vh-\langle h-M,\Phi \rangle \Phi \big]\ $ satisfies the equation
$$ \tilde Q [\mathcal{T}_{\lambda,\eta}(h)-\mathcal{T}_{\lambda',\eta}(h)] +  \mathrm{i} \eta v_1 [\mathcal{T}_{\lambda,\eta}(h)-\mathcal{T}_{\lambda',\eta}(h)] - \lambda\eta^{\frac{2}{3}} [\mathcal{T}_{\lambda,\eta}(h)-\mathcal{T}_{\lambda',\eta}(h)] = (\lambda-\lambda')\eta^{\frac{2}{3}} \mathcal{T}_{\lambda',\eta}(h) . $$
Then, by integrating the previous equality multiplied by $\overline{[\mathcal{T}_{\lambda,\eta}(h)- \mathcal{T}_{\lambda',\eta}(h)]}$, we obtain
$$ a_{\lambda,\eta}\big(\mathcal{T}_{\lambda,\eta}(h)-\mathcal{T}_{\lambda',\eta}(h),\mathcal{T}_{\lambda,\eta}(h)-\mathcal{T}_{\lambda',\eta}(h)\big) = (\lambda-\lambda')\eta^{\frac{2}{3}} \int_{\RR^d}\mathcal{T}_{\lambda',\eta}(h) \overline{[\mathcal{T}_{\lambda,\eta}(h)-\mathcal{T}_{\lambda',\eta}(h)]} \mathrm{d}v . $$
Now, by Cauchy-Schwarz inequality
$$ \bigg| (\lambda-\lambda')\eta^{\frac{2}{3}} \int_{\RR^d}\mathcal{T}_{\lambda',\eta}(h) \overline{[\mathcal{T}_{\lambda,\eta}(h)-\mathcal{T}_{\lambda',\eta}(h)]} \mathrm{d}v \bigg| \leqslant |\lambda-\lambda'| \eta^{\frac{2}{3}}\|\mathcal{T}_{\lambda',\eta}(h) \|_2 \|\mathcal{T}_{\lambda,\eta}(h)-\mathcal{T}_{\lambda',\eta}(h)\|_2 ,$$
and by Lemma \ref{Poincare1} we get
$$ \bigg| (\lambda-\lambda')\eta^{\frac{2}{3}} \int_{\RR^d}\mathcal{T}_{\lambda',\eta}(h) \overline{[\mathcal{T}_{\lambda,\eta}(h)-\mathcal{T}_{\lambda',\eta}(h)]} \mathrm{d}v \bigg| \leqslant C |\lambda-\lambda'| \|\mathcal{T}_{\lambda',\eta}(h) \|_{\mathcal{H}_\eta}  \|\mathcal{T}_{\lambda,\eta}(h)-\mathcal{T}_{\lambda',\eta}(h)\|_{\mathcal{H}_\eta} . $$
Therefore, 
$$ \big| a_{\lambda,\eta} \big(\mathcal{T}_{\lambda,\eta}(h)-\mathcal{T}_{\lambda',\eta}(h) ,\mathcal{T}_{\lambda,\eta}(h)-\mathcal{T}_{\lambda',\eta}(h)\big)\big| \leqslant C  |\lambda-\lambda'|  \| \mathcal{T}_{\lambda',\eta}(h) \|_{\mathcal{H}_\eta} \| \mathcal{T}_{\lambda,\eta}(h)-\mathcal{T}_{\lambda',\eta}(h) \|_{\mathcal{H}_\eta} .$$
Hence, by the Complementary Lemma \ref{lemme compl.}, we write
$$ \| \mathcal{T}_{\lambda,\eta}(h)-\mathcal{T}_{\lambda',\eta}(h) \|_{\mathcal{H}_0} \leqslant \| \mathcal{T}_{\lambda,\eta}(h)-\mathcal{T}_{\lambda',\eta}(h) \|_{\mathcal{H}_\eta} \leqslant C  |\lambda-\lambda'|  \| \mathcal{T}_{\lambda',\eta}(h) \|_{\mathcal{H}_\eta} .$$
That leads to 
$$ \| \mathcal{T}_{\lambda,\eta}(h)-\mathcal{T}_{\lambda',\eta}(h) \|_{\mathcal{L}(\mathcal{H}_0)} \leqslant  C  |\lambda-\lambda'| \big( \|h \|_{\mathcal{H}_0} +\|\langle v \rangle\Phi\|_2 \big) . $$

\noindent \textbf{Continuity of $\mathcal{T}_{\lambda,\eta}$ with respect to $\eta$.} Let $\eta' \in [0,\eta_0]$.  Without loss of generality, we can assume that $\eta \leqslant \eta'$. Then, as before, we have  for $h \in \mathcal{H}_0$
$$ \big[\tilde Q +  \mathrm{i} \eta v_1 - \lambda\eta^{\frac{2}{3}}\big]\big(T_{\lambda,\eta}\big[Vh-\langle h-M,\Phi \rangle \Phi \big]\big) =  Vh-\langle h-M,\Phi \rangle \Phi$$
and 
$$ \big[\tilde Q +  \mathrm{i} \eta' v_1 - \lambda{\eta'}^{\frac{2}{3}}\big]\big(T_{\lambda,\eta'}\big[Vh-\langle h-M,\Phi \rangle \Phi \big]\big) =  Vh-\langle h-M,\Phi \rangle \Phi .$$
Thus,  the function $\ \mathcal{T}_{\lambda,\eta}(h)-\mathcal{T}_{\lambda,\eta'}(h) = (T_{\lambda,\eta}-T_{\lambda,\eta'})\big[Vh-\langle h-M,\Phi \rangle \Phi \big] \ $ satisfies the equation 
$$ [\tilde Q +  \mathrm{i} \eta v_1 - \lambda\eta^{\frac{2}{3}}] (\mathcal{T}_{\lambda,\eta}(h)-\mathcal{T}_{\lambda,\eta'}(h)) = [ \mathrm{i} (\eta-\eta')v_1 - \lambda(\eta^{\frac{2}{3}}-{\eta'}^{\frac{2}{3}})] \mathcal{T}_{\lambda,\eta'}(h) ,$$
and integrating this equation against $\overline{[\mathcal{T}_{\lambda,\eta}(h) - \mathcal{T}_{\lambda,\eta'}(h)]}$ we get  
\begin{align*}
a\big(\mathcal{T}_{\lambda,\eta}(h)-\mathcal{T}_{\lambda,\eta'}(h),\mathcal{T}_{\lambda,\eta}(h)-\mathcal{T}_{\lambda,\eta'}(h)\big) &=  \mathrm{i} (\eta-\eta') \int_{\RR^d}  v_1 \mathcal{T}_{\lambda,\eta'}(h) \overline{[\mathcal{T}_{\lambda,\eta}(h)-\mathcal{T}_{\lambda,\eta'}(h)]} \mathrm{d}v \\
& - \lambda(\eta^{\frac{2}{3}}-{\eta'}^{\frac{2}{3}}) \int_{\RR^d}   \mathcal{T}_{\lambda,\eta'}(h) \overline{[\mathcal{T}_{\lambda,\eta}(h)-\mathcal{T}_{\lambda,\eta'}(h)]} \mathrm{d}v \\
&=: I^{\lambda,\eta,\eta'}_1 + I^{\lambda,\eta,\eta'}_2 .
\end{align*}
For $I^{\lambda,\eta,\eta'}_1$, we write
\begin{align*}
 \big| I^{\lambda,\eta,\eta'}_1 \big| &\leqslant \bigg|1-\frac{\eta'}{\eta} \bigg| \big\|\eta^\frac{1}{2} |v_1|^{\frac{1}{2}} \mathcal{T}_{\lambda,\eta'}(h) \big\|_2 \big\|\eta^\frac{1}{2} |v_1|^{\frac{1}{2}} [\mathcal{T}_{\lambda,\eta}(h)-\mathcal{T}_{\lambda,\eta'}(h)] \big\|_2  \\
 &\leqslant \bigg|1-\frac{\eta'}{\eta} \bigg| \| \mathcal{T}_{\lambda,\eta'}(h) \|_{\mathcal{H}_\eta} \| \mathcal{T}_{\lambda,\eta}(h)-\mathcal{T}_{\lambda,\eta'}(h) \|_{\mathcal{H}_\eta} .
\end{align*}
Now for $I^{\lambda,\eta,\eta'}_2$, by using Lemma \ref{Poincare1}, we write
\begin{align*}
\big| I^{\lambda,\eta,\eta'}_2 \big| &\leqslant \eta^{\frac{2}{3}} |\lambda|\bigg|1-\bigg|\frac{\eta'}{\eta}\bigg|^{\frac{2}{3}}\bigg|   \|\mathcal{T}_{\lambda,\eta'}(h) \|_2 \| \mathcal{T}_{\lambda,\eta}(h)-\mathcal{T}_{\lambda,\eta'}(h) \|_2 \\
& \leqslant C |\lambda|\bigg|1-\bigg|\frac{\eta'}{\eta}\bigg|^{\frac{2}{3}}\bigg|   \| \mathcal{T}_{\lambda,\eta'}(h) \|_{\mathcal{H}_\eta} \| \mathcal{T}_{\lambda,\eta}(h)-\mathcal{T}_{\lambda,\eta'}(h) \|_{\mathcal{H}_\eta} .
\end{align*}
Hence,
\begin{align*}
&\big|a_{\lambda,\eta}\big(\mathcal{T}_{\lambda,\eta}(h)-\mathcal{T}_{\lambda,\eta'}(h),\mathcal{T}_{\lambda,\eta}(h)-\mathcal{T}_{\lambda,\eta'}(h)\big)\big| \\
&\leqslant \left(\bigg|1-\frac{\eta'}{\eta} \bigg|+C|\lambda|\bigg|1-\bigg|\frac{\eta'}{\eta}\bigg|^{\frac{2}{3}}\bigg| \right)  \| \mathcal{T}_{\lambda,\eta'}(h) \|_{\mathcal{H}_\eta} \| \mathcal{T}_{\lambda,\eta}(h)-\mathcal{T}_{\lambda,\eta'}(h) \|_{\mathcal{H}_\eta} .
\end{align*}
Which implies, by inequality \eqref{psi < C F dans H_eta} of the complementary lemma, that
$$\| \mathcal{T}_{\lambda,\eta}(h)-\mathcal{T}_{\lambda,\eta'}(h) \|_{\mathcal{H}_\eta} \leqslant \left(\bigg|1-\frac{\eta'}{\eta} \bigg|+C|\lambda|\bigg|1-\bigg|\frac{\eta'}{\eta}\bigg|^{\frac{2}{3}}\bigg| \right)  \| \mathcal{T}_{\lambda,\eta'}(h) \|_{\mathcal{H}_\eta} . $$
Then, since $\| \mathcal{T}_{\lambda,\eta}(h)-\mathcal{T}_{\lambda,\eta'}(h) \|_{\mathcal{H}_0} \leqslant \| \mathcal{T}_{\lambda,\eta}(h)-\mathcal{T}_{\lambda,\eta'}(h) \|_{\mathcal{H}_\eta} $ and since $\eta \leqslant \eta'$ implies that $\| \mathcal{T}_{\lambda,\eta'}(h) \|_{\mathcal{H}_\eta} \leqslant \| \mathcal{T}_{\lambda,\eta'}(h) \|_{\mathcal{H}_{\eta'}} \leqslant C \big(\| h \|_{\mathcal{H}_0} + \| \langle v \rangle \Phi \|_2 \big)$, we get
$$ \| \mathcal{T}_{\lambda,\eta}(h)-\mathcal{T}_{\lambda,\eta'}(h) \|_{\mathcal{H}_0} \leqslant  C  \left(\bigg|1-\frac{\eta'}{\eta} \bigg|+C\lambda_0\bigg|1-\bigg|\frac{\eta'}{\eta}\bigg|^{\frac{2}{3}}\bigg| \right) \big(\| h \|_{\mathcal{H}_0} + \| \langle v \rangle \Phi \|_2 \big) . $$
Which ends of the proof.
\ep 
\begin{lemma}\label{compacite de T_0^l}
The map $\mathcal{T}^l_0$ is compact.
\end{lemma}

\noindent \bp  First, since the two functions $g_1:=\langle v\rangle^2V$ and $g_2:= \Phi$ belong to $C_0^1(\RR^d,\RR)$ and $H^1_{\langle v\rangle^2}(\RR^d,\RR)$ respectively, where $C_0^1(\RR^d,\RR)$ denote the space of $C^1$ functions converging to $0$ at infinity as well as their first derivatives, then for $\eps>0 $, there exists $g_1^\varepsilon, g_2^\varepsilon \in C_c^\infty(\RR^d,\RR)$ such that $\|g^\eps_1-g_1\|_{W^{1,\infty}} \leqslant \frac{\eps}{2C} $ and $\| g^\eps_2-g_2\|_{H^1_{\langle v \rangle^2}} \leqslant \frac{\eps}{2C} $, where $C$ is the constant of inequality \eqref{continuite de H_0 dans H_eta}. Now if we denote by $\mathcal{T}_0^{\eps}$ the operator $\mathcal{T}_0^{\eps}(h):= T_0\big[g^\eps_1\frac{h}{\langle v\rangle^2}-\langle h,\Phi\rangle g^\eps_2\big]$, then we can write:
\begin{align*}
\|\mathcal{T}_0^l(h)-\mathcal{T}_0^{\eps}(h)\|_{\mathcal{H}_0} &= \big\|T_0\big[(g^\eps_1-g_1)h/\langle v\rangle^2-\langle h,\Phi\rangle (g^\eps_2-g_2)\big]\big\|_{\mathcal{H}_0} \\
&\leqslant C\bigg(\|g^\eps_1-g_1\|_\infty + \big\|\langle v\rangle \Phi \big\|_2 \|g^\eps_2-g_2\|_{L^2_{\langle v\rangle^2}} \bigg) \|h\|_{\mathcal{H}_0} \\
&\leqslant \eps \|h\|_{\mathcal{H}_0} .
\end{align*}
Hence, $ \|\mathcal{T}_0^l-\mathcal{T}_0^{\eps}\|_{\mathcal{L}(\mathcal{H}_0)}\leqslant \eps $. Thus, the operator $\mathcal{T}_0^l$ can be seen as the limit of the operator $\mathcal{T}_0^{\eps}$ when $\eps$ goes to $0$. Indeed for $(h_n)_n \subset \mathcal{H}_0$ such that $\|h_n \|_{\mathcal{H}_0} \leqslant 1$ we have up to a subsequence, $h_n\rightharpoonup h$ in $\mathcal{H}_0$. Moreover, we have
\begin{align}\label{T^l(h_n)-T^l(h)}
\|\mathcal{T}^l_0(h_n)-\mathcal{T}^l_0(h)\|_{\mathcal{H}_0} &\leqslant \|\mathcal{T}^l_0(h_n)-\mathcal{T}^\eps_0(h_n)\|_{\mathcal{H}_0} + \|\mathcal{T}^\eps_0(h_n)-\mathcal{T}^\eps_0(h)\|_{\mathcal{H}_0} + \|\mathcal{T}^\eps_0(h)-\mathcal{T}^l_0(h)\|_{\mathcal{H}_0} \nonumber \\
&\leqslant \eps \| h_n\|_{\mathcal{H}_0} + \|\mathcal{T}^\eps_0(h_n)-\mathcal{T}^\eps_0(h)\|_{\mathcal{H}_0} + \eps \| h\|_{\mathcal{H}_0} \nonumber \\
& \leqslant 2\eps + \|\mathcal{T}^\eps_0(h_n)-\mathcal{T}^\eps_0(h)\|_{\mathcal{H}_0}.
 \end{align}
Let us now prove that we have the strong convergence $\|\mathcal{T}^\eps_0(h_n)-\mathcal{T}^\eps_0(h)\|_{\mathcal{H}_0}\rightarrow 0$.\\
For that purpose, we will use Rellich's theorem for the sequence $ \mathtt{H}_n^\eps$ defined by  $ \mathtt{H}_n^\eps:=g^\eps_1 \frac{h_n}{\langle v \rangle^2} - \langle h_n,\Phi \rangle g^\eps_2$. Indeed, it is uniformly bounded in $H^1_{\langle v \rangle^2}$ since we have:
\begin{align*}
\int_{\RR^d}\langle v \rangle^2 |\mathtt{H}_n^\eps|^2 \mathrm{d}v &\leqslant 2\int_{\RR^d} \bigg( |g_1^\eps|^2 \frac{|h_n|^2}{\langle v \rangle^2} + \bigg\|\frac{h_n}{\langle v \rangle}  \bigg\|_2^2 \|\langle v \rangle \Phi \|_2^2 \langle v \rangle^2 |g_2^\eps|^2 \bigg) \mathrm{d}v \\
&\leqslant 2 \bigg(\|g_1^\eps\|_\infty^2 + \|\Phi\|_{L^2_{\langle v \rangle^2}}^2 \|g_2^\eps\|_{L^2_{\langle v \rangle^2}}^2\bigg)\|h_n\|_{\mathcal{H}_0}^2 \lesssim 1
\end{align*}
and 
\begin{align*}
\int_{\RR^d} \langle v \rangle^2 \big| \nabla_v \mathtt{H}_n^\eps\big|^2 \mathrm{d}v &= \int_{\RR^d} \langle v \rangle^2 \bigg| \nabla_v g_1^\eps  \frac{h_n}{\langle v \rangle^2} +  \frac{g_1^\eps}{\langle v \rangle^2} \nabla_v h_n -2\frac{v}{\langle v \rangle^2}g_1^\eps \frac{h_n}{\langle v \rangle^2}- \langle h_n,\Phi \rangle \nabla_v g_2^\eps\bigg|^2 \mathrm{d}v \\
&\lesssim \bigg(\|g_1^\eps\|_{W^{1,\infty}}^2 + \|\Phi\|_{L^2_{\langle v \rangle^2}}^2 \|g_2^\eps\|_{H^1_{\langle v \rangle^2}}^2\bigg)\|h_n\|_{\mathcal{H}_0}^2 \lesssim 1  ,
\end{align*}
where $g_1^\eps$ and $g_2^\eps$ are uniformly bounded in $W^{1,\infty}$ and $H^1_{\langle v \rangle^2}$ respectively,  and since $\|\Phi\|_{L^2_{\langle v \rangle^2}}\leqslant 1$ and $\| h_n \|_{\mathcal{H}_0}\leqslant 1$.\\
 Then, there exists $\mathtt{H}^\eps\in H^1_{\langle v \rangle^2}$ such that $\langle v \rangle \mathtt{H}_n^\eps \longrightarrow \langle v \rangle \mathtt{H}^\eps$ in $L^2(K)$, up to a subsequence, for all $K \subset \RR^d$ bounded, in particular for $K=B(0,R_\eps)$, where $R_\eps>0$ is such that $$ \mathrm{supp}(g^\eps_1)\cup\mathrm{supp}(g^\eps_2)\subset B(0,R_\eps) .$$ 
The limit $\mathtt{H}^\eps$ can be identified as the unique limit in $\mathcal{D}'(\RR^d)$, $\mathtt{H}^\eps = g^\eps_1 \frac{h}{\langle v \rangle^2} - \langle h,\Phi \rangle \Phi$. So for all $\eps'>0$, there exists $N_{\eps'} \in \NN$ such that, for all $n \geqslant N_{\eps'}$ we have: $\| \mathtt{H}^\eps_n - \mathtt{H}^\eps\|_{L^2_{\langle v \rangle^2}} \leqslant \frac{\eps'}{3C}$. Therefore, for $\eps<\frac{\eps'}{3}$ and $n\geqslant N_{\eps'}$ we obtain, thanks to \eqref{T^l(h_n)-T^l(h)} and the inequality $\|\mathcal{T}^l_0(h)\|_{\mathcal{H}_0} \leqslant C\|\mathtt{H}\|_{L^2_{\langle v \rangle^2}}$, that:
\begin{align*}
\|\mathcal{T}^l_0(h_n)-\mathcal{T}^l_0(h)\|_{\mathcal{H}_0} \leqslant 2\eps + C \| \mathtt{H}^\eps_n - \mathtt{H}^\eps\|_{L^2_{\langle v \rangle^2}} \leqslant \eps'.
 \end{align*}
Hence the compactness of $\mathcal{T}^l_0$ holds.
\ep
\begin{proposition}[Assumptions of the implicit function theorem]\label{hypotheses de TFI}
\item \begin{enumerate}
\item The map $F(\lambda,\eta,\cdot)=Id-\mathcal{T}_{\lambda,\eta}$ is continuous in $\mathcal{H}_0$ uniformly with respect to $\lambda \mbox{ and } \eta$. Moreover, there exists $c>0$, independent of $\lambda$ and $\eta$ such that
$$ \|F(\lambda,\eta,h_1)-F(\lambda,\eta,h_2)\|_{\mathcal{H}_0} \leqslant c \|h_1-h_2\|_{\mathcal{H}_0}, \quad \forall h_1, h_2 \in \mathcal{H}_0, \forall\eta\in[0,\eta_0], \forall |\lambda|\leqslant \lambda_0 .$$
\item The map $F$ is continuous with respect to $\lambda$ and $\eta$ and we have
$$ \underset{\eta \to \eta'}{\lim}\|F(\lambda,\eta,h)-F(\lambda,\eta',h)\|_{\mathcal{H}_0} = \underset{\lambda \to \lambda'}{\lim} \|F(\lambda,\eta,h)-F(\lambda',\eta,h)\|_{\mathcal{H}_0} = 0,  \quad \forall h \in \mathcal{H}_0.$$
\item The map $F(\lambda,\eta,\cdot)$ is differentiable in $\mathcal{H}_0$. Moreover,
 $$ \frac{\partial F}{\partial h}(\lambda,\eta,\cdot)=Id-\mathcal{T}_{\lambda,\eta}^l, \quad \forall |\lambda|\leqslant\lambda_0,\forall \eta\in[0,\eta_0] .$$
\item We have $\ F(0,0,M)=0 \ $ and $\ \frac{\partial F}{\partial h}(0,0,M)$ is invertible.
\end{enumerate}
\end{proposition}

\noindent \bp 1. Let $h_1, h_2 \in \mathcal{H}_0$. Let $\eta\in[0,\eta_0]$ and $|\lambda|\leqslant \lambda_0$ with $\eta_0$ and $\lambda_0$ small enough. Then,
\begin{align*}
\|F(\lambda,\eta,h_1)-F(\lambda,\eta,h_2)\|_{\mathcal{H}_0} &\leqslant \|(h_1-h_2)+\mathcal{T}_{\lambda,\eta}^l(h_1-h_2)\|_{\mathcal{H}_0} \\
&\leqslant \big(1+C\big)\|h_1-h_2\|_{\mathcal{H}_0} .
\end{align*}
2. The proof of this point is a direct consequence of the second point of Lemma \ref{continuite de T_eta}.  \\
3. The third point is immediate since $\mathcal{T}_{\lambda,\eta}$ is an affine map with respect to $h$. \\
4.  Recall that $L_0=\tilde Q$ is the inverse of $T_0$ and $V:=\tilde W - W $. Thus, since we have
$$ L_0\big(F(0,0,M) \big) = L_0\big(M-T_0[VM]\big) = [\tilde Q - V](M)= Q(M)=0 .$$
Then, we obtain $\ F(0,0,M) = 0$, thanks to the injectivity of $L_0$. \\
For the differential, we have $\frac{\partial F}{\partial h}(0,0,M)= Id-\mathcal{T}_0^l$. By the Fredholm Alternative, this point is true if $\ \mathrm{Ker}(Id-\mathcal{T}_0^l)=\{0\}$. Let $h\in \mathcal{H}_0$ such that $h-\mathcal{T}_0^l(h) := h-T_0\big[Vh-\langle h,\Phi\rangle\Phi\big] = 0$. Applying the operator $L_0=\tilde Q$ to this last equality we obtain 
$$ \tilde Q(h) - Vh + \langle h,\Phi\rangle\Phi =  Q(h)+ \langle h,\Phi\rangle\Phi = 0 . $$ 
 Integrating this last equation against $M$ and using the fact that $\langle \Phi,M \rangle  =1$, we get
$$ 0 =  \big\langle Q(h)+ \langle h,\Phi\rangle\Phi,M \big\rangle = \langle h,\Phi\rangle \langle \Phi,M \rangle  = \langle h,\Phi\rangle . $$
Therefore, $h$ is solution to $Q(h)=0$. Then, there exists $c_1, c_2 \in \CC$ such that $h=c_1M+c_2Z$. Since $h \in \mathcal{H}_0$ and $Z \notin \mathcal{H}_0$ then, $c_2=0$ and $h=c_1 M$. Thus, $\langle h,\Phi\rangle =c_1=0$.  Hence, $h=0$. This completes the proof of the Proposition.
\ep
\begin{theorem}[Existence of solutions with constraint]\label{thm d'existence} 
 There is a unique function $M_{\lambda,\eta}$ in $\mathcal{H}_0$ solution to the penalized equation
\begin{equation}\label{eq de M_eta + cont.}
[-\Delta_v+ W(v)+  \mathrm{i}  \eta v -\lambda\eta^{\frac{2}{3}}]M_{\lambda,\eta}(v)= b(\lambda,\eta)\Phi(v), \quad v \in \RR^d .
\end{equation} 
where $\ds b(\lambda,\eta):=\langle N_{\lambda,\eta},\Phi \rangle$ with $N_{\lambda,\eta} := M_{\lambda,\eta}-M$. Moreover,
\begin{equation}\label{M_lambda,eta-M_0-->0 dans H_0}
 \| N_{\lambda,\eta} \|_{\mathcal{H}_0} = \| M_{\lambda,\eta} - M \|_{\mathcal{H}_0} \underset{\eta \rightarrow 0}{\longrightarrow}0.
\end{equation}
\end{theorem}

\noindent \bp By Proposition \ref{hypotheses de TFI}, $F$ satisfies the assumptions of the implicit function theorem around the point $(0,0,M)$. Then, there exists $\lambda_0, \eta_0>0$ small enough, there exists a unique function $\mathcal{M} : \{|\lambda|\leqslant\lambda_0\}\times [0,\eta_0]\longrightarrow \mathcal{H}_0$, continuous with respect to $\lambda$ and $\eta$ such that \begin{center}
$ F(\lambda,\eta,\mathcal{M}(\lambda,\eta))=0$, for all $(\lambda,\eta) \in \{|\lambda|<\lambda_0\}\times [0,\eta_0[$.
\end{center}
Let's denote $M_{\lambda,\eta}:=\mathcal{M}(\lambda,\eta)$. The function $M_{\lambda,0}$ does not depend on $\lambda$ and the continuity of $\mathcal{M}$ with respect to $\eta$ implies that 
$$
\underset{\eta\rightarrow 0}{\lim}\|M_{\lambda,\eta}-M_{\lambda,0}\|_{\mathcal{H}_0} = \underset{\eta\rightarrow 0}{\lim} \|M_{\lambda,\eta}-M\|_{\mathcal{H}_0} = 0 .
$$
\ep
\begin{remark}\label{rmq sur thm d'existence}
\item \begin{enumerate}
\item Since $\Phi(-v)=\Phi(v)$ for all $v\in \RR^d$ and the function $\bar{M_{\bar \lambda,\eta}}(-v_1,v')$ satisfies equation \eqref{eq de M_eta + cont.} then, by uniqueness, $\bar{M_{\bar \lambda,\eta}}(-v_1,v')$ is solution to \eqref{eq de M_eta + cont.} and the following symmetry 
\begin{equation}\label{bar M_eta(-v)=M_eta(v)}
\bar{M_{\bar \lambda,\eta}}(-v_1,v') = M_{\lambda,\eta}(v_1,v') 
\end{equation}
holds for all $(v_1,v') \in \RR\times\RR^{d-1}$, $\eta \in [0,\eta_0]$ and $|\lambda| \leqslant \lambda_0$.
\item The sequence $|b(\lambda,\eta)|$ is uniformly bounded with respect to $\lambda$ and $\eta$ since $ \ds |b(\lambda,\eta)| \underset{\eta \rightarrow 0}{\longrightarrow}0$, which we obtain by the Cauchy-Schwarz inequality and  limit \eqref{M_lambda,eta-M_0-->0 dans H_0}:
\begin{equation}
|b(\lambda,\eta)|= |\langle N_{\lambda,\eta},\Phi \rangle| \leqslant \bigg\|\frac{N_{\lambda,\eta}}{\langle v \rangle} \bigg\|_2 \| \langle v \rangle \Phi\|_2 \leqslant \| N_{\lambda,\eta} \|_{\mathcal{H}_0} \| \langle v \rangle \Phi\|_2 \underset{\eta \rightarrow 0}{\longrightarrow}0 .
\end{equation}
\end{enumerate}
\end{remark}
\section{Existence of the eigen-solution $\big(\mu(\eta),M_{\mu,\eta}\big)$}
The aim of this section is to prove Theorem \ref{main}. It is composed of three subsections. In the first one, we establish some $L^2$ estimates. The second one is devoted to the study of the constraint and the existence of the eigen-solution $(\mu(\eta),M_\eta)$. Finally, in the last subsection, we give an approximation of the eigenvalue and its relation with the diffusion coefficient.

\subsection{$L^2$ estimates for the solution $M_{\lambda,\eta}$}
In this subsection, we will establish some $L^2$ estimates for the solution of the penalised equation \eqref{eq penalisee2}. 

\begin{proposition}\label{propostion estimation L^2} Let $ \eta_0>0$ and $\lambda_0>0$ small enough. Let $M_{\lambda,\eta}$ be the solution of the penalised equation \eqref{eq de M_eta + cont.}. Then, for all $\eta \in [0,\eta_0]$ and for all $\lambda \in \CC$ such that $|\lambda|\leqslant \lambda_0$, one has
\begin{enumerate}
\item For all $\gamma > \frac{d}{2}$, the function $M_{\lambda,\eta}$ is uniformly bounded, with respect to $\lambda$ and $\eta$, in $L^2(\RR^d,\CC)$. Moreover, the following estimate holds
\begin{equation}\label{estimation de N_lambda,eta dans L^2}
\|N_{\lambda,\eta}\|_{L^2(\RR^d)}^2 = \|M_{\lambda,\eta}-M\|_{L^2(\RR^d)}^2 \lesssim |\lambda| +\nu_\eta ,
\end{equation}
where $\ds \nu_\eta \underset{\eta \to 0}{\longrightarrow} 0 $.
\item For all $\gamma > \frac{d+1}{2}$, the function $|v_1|^{\frac{1}{2}} M_{\lambda,\eta}$ is uniformly bounded, with respect to $\lambda$ and $\eta$, in $L^2(\RR^d,\CC)$.  
\end{enumerate}
\end{proposition}
\bp We are going to prove the first point, the second is done in a similar way. Let denote $v:=(v_1,v') \in \RR\times\RR^{d-1}$. The proof of this Proposition is given in four steps and the idea is as follows: first, we decompose $\RR^d$ into two parts, $\RR^d = \{ |v_1|\leqslant s_0\eta^{-\frac{1}{3}}\} \cup \{|v_1| \geqslant s_0\eta^{-\frac{1}{3}}\}$, small/medium and large velocities.  In the first step, using the equation of $M_{\lambda,\eta}$, we estimate the norm of $M_{\lambda,\eta}$ for large velocities to get
$$ \|M_{\lambda,\eta} \|^2_{L^2(\{|v_1| \geqslant s_0\eta^{-\frac{1}{3}}\})}  \leqslant \nu_1 \|M_{\lambda,\eta} \|^2_{L^2(\{ |v_1|\leqslant s_0\eta^{-\frac{1}{3}}\})}  + c_1  , $$
where $\nu_1$ and $c_1$ depend on $s_0$, $\lambda$ and $\eta$. To estimate $\|M_{\lambda,\eta} \|_{L^2(\{ |v_1|\leqslant s_0\eta^{-\frac{1}{3}}\})}$, it is enough to estimate $\| N_{\lambda,\eta} \|_{L^2(\{ |v_1|\leqslant s_0\eta^{-\frac{1}{3}}\})}$ since $M$ belongs to $L^2$, which is the purpose of steps two and three. In step 2, using a Poincar\'e type inequality, we show that
$$\|N_{\lambda,\eta} \|^2_{L^2(\{ |v_1| \leqslant s_0\eta^{-\frac{1}{3}}\leqslant |v'|\})}  \leqslant C_1 \|M_{\lambda,\eta} \|^2_{L^2(\{|v_1| \geqslant s_0\eta^{-\frac{1}{3}}\})}  + c_2  , $$
where $C_1$ is a positive constant and $c_2$ depends on $s_0$, $\lambda$ and $\eta$. Then, in the third step, using the Hardy-Poincar\'e inequality, we prove that
 $$\|N_{\lambda,\eta} \|^2_{L^2(\{|v|\leqslant s_0\eta^{-\frac{1}{3}}\})}  \leqslant \nu_2 \|N_{\lambda,\eta} \|^2_{L^2(\{ |v_1|\leqslant s_0\eta^{-\frac{1}{3}}\})} + \nu_3 \|M_{\lambda,\eta} \|^2_{L^2(\{|v_1| \geqslant s_0\eta^{-\frac{1}{3}}\})}  + c_3  ,$$  
with $\nu_2$, $\nu_3$ and $c_3$ depend on $s_0$, $\lambda$ and $\eta$. The last step is left for the conclusion: we first fix $s_0$ large enough, then $|\lambda|$ small enough, then $\eta$ small enough, we obtain $\nu_2\leqslant \frac{1}{4}$, $\nu_3\leqslant \frac{1}{4}$ and $\nu_1 \big(C_1+\frac{\nu_3}{1-\nu_2}\big) \leqslant \frac{1}{2}$, which allows us to conclude.\\

\noindent Before starting the proof, we will define some sets to simplify the notations and avoid long expressions. 
We set: $ A_\eta:=\{|v_1|\leqslant s_0\eta^{-\frac{1}{3}}\}$ (resp., $\tilde A_\eta:=\{|v_1| \leqslant 2s_0\eta^{-\frac{1}{3}}\}$), $B_\eta:=\{|v|\leqslant s_0\eta^{-\frac{1}{3}}\}$, $C_\eta:=\{|v_1| \leqslant s_0\eta^{-\frac{1}{3}}\leqslant |v'|\}$ (resp., $\tilde C_\eta:=\{| v_1|\leqslant 2s_0\eta^{-\frac{1}{3}}\leqslant 2|v'|\}$) and $D_\eta:=\{|v_1|\geqslant \frac{s_0}{2}\eta^{-\frac{1}{3}}\}$.  
\begin{figure}[!ht]
\centerline{{\includegraphics[scale=0.4]{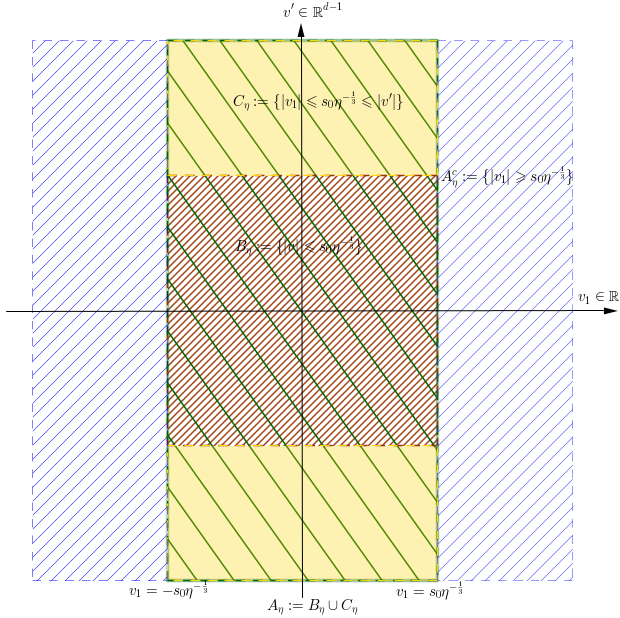}}}
\caption{Decomposition of $\RR^d$ into $A_\eta$ and $A_\eta^c$.}\label{figure}
\end{figure}

\ni  The part $A_\eta^c$ is represented by the blue zone, while the part $A_\eta$, in green stripes, is broken down into two other parts: the brown zone $B_\eta$ for $|v '|$ small, and the yellow zone $C_\eta$ for $|v'|$ large. \\

\ni The parts $\tilde A_\eta$, $\tilde C_\eta$ and $D_\eta$ are an extensions ``in the direction of $v_1$'' of the parts $A_\eta$, $C_\eta$ and $A_\eta^c$ respectively, and are not shown in the figure above.  \\
 
\noindent \textbf{Step 1: Estimation of $\|M_{\lambda,\eta}\|_{L^2(A_\eta^c)}$.} We summarize this step in the following inequality: 
\begin{equation}\label{M_eta sur A^c}
\|M_{\lambda,\eta}\|^2_{L^2(A_\eta^c)} \leqslant \frac{1}{s_0^2} \big\|\eta^{\frac{1}{3}} v_1 M_{\lambda,\eta} \big\|^2_{L^2(A_\eta^c)} \lesssim \frac{1}{s_0^3}\bigg( \|N_{\lambda,\eta}\|^2_{L^2(A_\eta)}+c_1^\eta\bigg) ,
\end{equation} 
where $c_1^\eta = c_1(\lambda,\eta,s_0) = s_0^{-2\delta} \eta^{\frac{2\delta}{3}}(|b(\lambda,\eta)|+1) \| |v_1|^{\delta}M\|^2_{L^2(\RR^d)}$ where $\delta$ can be chosen as follows $\delta := \frac{1}{2}(\gamma-\frac{d}{2})$ to ensure that $|v|^\delta M$ belongs to $L^2$.\\

\noindent $\bullet$ \textbf{Estimation of $\|\eta^{\frac{1}{3}} v_1 M_{\lambda,\eta}\|_{L^2(A_\eta^c)}$.} In order to localize the velocities on the part $A_\eta^c$ and to be able to use the equation of $M_{\lambda,\eta}$ and make integrations by part, we introduce the function $\chi_\eta$ defined by: $\chi_{\eta}(v_1) := \chi\big(\frac{v_1}{s_0\eta^{-1/3}}\big)$, where $\chi \in C^\infty(\RR)$ is such that $0 \leqslant \chi \leqslant 1$, $\chi \equiv 0$ on $B(0,\frac{1}{2})$ and $\chi \equiv 1$ outside of $B(0,1)$.  Then, one has: $\|\eta^{\frac{1}{3}} v_1 M_{\lambda,\eta}\|_{L^2(A_\eta^c)} \leqslant \|\eta^{\frac{1}{3}} v_1 \chi_\eta M_{\lambda,\eta}\|_{L^2(D_\eta)}$. Now, multiplying the equation of $M_{\lambda,\eta}$ by $v_1 \bar M_{\lambda,\eta} \chi_\eta^2$, integrating it over $D_\eta $ and taking the imaginary part, we get:
\begin{align*}
\big\| \eta^{\frac{1}{3}} v_1 \chi_\eta M_{\lambda,\eta} \big\|^2_{L^2(D_\eta)} &= -\eta^{-\frac{1}{3}} \mathrm{Im} \bigg( \int_{D_\eta} Q(M_{\lambda,\eta}) v_1 \bar M_{\lambda,\eta} \chi_\eta^2 \mathrm{d}v\bigg)+ \mathrm{Im} \bigg(  \lambda \int_{D_\eta}  \eta^{\frac{1}{3}} v_1 |M_{\lambda,\eta} \chi_\eta|^2\mathrm{d}v\bigg)  \\
& \quad - \eta^{-\frac{1}{3}} \mathrm{Im} \bigg( b(\lambda,\eta)\int_{D_\eta} \Phi v_1 \bar M_{\lambda,\eta} \chi_\eta^2\mathrm{d}v \bigg)  \\
&=: -E^\eta_1 + E^\eta_2 + E^\eta_3.
\end{align*}
Let's start with $E^\eta_2$ and $E^\eta_3$ which are simpler. \\
$\bullet$ \textbf{Estimation of $E^\eta_2$:} For this term, we just use the fact that on $D_\eta$: $\frac{s_0}{2}\leqslant\eta^{\frac{1}{3}}|v_1|$. Thus, 
\begin{equation}\label{E^eta_2}
|E^\eta_2|:=\bigg|\mathrm{Im} \bigg(  \lambda \int_{D_\eta}  \eta^{\frac{1}{3}} v_1 |M_{\lambda,\eta} \chi_\eta|^2\mathrm{d}v\bigg)\bigg| \leqslant \frac{2|\lambda|}{s_0} \big\| \eta^{\frac{1}{3}} v_1 \chi_\eta M_{\lambda,\eta} \big\|^2_{L^2(D_\eta)} .
\end{equation}
$\bullet$ \textbf{Estimation of $E^\eta_3$:} First of all,  since $\Phi(v) := \big(\int_{\RR^d} \langle v \rangle^{-2-2\gamma}\mathrm{d}v\big)^{-1}\frac{M(v)}{\langle v \rangle^2}$ then: 
\begin{equation}\label{Phi zeta_eta}
\Phi(v)\chi_{\{|v_i|\geqslant s_0\eta^{-\frac{1}{3}}\}}(v)\leqslant C s_0^{-2}\eta^{\frac{2}{3}}M(v) .
\end{equation}
In particular,
\begin{equation}\label{Phi zeta_eta(v_1)}
\Phi(v)\chi_\eta(v_1)\leqslant C s_0^{-2}\eta^{\frac{2}{3}}M(v) .
\end{equation}
Similarly, we have:
\begin{equation}\label{3.6}
\| M \|_{L^2(\{|v_i|\geqslant  s_0\eta^{-\frac{1}{3}}\})} \leqslant s_0^{-\delta}\eta^{\frac{\delta}{3}} \big\| |v|^\delta M \big\|_{L^2(\RR^d)} .
\end{equation}
Then, using \eqref{Phi zeta_eta(v_1)}, we get
$$
|E^\eta_3|:=\bigg|\eta^{-\frac{1}{3}} \mathrm{Im} \bigg( b(\lambda,\eta)\int_{D_\eta} \Phi v_1 \bar M_{\lambda,\eta} \chi_\eta^2\mathrm{d}v \bigg)\bigg|  \leqslant 4 \frac{|b(\lambda,\eta)|}{s_0^2} \|\chi_\eta M\|_{L^2(D_\eta)} \big\| \eta^{\frac{1}{3}} v_1 \chi_\eta M_{\lambda,\eta} \big\|_{L^2(D_\eta)} .
$$
Finally, by  inequality \eqref{3.6}
\begin{equation}\label{E^eta_3}
|E^\eta_3| \leqslant 2 \frac{|b(\lambda,\eta)|}{s_0^2} \bigg( \big\| \eta^{\frac{1}{3}} v_1 \chi_\eta M_{\lambda,\eta} \big\|^2_{L^2(D_\eta)} + 4s_0^{-2\delta}\eta^{\frac{2\delta}{3}} \big\| |v|^\delta M \big\|_{L^2(\RR^d)}^2 \bigg)  .
\end{equation}
$\bullet$ \textbf{Estimation of $E^\eta_1$:} By an integration by parts, we write 
$$ E^\eta_1:= \eta^{-\frac{1}{3}} \mathrm{Im}  \int_{D_\eta} Q(M_{\lambda,\eta}) v_1 \bar M_{\lambda,\eta} \chi_\eta^2  \mathrm{d}v = \eta^{-\frac{1}{3}} \mathrm{Im}  \int_{D_\eta} \big[\chi_\eta \bar M_{\lambda,\eta} + 2  v_1\chi_\eta' \bar M_{\lambda,\eta}\big] \pa_{v_1}\bigg(\frac{M_{\lambda,\eta}}{M}\bigg)M\chi_\eta \mathrm{d}v .$$
Thus, by Cauchy-Schwarz
$$ |E^\eta_1| \leqslant \eta^{-\frac{1}{3}}\bigg\|\pa_{v_1}\bigg(\frac{M_{\lambda,\eta}}{M}\bigg)M\chi_\eta\bigg\|_{L^2(D_\eta)}\bigg(\|\chi_\eta M_{\lambda,\eta}\|_{L^2(D_\eta)}+2\|v_1\chi_\eta' M_{\lambda,\eta}\|_{L^2(D_\eta)}\bigg) .$$
Since $\chi_\eta' \equiv 0$ except on: $D_\eta \setminus A_\eta^c =\{\frac{s_0}{2}\eta^{-\frac{1}{3}}\leqslant |v_1|\leqslant s_0\eta^{-\frac{1}{3}}\} \subset A_\eta:=\{ |v_1|\leqslant s_0\eta^{-\frac{1}{3}}\}$. Then, 
$$\|v_1\chi_\eta' M_{\lambda,\eta}\|_{L^2(D_\eta)}=\|v_1\chi_\eta' M_{\lambda,\eta}\|_{L^2(D_\eta\setminus A_\eta^c)}\leqslant C\| M_{\lambda,\eta}\|_{L^2(D_\eta\setminus A_\eta^c)} , $$
where $C=\underset{\frac{1}{2}\leqslant |t|\leqslant 1}{\sup} |t\chi'(t)|$. Also, we have: $\|\chi_\eta M_{\lambda,\eta}\|_{L^2(D_\eta)} \leqslant \frac{1}{s_0}\|\eta^{\frac{1}{3}} v_1 \chi_\eta M_{\lambda,\eta}\|_{L^2(D_\eta)}$. Thus,
$$ 
|E^\eta_1| \leqslant \eta^{-\frac{1}{3}}\bigg\|\pa_{v_1}\bigg(\frac{M_{\lambda,\eta}}{M}\bigg)M\chi_\eta\bigg\|_{L^2(D_\eta)}\bigg(\frac{1}{s_0}\|\eta^{\frac{1}{3}} v_1 \chi_\eta M_{\lambda,\eta}\|_{L^2(D_\eta)}+C \|M_{\lambda,\eta}\|_{L^2(D_\eta\setminus A_\eta^c)}\bigg) .
$$
Finally, by Young's inequality:
\begin{equation}\label{E^eta_11}
|E^\eta_1| \lesssim s_0\eta^{-\frac{2}{3}}\bigg\|\pa_{v_1}\bigg(\frac{M_{\lambda,\eta}}{M}\bigg)M\chi_\eta\bigg\|^2_{L^2(D_\eta)} + \frac{1}{s_0^3}\|\eta^{\frac{1}{3}} v_1 \chi_\eta M_{\lambda,\eta}\|^2_{L^2(D_\eta)}+\frac{1}{s_0} \|M_{\lambda,\eta}\|^2_{L^2(D_\eta\setminus A_\eta^c)} .
\end{equation}
It remains to estimate $\big\|\pa_{v_1}\big(\frac{M_{\lambda,\eta}}{M}\big)M\chi_\eta\big\|_{L^2( D_\eta)}$. For this, one has
\begin{align}\label{F_1 + F_2}
 \bigg\|\pa_{v_1}\bigg(\frac{M_{\lambda,\eta}}{M}\bigg)M\chi_\eta\bigg\|^2_{L^2(D_\eta)} &\leqslant \bigg\|\nabla_v\bigg(\frac{M_{\lambda,\eta}}{M}\bigg)M\chi_\eta\bigg\|^2_{L^2(D_\eta)}  \nonumber \\
&= \Re \int_{D_\eta} \bigg[ Q(M_{\lambda,\eta})\bar M_{\lambda,\eta}\chi_\eta^2  + 2\zeta_\eta'\zeta_\eta\frac{\bar M_{\lambda,\eta}}{M}\pa_{v_1}\bigg(\frac{M_{\lambda,\eta}}{M}\bigg)M^2 \bigg]\mathrm{d}v  \nonumber \\
&=: F^\eta_1 + F^\eta_2.
\end{align}
By integrating the equation of $M_{\lambda,\eta}$, multiplied by $\bar M_{\lambda,\eta} \chi_\eta^2$, over $D_\eta$, and using \eqref{Phi zeta_eta}, we obtain
\begin{equation}\label{F^eta_1}
|F^\eta_1|  \lesssim \frac{\eta^{\frac{2}{3}}}{s_0^2}\bigg[\bigg(|\lambda| + \frac{|b(\lambda,\eta)|}{s_0^2}\bigg)\big\| \eta^{\frac{1}{3}} v_1 \chi_\eta M_{\lambda,\eta} \big\|^2_{L^2(D_\eta)} + |b(\lambda,\eta)| \| M \|^2_{L^2(D_\eta)} \bigg]  .
\end{equation}
For $F^\eta_2$, by  inequality $(2Cab\leqslant 2C^2a^2+\frac{b^2}{2})$:
\begin{align}
|F^\eta_2| & \leqslant 2\| \chi_\eta' M_{\lambda,\eta}\|_{L^2(D_\eta)} \bigg\|\pa_{v_1}\bigg(\frac{M_{\lambda,\eta}}{M}\bigg)M\chi_\eta\bigg\|_{L^2(D_\eta)}    \nonumber \\
&\leqslant 2C \frac{\eta^{\frac{1}{3}}}{s_0} \|M_{\lambda,\eta}\|_{L^2(D_\eta\setminus A_\eta^c)}\bigg\|\pa_{v_1}\bigg(\frac{M_{\lambda,\eta}}{M}\bigg)M\chi_\eta\bigg\|_{L^2(D_\eta)} \nonumber \\
&\leqslant  C' \frac{\eta^{\frac{2}{3}}}{s_0^2} \|M_{\lambda,\eta}\|^2_{L^2(D_\eta\setminus A_\eta^c)}+\frac{1}{2}\bigg\|\pa_{v_1}\bigg(\frac{M_{\lambda,\eta}}{M}\bigg)M\chi_\eta\bigg\|^2_{L^2(D_\eta)} ,
\end{align}
where $C=\underset{\frac{1}{2}\leqslant |t|\leqslant 1}{\sup} |t\chi'(t)|$ and $C'=2C^2$. Then, we obtain by returning to \eqref{F_1 + F_2}
$$ \bigg\|\pa_{v_1}\bigg(\frac{M_{\lambda,\eta}}{M}\bigg)M\chi_\eta\bigg\|^2_{L^2(D_\eta)} \leqslant |F^\eta_1| + C' \frac{\eta^{\frac{2}{3}}}{s_0^2} \|M_{\lambda,\eta}\|^2_{L^2(D_\eta\setminus A_\eta^c)}+\frac{1}{2}\bigg\|\pa_{v_1}\bigg(\frac{M_{\lambda,\eta}}{M}\bigg)M\chi_\eta\bigg\|^2_{L^2(D_\eta)} .
$$
Therefore,
$$ \bigg\|\pa_{v_1}\bigg(\frac{M_{\lambda,\eta}}{M}\bigg)M\chi_\eta\bigg\|^2_{L^2(D_\eta)} \lesssim |F^\eta_1| + \frac{\eta^{\frac{2}{3}}}{s_0^2} \|M_{\lambda,\eta}\|^2_{L^2(D_\eta\setminus A_\eta^c)} . $$
Hence, from \eqref{F^eta_1}, \eqref{3.6} and the last inequality
\begin{align*}
\bigg\|\pa_{v_1}\bigg(\frac{M_{\lambda,\eta}}{M}\bigg)M\chi_\eta\bigg\|^2_{L^2(D_\eta)} \lesssim \frac{\eta^{\frac{2}{3}}}{s_0^2} &\bigg[\bigg(|\lambda| + \frac{|b(\lambda,\eta)|}{s_0^2}\bigg)\big\| \eta^{\frac{1}{3}} v_1 \chi_\eta M_{\lambda,\eta} \big\|^2_{L^2(D_\eta)} + \|M_{\lambda,\eta}\|^2_{L^2(D_\eta\setminus A_\eta^c)} \\
&+ s_0^{-2\delta}\eta^{\frac{2\delta}{3}} |b(\lambda,\eta)| \big\| |v_1|^\delta M \big\|^2_{L^2(\RR^d)} \bigg] .
\end{align*}
Which implies, by inequality \eqref{E^eta_11}, that
\begin{align}
|E^\eta_1| \lesssim \frac{1}{s_0}&\bigg[\bigg(|\lambda| + \frac{1+|b(\lambda,\eta)|}{s_0^2}\bigg)\big\| \eta^{\frac{1}{3}} v_1 \chi_\eta M_{\lambda,\eta} \big\|^2_{L^2(D_\eta)}  + \|M_{\lambda,\eta}\|^2_{L^2(D_\eta\setminus A_\eta^c)}   \nonumber \\
&+ s_0^{-2\delta} \eta^{\frac{2\delta}{3}}|b(\lambda,\eta)| \big\| |v_1|^\delta M \big\|^2_{L^2(\RR^d)} \bigg] . \label{E^eta_1}
\end{align}
Thus, by summing the inequalities obtained from $E^\eta_1$, $E^\eta_2$ and $E^\eta_3$, namely \eqref{E^eta_1}, \eqref{E^eta_2} and \eqref{E^eta_3} respectively, we obtain
\begin{align*}
\big\| \eta^{\frac{1}{3}} v_1 \chi_\eta M_{\lambda,\eta} \big\|^2_{L^2(D_\eta)} \lesssim \frac{1}{s_0}&\bigg[\bigg(|\lambda| + \frac{1+|b(\lambda,\eta)|}{s_0^2}\bigg)\big\| \eta^{\frac{1}{3}} v_1 \chi_\eta M_{\lambda,\eta} \big\|^2_{L^2(D_\eta)}  + \|M_{\lambda,\eta}\|^2_{L^2(D_\eta\setminus A_\eta^c)}    \\
&+ s_0^{-2\delta} \eta^{\frac{2\delta}{3}}|b(\lambda,\eta)| \big\| |v_1|^\delta M \big\|^2_{L^2(\RR^d)} \bigg] .
\end{align*}
Hence the following estimate
\begin{equation}\label{eta^(1/3) v_1 M_eta sur A^c}
\big\| \eta^{\frac{1}{3}} v_1 \chi_\eta M_{\lambda,\eta} \big\|^2_{L^2(D_\eta)} \lesssim \frac{1}{s_0} \bigg( \|M_{\lambda,\eta}\|^2_{L^2(D_\eta\setminus A_\eta^c)} + s_0^{-2\delta} \eta^{\frac{2\delta}{3}}|b(\lambda,\eta)| \big\| |v_1|^\delta M \big\|^2_{L^2(\RR^d)} \bigg) 
\end{equation}
holds true for $s_0>0$ large enough and for all $|\lambda|\leqslant \lambda_0$ and $\eta\in[0,\eta_0]$, with $\lambda_0$ and $\eta_0$ small enough.  Finally, \eqref{M_eta sur A^c} comes from the previous inequality \eqref{eta^(1/3) v_1 M_eta sur A^c}, and since $D_\eta\setminus A_\eta^c \subset A_ \eta$ and $D_\eta\setminus A_\eta^c \subset D_\eta$ implies that,
$$ \|M_{\lambda,\eta}\|^2_{L^2(D_\eta\setminus A_\eta^c)} \leqslant \| N_{\lambda,\eta}\|^2_{L^2(A_\eta)} + \|M\|^2_{L^2(D_\eta)} \leqslant \| N_{\lambda,\eta}\|^2_{L^2(A_\eta)} +4s_0^{-2\delta} \eta^{\frac{2\delta}{3}} \big\| |v_1|^\delta M \big\|^2_{L^2(\RR^d)} . $$

\noindent \textbf{Step 2: Estimation of $\|N_{\lambda,\eta}\|_{L^2(C_\eta)}$}. In this step, we will establish the following inequality:
\begin{equation}\label{N_eta sur C}
\|N_{\lambda,\eta}\|^2_{L^2(C_\eta)} \lesssim  \|N_{\lambda,\eta}\|^2_{L^2(A_\eta^c)} + c_2^\eta ,
\end{equation}
where $c_2^\eta := s_0^{-\delta} \eta^{\frac{\delta}{3}}\big(s_0^2 |\lambda|+s_0^3+|b(\lambda,\eta)| \big) \| |v|^\delta M \|^2_{L^2(\RR^d)}$, and where we recall that $\delta := \frac{1}{2}(\gamma-\frac{d}{2})$, $C_\eta := \{|v_1|\leqslant s_0\eta^{-\frac{1}{3}}\leqslant|v'|\}$ and $A_\eta^c := \{|v_1| \geqslant s_0\eta^{-\frac{1}{3}}\}$. 
We start with the following Lemma:
\begin{lemma}[Poincar\'e-type inequality]\label{lemme Poincare2}
Let $R>0$ be fixed and let $C_R$ be the set defined by: $C_R:=\{ v\in \RR^d; \ |v_1|\leqslant R \leqslant |v'| \}$. Then, there exists a constant $C>0$ such that, for any function $\psi$ in the space $\mathcal{H}:=\big\{\int_{C_R} \big|\pa_{v_1}\big(\frac{\psi} {M}\big)\big|^2 M^2\mathrm{d}v < \infty ; \ \psi(-R,\cdot)=\psi(R,\cdot)=0 \big\}$, the inequality 
\begin{equation}\label{Poincare2}
\|\psi\|^2_{L^2(C_R)} \leqslant C R^2 \big\|\pa_{v_1}\big(\frac{\psi}{M}\big)M\big\|^2_{L^2(C_R)}
\end{equation}
holds true.
\end{lemma}
\bpl \ref{lemme Poincare2}. We have for $\psi \in \mathcal{H}$:
$\ \ds\frac{\psi}{M} = \int_{-R}^{v_1} \pa_{w_1}\big(\frac{\psi}{M}\big) \mathrm{d}w_1 $. \\
Then, by taking the square and applying the Cauchy-Schwarz inequality, we get:
$$|\psi|^2 \leqslant M^2(v_1,v')\bigg(\int_{-R}^{v_1}\pa_{w_1}\bigg(\frac{\psi}{M}\bigg) \mathrm{d}w_1\bigg)^2  \leqslant \int_{-R}^{R} \frac{M^2(v_1,v')}{M^2(w_1,v')} \mathrm{d}w_1 \int_{-R}^{R} \bigg|\pa_{w_1}\bigg(\frac{\psi}{M}\bigg)M\bigg|^2 \mathrm{d}w_1 . $$
Now, we have for $v_1,w_1 \in [-R,R]$ and $|v'|\geqslant R$, $\frac{M^2(v_1,v')}{M^2(w_1,v')} \lesssim 1$. Therefore,
$$ |\psi|^2 \lesssim R \int_{-R}^{R} \bigg|\pa_{w_1}\bigg(\frac{\psi}{M}\bigg)M\bigg|^2 \mathrm{d}w_1 . $$
Thus, we obtain  inequality \eqref{Poincare2} by integrating the last one over $C_R$.

\epl
Now back to the estimate of $\|N_{\lambda,\eta}\|_{L^2(C_\eta)}$.
Let $\zeta \in C^\infty(\RR)$ such that $0\leqslant \zeta \leqslant 1$, $\zeta \equiv 1$ on $B(0,1)$ and $\zeta \equiv 0$ outside of $B(0,2)$. We define $\zeta_{\eta}$ by: $\zeta_{\eta}(v_1) := \zeta\big(\frac{v_1}{s_0\eta^{-1/3}}\big)$. Then, for $\eta>0$ and  $s_0>0$ fixed, by applying Lemma \ref{lemme Poincare2} for $R=s_0\eta^{-\frac{1}{3}}$, we obtain:
\begin{equation}\label{zeta N_eta sur tilde C}
\|N_{\lambda,\eta}\|^2_{L^2(C_\eta)} \leqslant \|\zeta_\eta N_{\lambda,\eta}\|^2_{L^2(\tilde C_\eta)} \lesssim \ s_0^2 \eta^{-\frac{2}{3}} \bigg\|\pa_{v_1}\bigg(\frac{\zeta_\eta N_{\lambda,\eta}}{M}\bigg)M\bigg\|^2_{L^2(\tilde C_\eta)} ,
\end{equation}
recalling that $\tilde C_\eta := \{|v_1|\leqslant 2s_0\eta^{-\frac{1}{3}} \leqslant 2|v'|\}$. Furthermore, 
$$ \bigg\|\pa_{v_1}\bigg(\frac{\zeta_\eta N_{\lambda,\eta}}{M}\bigg)M\bigg\|^2_{L^2(\tilde C_\eta)} \leqslant \bigg\|\nabla_v\bigg(\frac{\zeta_\eta N_{\lambda,\eta}}{M}\bigg)M\bigg\|^2_{L^2(\tilde C_\eta)} = \Re \int_{\tilde C_\eta} Q(\zeta_\eta N_{\lambda,\eta}) \zeta_\eta \bar N_{\lambda,\eta}\mathrm{d}v . $$
However, 
\begin{align}\label{Re int Q(zetaN)zetaN}
\Re \int_{\tilde C_\eta} Q(\zeta_\eta N_{\lambda,\eta}) \zeta_\eta \bar N_{\lambda,\eta}\mathrm{d}v &= \Re \int_{\tilde C_\eta} \big[Q(N_{\lambda,\eta}) \bar N_{\lambda,\eta} \zeta_\eta^2 - \zeta_\eta \zeta_\eta''|N_{\lambda,\eta}|^2 - 2\zeta_\eta\zeta_\eta'\bar N_{\lambda,\eta}\pa_{v_1}N_{\lambda,\eta}\big]\mathrm{d}v  \nonumber \\
&= \Re \int_{\tilde C_\eta} Q(N_{\lambda,\eta}) \bar N_{\lambda,\eta} \zeta_\eta^2 \mathrm{d}v + \int_{\tilde C_\eta} |\zeta_\eta' N_{\lambda,\eta}|^2 \mathrm{d}v   ,
\end{align}
where we used the fact that $Q(\zeta_\eta N_{\lambda,\eta})=Q(N_{\lambda,\eta})\zeta_\eta-\zeta_\eta''N_{\lambda,\eta}-2\zeta_\eta'\pa_{v_1}N_{\lambda,\eta}$ in the first line, since $Q:=-\frac{1}{M}\nabla_v\big(M^2\nabla\big(\frac{\cdot}{M}\big)\big)=-\Delta_v+W(v)$, and did an integration by parts for the term $\int_{\tilde C_\eta} \zeta_\eta \zeta_\eta''|N_{\lambda,\eta}|^2 \mathrm{d}v$, and used the identity: $\Re \big(\bar f \pa_{v_1}f\big) = \frac{1}{2}\pa_{v_1}|f|^2$ in the second line.\\
To handle $\int_{\tilde C_\eta} |\zeta_\eta' N_{\lambda,\eta}|^2 \mathrm{d}v$, we have:
$$ \int_{\tilde C_\eta} |\zeta_\eta' N_{\lambda,\eta}|^2 \mathrm{d}v = \int_{\tilde C_\eta\setminus B_\eta^c} |\zeta_\eta' N_{\lambda,\eta}|^2 \mathrm{d}v \leqslant \|\zeta'\|^2_{L^\infty(\tilde C_\eta\setminus B_\eta^c)} \|N_{\lambda,\eta}\|^2_{L^2({\tilde C_\eta\setminus B_\eta^c})} , $$
since $\zeta_\eta' \equiv 0$ except on: $\tilde C_\eta\setminus B_\eta^c=\{s_0\eta^{-\frac{1}{3}}\leqslant|v_1|\leqslant 2s_0 \eta^{-\frac{1}{3}}\leqslant 2 |v'|\} \subset A_\eta^c$, and that on $\tilde C_\eta\setminus B_\eta^c$ we have: $|\zeta_\eta'(v_1)|\lesssim \frac {\eta^\frac{1}{3}}{s_0}$. Then,
\begin{equation}\label{zeta' N_eta sur tilde C}
\int_{\tilde C_\eta} |\zeta_\eta' N_{\lambda,\eta}|^2 \mathrm{d}v \lesssim \frac{\eta^\frac{2}{3}}{s_0^2} \|N_{\lambda,\eta}\|^2_{L^2(A_\eta^c)} .
\end{equation}
To handle $\Re \int_{\tilde C_\eta} Q(N_{\lambda,\eta}) \bar N_{\lambda,\eta} \zeta_\eta^2 \mathrm{d}v$, we will proceed as in $E^\eta_1$. Indeed, recall that $N_{\lambda,\eta}$ satisfies the equation:
$$ Q(N_{\lambda,\eta})= (\lambda \eta^{\frac{2}{3}}-  \mathrm{i} \eta v_1) N_{\lambda,\eta} + (\lambda \eta ^{\frac{2}{3}}-  \mathrm{i} \eta v_1) M - b(\lambda,\eta) \Phi . $$
Then, multiplying this equation by $\bar N_{\lambda,\eta} \zeta_\eta^2$ and integrating it over $\tilde C_\eta$, we get
\begin{align}
\bigg|\Re \int_{\tilde C_\eta} Q(N_{\lambda,\eta}) \bar N_{\lambda,\eta} \zeta_\eta^2 \mathrm{d}v\bigg| &\lesssim  |\lambda| \eta^{\frac{2}{3}} \big( \| \zeta_\eta N_{\lambda,\eta} \|^2_{L^2(\tilde C_\eta)} +\|M\|^2_{L^2({\tilde C_\eta})} \big)+ \int_{\tilde C_\eta} \big|\eta v_1 M \bar N_{\lambda,\eta} \zeta_\eta^2\big| \mathrm{d}v   \nonumber \\
&+  |b(\lambda,\eta)| \int_{\tilde C_\eta} \big|\Phi \bar N_{\lambda,\eta} \zeta_\eta^2\big| \mathrm{d}v . \label{Re int_C Q(N_eta) bar N_eta}
\end{align}
Note that: $\ds \Re \int_{\tilde C_\eta} \mathrm{i} \eta v_1 | N_{\lambda,\eta} \zeta_\eta|^2\mathrm{d}v = 0$.  Since inequality \eqref{Phi zeta_eta} remains true on $\tilde C_\eta$: $|v'|\geqslant s_0 \eta^{-\frac{1}{3}}$, we have 
\begin{align} 
 \int_{\tilde C_\eta} \big|\Phi \bar N_{\lambda,\eta} \zeta_\eta^2\big| \mathrm{d}v &\lesssim s_0^{-2}\eta^{\frac{2}{3}} \|M\|_{L^2({\tilde C_\eta})} \| N_{\lambda,\eta} \zeta_\eta \|_{L^2({\tilde C_\eta})} \\
 &\lesssim s_0^{-2-\delta}\eta^{\frac{2+\delta}{3}} \big\| |v'|^\delta M \big\|_{L^2(\RR^d)} \| N_{\lambda,\eta} \zeta_\eta \|_{L^2({\tilde C_\eta})} \nonumber \\
 &\lesssim s_0^{-2-\delta}\eta^{\frac{2+\delta}{3}} \bigg( \| \zeta_\eta N_{\lambda,\eta} \|^2_{L^2(\tilde C_\eta)} +\big\| |v|^\delta M \big\|^2_{L^2(\RR^d)} \bigg) .\label{int_C Phi N_eta zeta^2}
\end{align}
Now, the right hand term of the first line in \eqref{Re int_C Q(N_eta) bar N_eta} is treated as follows:
\begin{align}
\int_{\tilde C_\eta} \big|\eta v_1 M \bar N_{\lambda,\eta} \zeta_\eta^2\big| \mathrm{d}v &\leqslant 2s_0^{1-\delta}\eta^{\frac{2+\delta}{3}}\| \zeta_\eta N_{\lambda,\eta} \|_{L^2(\tilde C_\eta)} \big\| |v|^\delta \zeta_\eta  M \big\|_{L^2({\tilde C_\eta})} \nonumber \\
&\leqslant s_0^{1-\delta}\eta^{\frac{2+\delta}{3}} \bigg( \| \zeta_\eta N_{\lambda,\eta} \|^2_{L^2(\tilde C_\eta)} + \big\| |v|^\delta  M \big\|^2_{L^2(\RR^d)} \bigg) .
\end{align}
Hence, from \eqref{zeta N_eta sur tilde C}, \eqref{zeta' N_eta sur tilde C} and the estimates obtained for the terms of \eqref{Re int_C Q(N_eta) bar N_eta} we obtain
\begin{align*}
\| \zeta_\eta N_{\lambda,\eta} \|^2_{L^2(\tilde C_\eta)} &\lesssim \big(s_0^2 |\lambda|+s_0^{3-\delta}\eta^{\frac{\delta}{3}}+s_0^{-\delta}\eta^\frac{\delta}{3}|b(\lambda,\eta)| \big)  \| \zeta_\eta N_{\lambda,\eta} \|^2_{L^2(\tilde C_\eta)} + \|N_{\lambda,\eta}\|^2_{L^2(A_\eta^c)} \\
&+\eta^{\frac{\delta}{3}}\big(s_0^{2-\delta} |\lambda|+s_0^{3-\delta}+s_0^{-\delta}|b(\lambda,\eta)| \big) \big\| |v|^\delta M \big\|^2_{L^2(\RR^d)} .
\end{align*}
So, for $s_0$ fixed and $|\lambda|$ and $\eta$ small enough, $\big(s_0^2 |\lambda|+s_0^{3-\delta}\eta^{\frac{\delta}{3}}+s_0^{-\delta}\eta^\frac{\delta}{3}|b(\lambda,\eta)| \big) \leqslant \frac{1}{2}$ and the term $\| \zeta_\eta N_{\lambda,\eta} \|^2_{L^2(\tilde C_\eta)}$ in the right side of the previous inequality is absorbed. Thus,
$$\| \zeta_\eta N_{\lambda,\eta} \|^2_{L^2(\tilde C_\eta)} \lesssim  \|N_{\lambda,\eta}\|^2_{L^2(A_\eta^c)} +\eta^{\frac{\delta}{3}}\big(s_0^{2-\delta} |\lambda|+s_0^{3-\delta}+s_0^{-\delta}|b(\lambda,\eta)| \big) \big\| |v|^\delta M \big\|^2_{L^2(\RR^d)}. $$
Hence  inequality \eqref{N_eta sur C} holds true.\\

\noindent \textbf{Step 3: Estimation of $\|N_{\lambda,\eta}\|_{L^2(B_\eta)}$. }  Recall that $B_\eta :=\{|v| \leqslant s_0 \eta^{-\frac{1}{3}}\}$.We claim that: 
\begin{align}\label{N_eta sur B}
\| N_{\lambda,\eta} \|^2_{L^2(B_\eta)} &\lesssim \nu_1  \| N_{\lambda,\eta} \|^2_{L^2(A_\eta)} + s_0^2 |\lambda| \| N_{\lambda,\eta} \|^2_{L^2(A_\eta^c)} + s_0^{2-\delta} \eta^{\frac{\delta}{3}} \big\| \eta^{\frac{1}{3}} v_1 M_{\lambda,\eta} \big\|^2_{L^2(A_\eta^c)} \nonumber \\
&+ ( s_0^2 |\lambda| + |c_\eta - 1| ) \|M\|^2_2 + s_0^{3-\delta} \eta^{\frac{\delta}{3}} \big\| |v_1|^\delta M \big\|^2_2 
\end{align}
where $\nu_1:=\nu_1(\lambda,\eta,s_0) = s_0^2 |\lambda| +s_0^{3-\delta} \eta^{\frac{\delta}{3}}$ and $\ds c_\eta := \bigg(\int_{\RR^d}\frac{M^2}{\langle v \rangle^2} \mathrm{d}v\bigg)^{-1}\int_{\RR^d} \frac{M M_{\lambda,\eta}}{\langle v \rangle^2}\mathrm{d}v$.\\

\noindent Let us denote $\ds \tilde N_{\lambda,\eta} := M_{\lambda,\eta}-c_\eta M$ the orthogonal projection of $M_{\lambda,\eta}$ to $M$ for the weighted scalar product $\int \frac{\cdot}{\langle v \rangle^2}$. On the one hand, we have: 
$$ \|N_{\lambda,\eta}\|^2_{L^2(B_\eta)} \lesssim \| \tilde N_{\lambda,\eta}\|^2_{L^2(B_\eta)} + |c_\eta - 1| \| M \|^2_{L^2(B_\eta)} $$ 
and $$ \| \tilde N_{\lambda,\eta}\|^2_{L^2(B_\eta)} \lesssim s_0^2 \eta^{-\frac{2}{3}} \bigg\| \frac{\tilde N_{\lambda, \eta}}{\langle v \rangle} \bigg\|^2_{L^2(\RR^d)} , $$ 
since $ \langle v \rangle \lesssim s_0 \eta^{-\frac{1}{3}}$ on $B_\eta$. On the other hand, applying  inequality \eqref{ineg.  Hardy-Poincare} to $\tilde N_{\lambda,\eta}$ which satisfies condition \eqref{condition d'orthogonalite}, we obtain:
$$
\int_{\RR^d}  \frac{|\tilde N_{\lambda,\eta}|^2}{\langle v \rangle^2} \mathrm{d}v \leqslant C_{\gamma,d} \int_{\RR^d} \bigg| \nabla_v \bigg(\frac{N_{\lambda,\eta}}{M}\bigg)\bigg|^2 M^2 \mathrm{d}v .
$$
Therefore,
\begin{equation}
 \|N_{\lambda,\eta}\|^2_{L^2(B_\eta)}  \lesssim s_0^2 \eta^{-\frac{2}{3}} \bigg\| \nabla_v \bigg(\frac{N_{\lambda,\eta}}{M}\bigg)M \bigg\|^2_{L^2(\RR^d)} + |c_\eta-1| \|M\|^2_{L^2(\RR^d)} .
\end{equation}
We have moreover,
$$ \bigg\| \nabla_v \bigg(\frac{N_{\lambda,\eta}}{M}\bigg)M \bigg\|^2_{L^2(\RR^d)} =  \int_{\RR^d} Q(N_{\lambda,\eta})\bar N_{\lambda,\eta} \mathrm{d}v .$$
Then, by integrating the equation of $N_{\lambda,\eta}$ multiplied by $\bar N_{\lambda,\eta}$, we obtain
$$ \int_{\RR^d} Q(N_{\lambda,\eta})\bar N_{\lambda,\eta} \mathrm{d}v + |\langle N_{\lambda,\eta},\Phi \rangle|^2 = \Re \int_{\RR^d} \bigg( \lambda \eta^{\frac{2}{3}} \big(|N_{\lambda,\eta}|^2  + M \bar N_{\lambda,\eta} \big) -  \mathrm{i} \eta v_1 M \bar N_{\lambda,\eta} \bigg) \mathrm{d}v . $$
From where,
\begin{equation}\label{M grad N_eta/M}
\bigg\| \nabla_v \bigg(\frac{N_{\lambda,\eta}}{M}\bigg)M \bigg\|^2_{L^2(\RR^d)} \lesssim |\lambda| \eta^{\frac{2}{3}} \big( \| N_{\lambda,\eta} \|^2_2 + \|M\|^2_2 \big)  + \eta \int_{\RR^d}  |v_1 M \mathrm{Im} N_{\lambda,\eta} | \mathrm{d}v .
\end{equation}
For the last term, using the fact that $\ \mathrm{Im} N_{\lambda,\eta} = \mathrm{Im} M_{\lambda,\eta}$, we write
\begin{equation}\label{eta^(1/3)v_1 M_eta M 1}
 \eta^{\frac{1}{3}} \int_{\RR^d}  |v_1 M \mathrm{Im} N_{\lambda,\eta} | \mathrm{d}v = \int_{A_\eta}  \eta^{\frac{1}{3}} |v_1|^{1-\delta} |v_1|^\delta M |\mathrm{Im} N_{\lambda,\eta} | \mathrm{d}v + \int_{A_\eta^c}  \eta^{\frac{1}{3}} |v_1  \mathrm{Im} M_{\lambda,\eta} | M \mathrm{d}v ,
\end{equation}
and since on $A_\eta:=\{|v_1| \leqslant s_0\eta^{-\frac{1}{3}}\}$ we have: $|v_1|^{1-\delta}\leqslant s_0^{1-\delta}\eta^{\frac{\delta-1}{3}}$, and on $A_\eta^c:=\{|v_1| \geqslant s_0\eta^{-\frac{1}{3}}\}$, $M(v) \leqslant s_0^{-\delta}\eta^{\frac{\delta}{3}} |v_1|^{\delta}M(v)$ then,   we obtain
\begin{align}
\eta^{\frac{1}{3}} \int_{\RR^d}  |v_1 M \mathrm{Im} N_{\lambda,\eta} | \mathrm{d}v &\leqslant \frac{\eta^{\frac{\delta}{3}}}{s_0^\delta} \big\| |v_1|^\delta M \big\|_2 \bigg( s_0 \|N_{\lambda,\eta}\|_{L^2(A_\eta)}  + \big\| \eta^{\frac{1}{3}} v_1 M_{\lambda,\eta} \big\|_{L^2(A_\eta^c)} \bigg)  \nonumber  \\
&\leqslant \frac{1}{2} \frac{\eta^{\frac{\delta}{3}}}{s_0^\delta} \bigg( s_0\|N_{\lambda,\eta}\|^2_{L^2(A_\eta)} + \big\| \eta^{\frac{1}{3}} v_1 M_{\lambda,\eta} \big\|^2_{L^2(A_\eta^c)} +2s_0 \big\| |v_1|^\delta M \big\|^2_2  \bigg) . \label{eta^(1/3)v_1 M_eta M 2}
\end{align}
Thus, returning to \eqref{M grad N_eta/M} we get
\begin{align*}
 \bigg\| \nabla_v \bigg(\frac{N_{\lambda,\eta}}{M}\bigg)M \bigg\|^2_{L^2(\RR^d)} &\lesssim \eta^{\frac{2}{3}} (|\lambda| +s_0^{1-\delta} \eta^{\frac{\delta}{3}})  \| N_{\lambda,\eta} \|^2_{L^2(A_\eta)} + \eta^{\frac{2}{3}} |\lambda| \| N_{\lambda,\eta} \|^2_{L^2(A_\eta^c)} \\
 &+ s_0^{-\delta} \eta^{\frac{2+\delta}{3}} \big\| \eta^{\frac{1}{3}} v_1 M_{\lambda,\eta} \big\|^2_{L^2(A_\eta^c)} +  \eta^{\frac{2}{3}} |\lambda| \|M\|^2_2 + s_0^{1-\delta} \eta^{\frac{2+\delta}{3}} \big\| |v_1|^\delta M \big\|^2_{L^2(\RR^d)} .
\end{align*}
Hence inequality \eqref{N_eta sur B} holds by multiplying the previous one by $s_0^2\eta^{-\frac{2}{3}}$ and adding the term $|c_\eta-1| \|M\|^2_{L^2(\RR^d)}$.     \\

\noindent \textbf{Step 4: Conclusion.} In this step, we will combine all the estimates obtained in the previous steps in order to conclude. First, by summing  inequalities \eqref{N_eta sur C} and \eqref{N_eta sur B} obtained in steps 2 and 3 respectively, and since $A_\eta = B_\eta \cup C_\eta$, we obtain
\begin{align}\label{N_eta sur A_eta}
\|N_{\lambda,\eta}\|^2_{L^2(A_\eta)} &\lesssim \nu_1 \|N_{\lambda,\eta}\|^2_{L^2(A_\eta)} + (s_0^2|\lambda|+1)\|N_{\lambda,\eta}\|^2_{L^2(A_\eta^c)} +  \frac{\eta^{\frac{\delta}{3}}}{s_0^\delta} \big\| \eta^{\frac{1}{3}} v_1 M_{\lambda,\eta} \big\|^2_{L^2(A_\eta^c)} \nonumber \\
&+ ( s_0^2|\lambda| +|c_\eta-1| ) \|M\|_{L^2(\RR^d)}^2  +  c_2^\eta ,
\end{align}
where $\nu_1:=s_0^2|\lambda|+s_0^{3-\delta}\eta^{\frac{\delta}{3}}$ and $c_2^\eta := s_0^{-\delta}\eta^{\frac{\delta}{3}}\big(s_0^3+s_0^2 |\lambda| +  |b(\lambda,\eta)|\big)\big\| |v|^\delta M \big\|^2_{L^2(\RR^d)}$. Now, since 
$$\|N_{\lambda,\eta}\|^2_{L^2(A_\eta^c)} \lesssim \|M_{\lambda,\eta}\|^2_{L^2(A_\eta^c)} + s_0^{-2\delta} \eta^{\frac{2\delta}{3}} \big\| |v_1|^\delta M \big\|^2_{L^2(A_\eta^c)} ,$$
then, using inequality \eqref{M_eta sur A^c} for the two terms $\|M_{\lambda,\eta}\|^2_{L^2(A_\eta^c)}$ (in the previous inequality) and $\big\| \eta^{\frac{1}{3}} v_1 M_{\lambda,\eta} \big\|^2_{L^2(A_\eta^c)}$ (in \eqref{N_eta sur A_eta}), returning to inequality \eqref{N_eta sur A_eta} we obtain
\begin{align*}
\|N_{\lambda,\eta}\|^2_{L^2(A_\eta)} &\lesssim \big(\nu_1+\frac{1}{s_0^3}\big)\|N_{\lambda,\eta}\|^2_{L^2(A_\eta)} + (s_0^2|\lambda| +|c_\eta-1| ) \|M\|_{L^2(\RR^d)}^2 + c_2^\eta.
\end{align*} 
Therefore, we first set $s_0$ large enough so that $\frac{1}{s_0^3} \leqslant \frac{1}{4}$, then for $|\lambda|$ and $\eta$ small enough so that $\nu_1 := s_0^2|\lambda|+s_0^{3-\delta}\eta^{\frac{\delta}{3}} \leqslant \frac{1}{4}$, we get:
\begin{equation}\label{N_eta sur A}
\|N_{\lambda,\eta}\|^2_{L^2(A_\eta)} \lesssim (s_0^2|\lambda| + |c_\eta-1| ) \|M\|_{L^2(\RR^d)}^2 + c_2^\eta \lesssim 1. 
\end{equation}
The right-hand side of the inequality above is uniformly bounded since $s_0^2|\lambda| \leqslant \frac{1}{4}$, $|c_\eta-1| \to 0$ and $c_2^\eta \to 0$ when $\eta$ goes $0$. Indeed, we have
\begin{align}\label{c_eta-1}
| c_\eta - 1 | &=  \bigg(\int_{\RR^d}\frac{M^2}{\langle v \rangle^2}\mathrm{d}v\bigg)^{-1}\bigg|\int_{\RR^d} \frac{M(M_{\lambda,\eta}-M)}{\langle v \rangle^2}\mathrm{d}v\bigg| \nonumber \\
&\leqslant \bigg\| \frac{M}{\langle v \rangle} \bigg\|^{-1}_2 \bigg\| \frac{N_{\lambda,\eta}}{\langle v \rangle} \bigg\|_2 \leqslant C \|N_{\lambda,\eta}\|_{\mathcal{H}_0} \underset{\eta \to 0}{\longrightarrow} 0.
\end{align}
For $c^\eta_2$, we have $c^\eta_2 \lesssim \eta^{\frac{\delta}{3}}$ since $|v|^\delta M \in L^2$ for all $\gamma > \frac{d}{2}$ and since $|b(\lambda,\eta)| \lesssim 1$ thanks to the second point of Remark \ref{rmq sur thm d'existence}.\\ 

\ni Now, we resume all the assumptions we did on $s_0$, $\lambda$ et $\eta$: 
$$ \frac{C_1}{s_0}\bigg(|\lambda| + \frac{1+|b(\lambda,\eta)|}{s_0^2}\bigg) \leqslant \frac{1}{2}, \quad \frac{C_2}{s_0^3} \leqslant \frac{1}{4}, \quad  C_3\big(s_0^2|\lambda|+s_0^{3-\delta}\eta^{\frac{\delta}{3}}\big) \leqslant \frac{1}{4} 
$$
and $$ C_4\big(s_0^2|\lambda|+s_0^{3-\delta}\eta^{\frac{\delta}{3}} +s_0^{-\delta}\eta^{\frac{\delta}{3}}|b(\lambda,\eta)|\big) \leqslant \frac{1}{2} . $$
Recall that $\delta:=\frac{1}{2}(\gamma-\frac{d}{2})>0$ for all $\gamma>\frac{d}{2}$. So, if we start by setting $s_0$ large enough, then $\lambda$ small enough, then $\eta$ small enough, we recover all the previous inequalities. \\

\noindent Finally, by injecting  inequality \eqref{N_eta sur A_eta} into \eqref{M_eta sur A^c}, we obtain:
\begin{align}\label{M_eta sur A^c fin}
\|M_{\lambda,\eta}\|^2_{L^2(A_\eta^c)} \leqslant \frac{1}{s_0^2} \big\| \eta^{\frac{1}{3}} v_1 M_{\lambda,\eta} \big\|_{L^2(A_\eta^c)} \lesssim \frac{1}{s_0^3}\big( \nu_1 + c_2^\eta \big) \lesssim 1.
\end{align}
Hence, $N_{\lambda,\eta}$ as well as $M_{\lambda,\eta}$ are uniformly bounded in $L^2(\RR^d)$.   
Now, from \eqref{N_eta sur A_eta} and \eqref{M_eta sur A^c fin} we obtain: 
$$ \|N_{\lambda,\eta}\|_{L^2(\RR^d)}^2 \lesssim  |\lambda| + |c_\eta-1| + c_2^\eta .$$
Hence the inequality \eqref{estimation de N_lambda,eta dans L^2} holds with $\nu_\eta := |c_\eta-1| + c_2^\eta \underset{\eta \to 0}{\longrightarrow} 0$. 
\ep
\subsection{Study of the constraint}
In this subsection, we will show the existence of a $\mu$, a function of $\eta$, such that the constraint $\langle M_{\mu(\eta),\eta}-M,\Phi\rangle =0$ is satisfied.
Let us start by giving the following result, which is a corollary of Proposition \ref{propostion estimation L^2}.
\begin{corollary}\label{lim int M_eta ...}
Let $M_{\lambda,\eta}$ be the solution to equation \eqref{eq de M_eta + cont.}. Then, for all $\lambda \in \CC$ such that, $|\lambda| \leqslant \lambda_0$ with $\lambda_0$ small enough, the following limit holds:
\begin{equation}\label{lim int eta^(1/3) v_1 M_eta M}
\underset{\eta \to 0}{\lim } \int_{\RR^d} \eta^{\frac{1}{3}} v_1 M_{\lambda,\eta}(v)M(v) \mathrm{d}v = 0 .
\end{equation}
For $\lambda =0$, one has
\begin{equation}\label{lim int  M_eta M}
\underset{\eta \to 0}{\lim } \int_{\RR^d} M_{0,\eta}(v)M(v) \mathrm{d}v = \int_{\RR^d} M^2(v) \mathrm{d}v .
\end{equation}
\end{corollary}
\bp For the first point, we proceed exactly as in \eqref{eta^(1/3)v_1 M_eta M 1}, i.e. cutting the integral into two parts $A_\eta:= \{|v_1| \leqslant s_0\eta^{-\frac{1}{3}}\}$ and $A_\eta^c$, we write: 
\begin{align*}
\bigg| \int_{\RR^d} \eta^{\frac{1}{3}} v_1 M_{\lambda,\eta}(v)M(v) \mathrm{d}v \bigg| &\leqslant \eta^{\frac{1}{3}} \int_{A_\eta}  |v_1|^{1-\delta}  |v_1|^\delta M(v) |M_{\lambda,\eta}(v)| \mathrm{d}v   \\
& \hspace{2cm}+ \int_{A_\eta^c} |v_1|^{-\delta} |v_1|^\delta M(v) |\eta^{\frac{1}{3}} v_1 M_{\lambda,\eta}(v)| \mathrm{d}v \\
&\leqslant   s_0^{-\delta}\eta^{\frac{\delta}{3}}\big\| |v_1|^\delta M \big\|_2\bigg( s_0 \| M_{\lambda,\eta} \|_{L^2(A_\eta)} + \big\| \eta^{\frac{1}{3}} v_1 M_{\lambda,\eta} \big\|_{L^2(A_\eta^c)}\bigg) \\
&\lesssim \eta^{\frac{\delta}{3}} \underset{\eta \to 0}{\longrightarrow}0 , \quad \mbox{ thanks to }\eqref{N_eta sur A} \mbox{ and } \eqref{M_eta sur A^c fin}.
\end{align*}
For the second point, for $\lambda =0$, we write 
$$
 \bigg| \int_{\RR^d} \big[M_{0,\eta}(v)-M(v)\big]M(v) \mathrm{d}v \bigg| \leqslant \| N_{0,\eta} \|_{L^2(\RR^d)} \|M\|_{L^2(\RR^d)} ,
$$
and the limit \eqref{lim int  M_eta M} holds true thanks to  inequality \eqref{estimation de N_lambda,eta dans L^2} of Proposition \ref{propostion estimation L^2}.
\ep

\ni \begin{proposition}[Constraint]\label{contrainte}
Define $$ B(\lambda,\eta):= \eta^{-\frac{2}{3}} b(\lambda,\eta)  .   $$
\begin{enumerate}
\item The expression of $B(\lambda,\eta)$ is given by
\begin{equation}\label{b(lambda,eta)}
B(\lambda,\eta)= \eta^{-\frac{2}{3}}\langle N_{\lambda,\eta},\Phi\rangle = \int_{\RR^d}(\lambda- \mathrm{i} \eta^{\frac{1}{3}} v)M_{\lambda,\eta}(v)M(v)\mathrm{d}v .
\end{equation}
\item The $\eta$ order of $B(\lambda,\eta)$ in its expansion with respect to $\lambda$ is given by
\begin{equation}\label{B(lambda,0)}
\underset{\eta \rightarrow 0}{\lim} \ \frac{\pa B}{\pa \lambda}(0,\eta) = \int_{\RR^d} M^2(v)\mathrm{d}v .
\end{equation}
\item There exists $\tilde\eta_0, \tilde\lambda_0>0$ small enough, a function $\tilde{\lambda}: \{|\eta|\leqslant\tilde\eta_0\} \longrightarrow \{|\lambda|\leqslant\tilde\lambda_0\}$ such that,\\
for all $(\lambda,\eta)\in [0,\tilde\eta_0[\times\{|\lambda|<\tilde\lambda_0\}$,  $\lambda = \tilde{\lambda}(\eta)$ and the constraint is satisfied: $$B(\lambda,\eta)= B(\tilde{\lambda}(\eta),\eta)=0.$$
\end{enumerate}
Consequently, $\mu(\eta) = \eta^{\frac{2}{3}} \tilde{\lambda} (\eta)$ is the eigenvalue associated to the eigenfunction $M_\eta: = M_{\tilde{\lambda}(\eta),\eta}$ for the operator $\mathcal{L}_\eta$, and the couple $\big(\mu(\eta),M_\eta \big)$ is solution to the spectral problem \eqref{M_mu,eta}.
\end{proposition} 

\noindent \bp 1. The first point is obtained by integrating the equation of $M_{\lambda,\eta}$ multiplied by $M$, and using the assumption $\langle M,\Phi \rangle = 1$.  \\
2. This point is exactly limit \eqref{lim int  M_eta M} of Corollary \ref{lim int M_eta ...}.  \\
3. The third point follows from the implicit function theorem applied to the function $B$ around the point $(\lambda,\eta)=(0,0)$. 
\ep

\subsection{Approximation of the eigenvalue}
In this subsection, we will give an approximation for the eigenvalue $\mu(\eta)$ given in Proposition \ref{contrainte}. The study of this limit is based on some estimates on $M_{0,\eta}$, the solution of equation \eqref{eq penalisee2} for $\lambda=0$, as well as the solution of the rescaled equation.\\

Before giving the proposition which summarizes the essential points of this subsection, we will first start by introducing the rescaled function of $M_{0,\eta}$ as well as the equation satisfied by this function. Recall that $M_{0,\eta}$ satisfies the equation:
$$ \big[ Q + \mathrm{i} \eta v_1 \big] M_{0,\eta}(v) = - b(0,\eta) \Phi(v) ,  \quad v \in \RR^d ,$$
with $Q = -\frac{1}{M}\nabla_v \cdot \big(M^2\nabla_v\big(\frac{\cdot}{M}\big)\big)$ and $b(0,\eta) = \langle M_{0,\eta}-M,\Phi \rangle$. 
Then, the rescaled function $H_\eta$ defined by $H_\eta(s):=\eta^{-\frac{\gamma}{3}}M_{0,\eta}(\eta^{-\frac{1}{3}}s)$ is solution to the rescaled equation
\begin{equation}\label{eq H_eta}
\big[ Q_\eta + \mathrm{i} s_1 \big] H_\eta(s) = - \eta^{-\frac{\gamma+2}{3}} b(0,\eta) \Phi_\eta(s) , \quad s \in \RR^d ,
\end{equation}
where 
$$ Q_\eta := -\frac{1}{|s|_\eta^{-\gamma}}\nabla_s \cdot \bigg(|s|_\eta^{-2\gamma}\nabla_s\bigg(\frac{\cdot}{|s|_\eta^{-\gamma}}\bigg)\bigg) \ ,  \quad  |s|_\eta^{-\gamma} :=\eta^{-\frac{\gamma}{3}}M(\eta^{-\frac{1}{3}}s)= \big(\eta^{\frac{2}{3}}+|s|^2\big)^{-\frac{\gamma}{2}} $$
and 
\begin{equation}\label{Phi_eta < ...}
\Phi_\eta(s) :=\Phi(\eta^{-\frac{1}{3}}s) = c_{\gamma,d}\ \eta^{\frac{\gamma+2}{3}} |s|_\eta^{-\gamma-2} .
\end{equation}
Note that: $Q(M)=0$ implies that $Q_\eta(|s|_\eta^{-\gamma})=0$.
\begin{proposition}[Approximation of the eigenvalue]\label{vp} 
Let $\alpha:=\frac{2\gamma-d+2}{3}$ for all $\gamma \in (\frac{d}{2},\frac{d+4}{2})$.
The eigenvalue $\mu(\eta)$ satisfies
\begin{equation}\label{mu(eta)}
\mu(\eta)=\overline \mu (-\eta)= \kappa |\eta|^{\alpha}\big(1+O(|\eta|^\alpha)\big) ,
\end{equation}
where  $\kappa$ is a positive constant given by
\begin{equation}\label{kappa}
\kappa := - 2 C_\beta^{2} \int_{\{s_1>0\}} s_1|s|^{-\gamma} \mathrm{Im} H_0(s) \mathrm{d}s ,
\end{equation}
and where $H_0$ is the unique solution to
\begin{equation}\label{eq de H_0}
\bigg[-\Delta_s + \frac{\gamma(\gamma-d+2)}{|s|^2} + \mathrm{i} s_1 \bigg]H_0(s)=0, \ s\in \RR^d\setminus\{0\}  ,
\end{equation}
satisfying \begin{equation}\label{condition H_0}
\int_{\{|s_1|\geqslant 1\}}|H_0(s)|^2\mathrm{d}s <+\infty \  \mbox{ and } \  H_0(s)\underset{0}{\sim} |s|^{-\gamma}  .
\end{equation}
\end{proposition}
\begin{remark}
Note that the existence of solutions for equation \eqref{eq de H_0} is obtained by passing to the limit in the rescaled equation \eqref{eq H_eta}, while the uniqueness is obtained by an integration by part on $\RR^d\setminus \{0\}$,  using the two conditions of \eqref{condition H_0}.
\end{remark}

\noindent In order to get  Proposition \ref{vp}, we need to prove the following series of lemmas. \\

\noindent The first one show that the small velocities in the first direction do not participate in the limit of the diffusion coefficient.
\begin{lemma}[Small velocities]\label{petites vitesses}
\item \begin{enumerate}
\item For all $\gamma \in(\frac{d+1}{2},\frac{d+4}{2})$, one has
\begin{equation}\label{Im M_0,eta/<v> v_1 petit}
\int_{\{|v_1|\leqslant R\}}  \bigg|\frac{\mathrm{Im} M_{0,\eta}(v)}{\langle v \rangle}\bigg|^2 \mathrm{d}v \lesssim \eta .
\end{equation}
\item For all $\gamma \in (\frac{d}{2},\frac{d+4}{2})$
\begin{equation}\label{lim int Im M_0,eta/<v> v_1 petit}
\underset{\eta \to 0}{\lim} \ \eta^{1-\alpha}  \int_{\{|v_1|\leqslant R\}} v_1 M_{0,\eta}(v)M(v) \mathrm{d}v = 0 .
\end{equation}
\end{enumerate}
\end{lemma}

\noindent  The second one contains some important estimates on the rescaled solution for large velocities.
\begin{lemma}[Large velocities]\label{grandes vitesses} Let $s_0>0$ be fixed, large enough.  We have the following estimates, uniform with respect to $\eta$,  for the rescaled solution:
\begin{enumerate}
\item For all $\gamma \in (\frac{d}{2},\frac{d+1}{2})$, one has
\begin{equation}\label{estimation de  H_eta pour gamma<(d+1)/2}
\big\| |s_1|^{\frac{1}{2}} \mathrm{Im} H_\eta \big\|_{L^2(\{|s_1|\leqslant s_0\})} + \| s_1 \mathrm{Im} H_\eta \|_{L^2(\{|s_1|\geqslant s_0\})} \lesssim 1 .
\end{equation}
\item For all $\gamma \in (\frac{d+1}{2},\frac{d+4}{2})$, one has
\begin{equation}\label{estimation de H_eta pour gamma>(d+1)/2}
\bigg\|  \frac{\mathrm{Im} H_\eta}{|s|_\eta} \bigg\|_{L^2(\{|s_1|\leqslant s_0\})} + \| s_1 \mathrm{Im} H_\eta \|_{L^2(\{|s_1|\geqslant s_0\})} \lesssim 1 .
\end{equation}
\end{enumerate}
\end{lemma}

\noindent \bpl \ref{petites vitesses}. \textbf{1.} By Remark \ref{rmq sur thm d'existence}, since $M_{0,\eta}$ is symmetric with respect to $v_1$ in the following sense: $\bar M_{\bar \lambda,\eta}(-v_1,v')=M_{\lambda,\eta}(v_1,v')$ then,  $\mathrm{Im} M_{0,\eta}$ is odd with respect to $v_1$. Therefore,
\begin{align*}
\int_{\RR^d} \frac{M\mathrm{Im} M_{0,\eta} }{\langle v \rangle^2}\mathrm{d}v &=  \int_{\RR^{d-1}} \bigg[ \int_{-\infty}^0 M \mathrm{Im} M_{0,\eta}(v_1,v') \frac{\mathrm{d}v_1}{\langle v \rangle^2}+\int_0^\infty M \mathrm{Im} M_{0,\eta}(v_1,v') \frac{\mathrm{d}v_1}{\langle v \rangle^2}\bigg] \mathrm{d}v' \\
&=  \int_{\RR^{d-1}} \int_0^\infty \bigg[ \mathrm{Im} \bar M_{0,\eta}(v_1,v') + \mathrm{Im} M_{0,\eta}(v_1,v') \bigg]M \frac{\mathrm{d}v_1}{\langle v \rangle^2} \mathrm{d}v' \\
&= 0.
\end{align*}
Note that we used the \emph{symmetry} of $M$ in the previous equalities.  Thus, the function $\mathrm{Im} M_{0,\eta}$ satisfies condition \eqref{condition d'orthogonalite}. Then, by the Hardy-Poincar\'e inequality \eqref{ineg.  Hardy-Poincare},  there exists a positive constant $C_{\gamma,d}$ such that: $$ \int_{\{|v_1|\leqslant R\}}  \bigg|\frac{\mathrm{Im} M_{0,\eta}(v)}{\langle v \rangle}\bigg|^2 \mathrm{d}v \leqslant \bigg\|\frac{\mathrm{Im} M_{0,\eta}}{\langle v \rangle}\bigg\|^2_2 \leqslant C_{\gamma,d} \bigg\|\nabla_v\bigg(\frac{\mathrm{Im} M_{0,\eta}}{M}\bigg)M\bigg\|^2_{L^2(\RR^d)} .$$
Now, as in step 3 of the proof of Proposition \ref{propostion estimation L^2}, we have on the one hand,
$$ 
\bigg\|\nabla_v\bigg(\frac{\mathrm{Im} M_{0,\eta}}{M}\bigg)M\bigg\|^2_2 = \int_{\RR^d} Q(\mathrm{Im} M_{0,\eta})\mathrm{Im} M_{0,\eta} \mathrm{d}v .
$$
On the other hand,
$$  Q(\mathrm{Im} M_{0,\eta}) = \eta v_1 \Re M_{0,\eta} - \eta \bigg(\int_{\RR^d} v_1 M \Re M_{0,\eta} \mathrm{d}v\bigg) \Phi . $$ 
Which implies that, 
$$
 \bigg| \int_{\RR^d} Q(\mathrm{Im} M_{0,\eta})\mathrm{Im} M_{0,\eta} \mathrm{d}v  \bigg|  \leqslant \eta \big\| |v_1|^\frac{1}{2} M_{0,\eta} \big\|_2\bigg( \big\| |v_1|^\frac{1}{2} M_{0,\eta} \big\|_2 + \big\| |v_1|^\frac{1}{2} M \big\|_2 \| M_{0,\eta} \|_2 \| \Phi \|_2 \bigg) .
$$
Hence inequality \eqref{Im M_0,eta/<v> v_1 petit} holds thanks to  $\big(1+ |v_1|^\frac{1}{2}\big) M_{0,\eta} \in L^2(\RR^d)$ for $\gamma > \frac{d+1}{2}$ (Proposition \ref{propostion estimation L^2}).\\
\textbf{2.}  First, since $v_1 \mathrm{Im} M_{0,\eta}$ and $M$ are even functions with respect to $v_1$, then
\begin{equation}\label{int_(|v_1|<R) .. .= 2 Im int }
  \int_{\{|v_1|\leqslant R\}} v_1 M_{0,\eta}(v)M(v) \mathrm{d}v = 2  \int_0^R \int_{\RR^{d-1}} v_1 \mathrm{Im} M_{0,\eta}(v)M(v) \mathrm{d}v'\mathrm{d}v_1 .
\end{equation}
\textbf{Case 1: $\gamma \in ]\frac{d}{2},\frac{d+1}{2}]$.} We have by Cauchy-Schwarz,
$$
 \eta^{1-\alpha}\bigg| \int_0^R \int_{\RR^{d-1}} v_1 \mathrm{Im} M_{0,\eta}(v)M(v) \mathrm{d}v \bigg| \leqslant  R \eta^{1-\alpha} \| \mathrm{Im} M_{0,\eta} \|_2 \| M \|_2 \leqslant R \eta^{1-\alpha} \| N_{0,\eta} \|_2 \| M \|_2 \underset{\eta \to 0}{\longrightarrow} 0,
$$
since $1-\alpha = \frac{1+d-2\gamma}{3} \geqslant 0$ for all $\gamma \leqslant\frac{d+1}{2}$ and $\|\mathrm{Im} M_{0,\eta}\|_2 = \|\mathrm{Im} N_{0,\eta}\|_2 \leqslant \|N_{0,\eta}\|_2 \underset{\eta \to 0}{\longrightarrow} 0$.  \\
\textbf{Case 2: $\gamma \in ]\frac{d+1}{2},\frac{d+4}{2}[$.} Similary, we have by Cauchy-Schwarz,
\begin{align*}
\eta^{1-\alpha} \bigg| \int_{\{|v_1|\leqslant R\}} v_1 M_{0,\eta}(v)M(v) \mathrm{d}v \bigg|  & \leqslant \eta^{1-\alpha}  \big\| v_1 \langle v \rangle M \big\|_{L^2(\{|v_1|\leqslant R\})} \bigg\| \frac{\mathrm{Im} M_{0,\eta}}{\langle v \rangle } \bigg\|_{L^2(\{|v_1|\leqslant R\})} \\
&\lesssim \eta^{2-\alpha} \underset{\eta \to 0}{\longrightarrow} 0,
\end{align*}
thanks to  inequality \eqref{Im M_0,eta/<v> v_1 petit} and since $\alpha < 2$ for $\gamma < \frac{d+4}{2}$ and $v_1 \langle v \rangle M \in L^2(\{|v_1|\leqslant R\})$ for $\gamma<\frac{d+4}{2}$.
\ep

\noindent \bpl \ref{grandes vitesses}. We will establish estimates on different ranges of (rescalated) velocities, and in order to avoid long expressions in the proof, we will fix some notations of ``sets" as in the proof of Proposition \ref{propostion estimation L^2}. Let denote $s:=(s_1,s') \in \RR\times\RR^{d-1}$. Let $s_0>0$. We set: $A :=\{|s_1|\leqslant s_0\}$ (resp., $\tilde A:=\{|s_1| \leqslant 2s_0\}$), $B :=\{|s| \leqslant s_0\}$, $C:=\{|s_1| \leqslant s_0 \leqslant |s'|\}$ (resp., $\tilde C:=\{|s_1|\leqslant 2s_0\leqslant 2|s'|\}$) and finally $D :=\{|s_1|\geqslant \frac{s_0}{2}\}$. Also, for $\eta>0$, we denote by $\tilde K_\eta$ the function defined by $\tilde K_\eta := H_\eta - c_\eta |s|_\eta^{-\gamma}$, with $c_\eta$ given by 
$$
 c_\eta :=  \bigg(\int_{\RR^d} |s|_\eta^{-2\gamma-2} \mathrm{d}s \bigg)^{-1}\int_{\RR^d}   \frac{|s|_\eta^{-\gamma} H_\eta(s)}{|s|_\eta^2} \mathrm{d}s , \quad \forall \eta>0.
$$
Note that $$ \int_{\RR^d} \frac{|s|_\eta^{-\gamma} \tilde K_\eta(s)}{|s|_\eta^2} \mathrm{d}s =0 .$$ Thus, $\tilde K_\eta$ satisfies the orthogonality condition \eqref{condition d'orthogonalite} of the Hardy-Poincar\'e Lemma \ref{Hardy-Poincare}.
\begin{remark}\label{remarque c_eta}
\item \begin{enumerate}
\item Observe that 
$$ c_\eta := \bigg(\int_{\RR^d}\frac{M^2}{\langle v \rangle^2}\mathrm{d}v \bigg)^{-1}\int_{\RR^d} \frac{M(v)M_{0,\eta}(v)}{\langle v \rangle^2} \mathrm{d}v .$$
Hence, $\ \ds c_\eta \underset{\eta \to 0}{\longrightarrow} 1$ by \eqref{c_eta-1}.
\item Since $\bar M_{0,\eta}(-v_1,v')=M_{0,\eta}(v_1,v')$ for all $(v_1,v')\in\RR\times\RR^{d-1}$, then
$$ \mathrm{Im} H_\eta(-s_1,s')=-\mathrm{Im} H_\eta(s_1,s') , \quad \forall (s_1,s')\in \RR\times\RR^{d-1} . $$ Therefore, $\mathrm{Im} \ c_\eta = 0$. 
\end{enumerate}
\end{remark} 

\noindent \textbf{1.} Let $\gamma \in (\frac{d}{2},\frac{d+1}{2})$. To prove the first point of this lemma, we will proceed exactly as in the proof of Proposition \ref{propostion estimation L^2}. We estimate $\big\| |s_1|^{\frac{1}{2}} K_\eta \big\|_{L^2(B)}$ using the Hardy-Poincar\'e inequality, $\big\| |s_1|^{\frac{1}{2}} K_\eta \big\|_{L^2(C)}$ using the weighted Poincar\'e inequality, Lemma \ref{lemme Poincare2}, and estimate $\|s_1 H_\eta\|_{L^2(A^c)}$ using the equation of $H_\eta$. Thus, we obtain  inequality \eqref{estimation de  H_eta pour gamma<(d+1)/2} by combining these estimates and since $|s_1|^{\frac{1}{2}} |s|_\eta^{-\gamma} \leqslant |s_1|^{\frac{1}{2}} |s|^{-\gamma} \in L^2(A)$ for $\gamma < \frac{d+1}{2}$.\\

\noindent \textbf{Estimation of $\big\| |s_1|^{\frac{1}{2}} K_\eta \big\|_{L^2(B)}$.} Recall that $B :=\{|s| \leqslant s_0\}$. On the one hand, we have: 
$$\big\| |s_1|^{\frac{1}{2}} K_\eta \big\|^2_{L^2(B)} \lesssim \big\| |s_1|^{\frac{1}{2}} \tilde K_\eta \big\|^2_{L^2(B)} + |c_\eta-1| \big\| |s_1|^{\frac{1}{2}} |s|_\eta^{-\gamma} \big\|^2_{L^2(B)} , $$ and by the Hardy-Poincar\'e inequality \eqref{ineg.  Hardy-Poincare} we get
$$
 \big\| |s_1|^{\frac{1}{2}} \tilde K_\eta \big\|^2_{L^2(B)} \leqslant s_0^3 \bigg\| \frac{\tilde K_\eta}{|s|_\eta} \bigg\|^2_{L^2(B)}  \lesssim_{\gamma,d} s_0^3 \bigg\| \nabla_s \bigg(\frac{\tilde K_\eta}{|s|_\eta^{-\gamma}}\bigg)|s|_\eta^{-\gamma} \bigg\|^2_2 .
$$
Therefore,
\begin{equation}\label{s_1^(1/2)K_eta sur B 1}
 \big\| |s_1|^{\frac{1}{2}} K_\eta \big\|^2_{L^2(B)}  \lesssim s_0^3 \bigg\| \nabla_s \bigg(\frac{\tilde K_\eta}{|s|_\eta^{-\gamma}}\bigg)|s|_\eta^{-\gamma} \bigg\|^2_2 + |c_\eta-1| \big\| |s_1|^{\frac{1}{2}} |s|_\eta^{-\gamma} \big\|^2_{L^2(B)} .
\end{equation}
On the other hand,  since $\nabla_s \big(\frac{\tilde K_\eta}{|s|_\eta^{-\gamma}}\big) = \nabla_s \big(\frac{ K_\eta}{|s|_\eta^{-\gamma}}\big)$, then
$$  \bigg\| \nabla_s \bigg(\frac{\tilde K_\eta}{|s|_\eta^{-\gamma}}\bigg)|s|_\eta^{-\gamma} \bigg\|^2_2 =  \Re \int_{\RR^d} Q_\eta(K_\eta)\bar K_\eta \mathrm{d}s \leqslant \bigg| \Re \int_{\RR^d} \mathrm{i} s_1 |s|_\eta^{-\gamma} \bar K_\eta \mathrm{d}s  \bigg| =  \bigg| \Re \int_{\RR^d} s_1 |s|_\eta^{-\gamma} \mathrm{Im} K_\eta  \mathrm{d}s\bigg| .$$
Now, since $\mathrm{Im} \ c_\eta = 0$, by the second item of Remark \ref{remarque c_eta}, we write $|\mathrm{Im} K_\eta|=| \mathrm{Im} H_\eta| \leqslant |H_\eta |$. Thus, by splitting the integral above into two parts, on $A :=\{|s_1|\leqslant s_0\}$ and on $A^c$,  we obtain:
$$ \bigg\| \nabla_s \bigg(\frac{\tilde K_\eta}{|s|_\eta^{-\gamma}}\bigg)|s|_\eta^{-\gamma} \bigg\|^2_2 \leqslant \big\| |s_1|^\frac{1}{2} |s|_\eta^{-\gamma} \big\|_{L^2(A)} \big\| |s_1|^\frac{1}{2} K_\eta \big\|_{L^2(A)} + \big\| |s|_\eta^{-\gamma} \big\|_{L^2(A^c)} \big\| s_1 H_\eta \big\|_{L^2(A^c)} . $$
Hence, returning to \eqref{s_1^(1/2)K_eta sur B 1}, we get:
\begin{align}\label{s_1^(1/2)K_eta sur B}
 \big\| |s_1|^{\frac{1}{2}} K_\eta \big\|^2_{L^2(B)} &\leqslant  \frac{1}{4} \big\| |s_1|^\frac{1}{2} K_\eta \big\|^2_{L^2(A)} + C s_0^3 \big\| |s|^{-\gamma} \big\|_{L^2(A^c)} \big\| s_1 H_\eta \big\|_{L^2(A^c)} \nonumber \\
 & + C \bigg( s_0^6 \big\| |s_1|^\frac{1}{2} |s|^{-\gamma} \big\|^2_{L^2(A)}+ |c_\eta-1| \big\| |s_1|^{\frac{1}{2}} |s|^{-\gamma} \big\|^2_{L^2(B)}\bigg)
\end{align}

\noindent \textbf{Estimation of $\big\| |s_1|^{\frac{1}{2}} K_\eta \big\|_{L^2(C)}$.} Recall that $C:=\{|s_1|\leqslant s_0\leqslant |s'|\}$. This step is identical to step 2 of the proof of Proposition \ref{propostion estimation L^2}. We start with estimate on $\big\| \zeta_{s_0} K_\eta \big\|^2_{L^2(\tilde C)}$. We have by using  inequality \eqref{Poincare2}
\begin{equation} \label{lader}\big\| \zeta_{s_0} K_\eta \big\|^2_{L^2(\tilde C)} \lesssim s_0^2 \bigg\| \pa_{s_1} \bigg(\frac{\zeta_{s_0}K_\eta}{|s|_\eta^{-\gamma}}\bigg)|s|_\eta^{-\gamma} \bigg\|^2_{L^2(\tilde C)} ,\end{equation}
where $\zeta_{s_0}(s_1):= \zeta\big( \frac{s_1}{s_0}\big)$, with $\zeta \in C^\infty(\RR)$ such that:  $0\leqslant \zeta \leqslant 1$,  $\zeta \equiv 1$ on $B(0,1)$ and $\zeta \equiv 0$ outside of $B(0,2)$, and where $\tilde C:=\{|s_1|\leqslant 2s_0\leqslant 2|s'|\}$. We have 
$$ \bigg\| \pa_{s_1} \bigg(\frac{\zeta_{s_0}K_\eta}{|s|_\eta^{-\gamma}}\bigg)|s|_\eta^{-\gamma} \bigg\|^2_{L^2(\tilde C)} \leqslant \bigg\| \nabla_s \bigg(\frac{\zeta_{s_0}K_\eta}{|s|_\eta^{-\gamma}}\bigg)|s|_\eta^{-\gamma} \bigg\|^2_{L^2(\tilde C)} =   \Re \int_{\tilde C} Q_\eta(\zeta_{s_0} K_\eta)\zeta_{s_0} \bar K_\eta \mathrm{d}s .$$
On the other hand, as in \eqref{Re int Q(zetaN)zetaN}
\begin{align*}
\Re \int_{\tilde C} Q_\eta(\zeta_{s_0} K_\eta)\zeta_{s_0} \bar K_\eta \mathrm{d}s &= \Re \int_{\tilde C} Q_\eta(K_\eta) \bar K_\eta \zeta_{s_0}^2  \mathrm{d}s +  \int_{\tilde C} |\zeta_{s_0} ' \bar K_\eta|^2 \mathrm{d}s  \\
&\leqslant \bigg| \Re \int_{\tilde C} -\mathrm{i} s_1 |s|_\eta^{-\gamma} \bar K_\eta \zeta_{s_0}^2 \mathrm{d}s \bigg| + \eta^{-\frac{2+\gamma}{3}}|b(0,\eta)|  \int_{\tilde C} \big| \Phi_\eta \bar K_\eta \zeta_{s_0}^2  \big| \mathrm{d}s \\
& \hspace{2cm} + \int_{\tilde C} |\zeta_{s_0}' K_\eta|^2 \mathrm{d}s .
\end{align*}
For the first term, we get
$$ \bigg| \Re \int_{\tilde C} -\mathrm{i} s_1 |s|_\eta^{-\gamma} \bar K_\eta \zeta_{s_0}^2 \mathrm{d}s \bigg| \leqslant \big\| |s_1|^{\frac{1}{2}} |s|_\eta^{-\gamma} \big\|_{L^2(\tilde C)} \big\||s_1|^{\frac{1}{2}} \zeta_{s_0} K_\eta \big\|_{L^2(\tilde C)} .$$
For the second term, recall that $\Phi_\eta(s):=c_{\gamma,d}\ \eta^{\frac{\gamma+2}{3}}|s|_\eta^{-\gamma-2}$, \eqref{Phi_eta < ...}, we get
\begin{align*}
\eta^{-\frac{2+\gamma}{3}}|b(0,\eta)|  \int_{\tilde C} \big| \Phi_\eta \bar K_\eta \zeta_{s_0}^2 \big| \mathrm{d}s &\lesssim s_0^{-\frac{5}{2}} |b(0,\eta)| \big\| |s_1|^{\frac{1}{2}} |s|_\eta^{-\gamma} \big\|_{L^2(\tilde C)} \| \zeta_{s_0} K_\eta \|_{L^2(\tilde C)} \\
&\lesssim s_0^{-\frac{5}{2}} |b(0,\eta)| \bigg(\big\| |s_1|^{\frac{1}{2}} |s|^{-\gamma} \big\|^2_{L^2(\tilde C)} + \| \zeta_{s_0} K_\eta \|^2_{L^2(\tilde C)}\bigg) .
\end{align*}
For the last term, since $\zeta_{s_0}'\equiv 0$ except on $\tilde C \setminus B^c := \{s_0\leqslant|s_1|\leqslant 2s_0 \leqslant 2|s'|\}$ where $|\zeta_{s_0}'(s_1)|\lesssim \frac{1}{s_0}$, and since $\tilde C \setminus B^c \subset A^c$, then 
$$
\int_{\tilde C} |\zeta_{s_0}' K_\eta|^2 \mathrm{d}s = \int_{\tilde C \setminus B^c} |\zeta_{s_0} ' K_\eta|^2 \mathrm{d}s \lesssim \frac{1}{s_0^2} \| K_\eta \|^2_{L^2(\tilde C \setminus B^c)} \lesssim \frac{1}{s_0^4} \big\| s_1 H_\eta \big\|^2_{L^2(A^c)} + \frac{1}{s_0^2} \| |s|_\eta^{-\gamma} \|^2_{L^2(A^c)} .
$$
Therefore,
\begin{align*}
\| \zeta_{s_0} K_\eta \|^2_{L^2(\tilde C)} &\lesssim s_0^2 \big\| |s_1|^{\frac{1}{2}} |s|^{-\gamma} \big\|_{L^2(\tilde C)} \big\||s_1|^{\frac{1}{2}} \zeta_{s_0} K_\eta \big\|_{L^2(\tilde C)} + s_0^{-\frac{1}{2}} |b(0,\eta)| \| \zeta_{s_0} K_\eta \|^2_{L^2(\tilde C)} \\
&+ s_0^{-\frac{1}{2}} |b(0,\eta)| \big\| |s_1|^{\frac{1}{2}} |s|^{-\gamma} \big\|^2_{L^2(\tilde C)} + \frac{1}{s_0^2} \big\| s_1 H_\eta \big\|^2_{L^2(A^c)} +  \| |s|^{-\gamma} \|^2_{L^2(A^c)}
\end{align*}
Since $|b(0,\eta)| \lesssim \eta^{\frac{2}{3}}$, thanks to \eqref{lim int eta^(1/3) v_1 M_eta M} and \eqref{b(lambda,eta)}, then for $\eta$ small enough we get
\begin{align*}
\| \zeta_{s_0} K_\eta \|^2_{L^2(\tilde C)} &\lesssim s_0^2 \big\| |s_1|^{\frac{1}{2}} |s|^{-\gamma} \big\|_{L^2(\tilde C)} \big\||s_1|^{\frac{1}{2}} \zeta_{s_0} K_\eta \big\|_{L^2(\tilde C)} + \frac{1}{s_0^2} \big\| s_1 H_\eta \big\|^2_{L^2(A^c)} \\
&+ s_0^{-\frac{1}{2}} |b(0,\eta)| \big\| |s_1|^{\frac{1}{2}} |s|^{-\gamma} \big\|^2_{L^2(\tilde C)} +  \| |s|^{-\gamma} \|^2_{L^2(A^c)} .
\end{align*}
Now, by \eqref{lader}, we get
$$ \big\| |s_1|^{\frac{1}{2}} \zeta_{s_0} K_\eta \big\|^2_{L^2(\tilde C)} \leqslant 2s_0 \big\| \zeta_{s_0} K_\eta \big\|^2_{L^2(\tilde C)} \lesssim s_0^3 \bigg\| \pa_{s_1} \bigg(\frac{\zeta_{s_0}K_\eta}{|s|_\eta^{-\gamma}}\bigg)|s|_\eta^{-\gamma} \bigg\|^2_{L^2(\tilde C)} .$$
Then,
\begin{align*}
\big\||s_1|^{\frac{1}{2}} \zeta_{s_0} K_\eta \big\|^2_{L^2(\tilde C)} &\lesssim s_0^3 \big\| |s_1|^{\frac{1}{2}} |s|^{-\gamma} \big\|_{L^2(\tilde C)} \big\||s_1|^{\frac{1}{2}} \zeta_{s_0} K_\eta \big\|_{L^2(\tilde C)} + \frac{1}{s_0} \big\| s_1 H_\eta \big\|^2_{L^2(A^c)} \\
&+ s_0^{\frac{1}{2}} |b(0,\eta)| \big\| |s_1|^{\frac{1}{2}} |s|^{-\gamma} \big\|^2_{L^2(\tilde C)} +  s_0\| |s|^{-\gamma} \|^2_{L^2(A^c)}
\end{align*}
Finally, since $\big\||s_1|^{\frac{1}{2}} K_\eta \big\|^2_{L^2(C)} \leqslant \big\||s_1|^{\frac{1}{2}} \zeta_{s_0} K_\eta \big\|^2_{L^2(\tilde C)}$, we get:
\begin{equation}\label{s_1^(1/2)K_eta sur C}
\big\||s_1|^{\frac{1}{2}} K_\eta \big\|^2_{L^2(C)} \lesssim \frac{1}{s_0} \big\| s_1 H_\eta \big\|^2_{L^2(A^c)} + s_0^6 \big\| |s_1|^{\frac{1}{2}} |s|^{-\gamma} \big\|^2_{L^2(\tilde C)}  + s_0 \| |s|^{-\gamma} \|^2_{L^2(A^c)}
\end{equation}
\textbf{Conclusion:} Since $A=B \cup C$ then, by summing the two inequalities \eqref{s_1^(1/2)K_eta sur B} and \eqref{s_1^(1/2)K_eta sur C} we find
\begin{align*}
\big\||s_1|^{\frac{1}{2}} K_\eta \big\|^2_{L^2(A)} &\leqslant \frac{1}{4} \big\| |s_1|^\frac{1}{2} K_\eta \big\|^2_{L^2(A)} +  \frac{C}{s_0} \big\| s_1 H_\eta \big\|^2_{L^2(A^c)} + C(s_0^6+|c_\eta-1|) \big\| |s_1|^{\frac{1}{2}} |s|^{-\gamma} \big\|^2_{L^2(\tilde A)}   \\
& + C s_0^3 \| |s|^{-\gamma} \|^2_{L^2(A^c)} 
\end{align*}
Hence,
\begin{equation}\label{est. sur A en fonction de A^c}
\big\||s_1|^{\frac{1}{2}} K_\eta \big\|^2_{L^2(A)} \lesssim \frac{1}{s_0} \big\| s_1 H_\eta \big\|^2_{L^2(A^c)} + (s_0^6+|c_\eta-1|) \big\| |s_1|^{\frac{1}{2}} |s|^{-\gamma} \big\|^2_{L^2(\tilde A)} + s_0^3 \big\| |s|^{-\gamma} \big\|^2_{L^2(A^c)} ,
\end{equation}
where $\tilde A:=\{|s_1|\leqslant 2s_0\} $. So it remains to estimate $\big\| s_1 H_\eta \big\|_{L^2(A^c)}$, where $A^c:=\{|s_1|\geqslant s_0\}$. \\

\noindent \textbf{Estimation of $\big\| s_1 H_\eta \big\|_{L^2(A^c)}$.} We have $\big\| s_1 H_\eta \big\|_{L^2(A^c)} \leqslant \big\| s_1 \chi_{s_0} H_\eta \big\|_{L^2(D)}$, where $\chi_{s_0}(s_1):= \chi\big( \frac{s_1}{s_0}\big)$, with $\chi \in C^\infty(\RR)$ such that $0 \leqslant \chi \leqslant 1$, $\chi \equiv 0$ on $B(0,\frac{1}{2})$ and $\chi \equiv 1$ outside $B(0,1)$ and where $D:=\{|s_1|\geqslant \frac{s_0}{2}\}$. Now, integrating the equation of $H_\eta$ against $ s_1 \bar H_\eta \chi_{s_0}^2$ and take the imaginary part, we obtain:
\begin{equation}\label{int_D s_1 zeta H_eta}
\int_D \big| s_1 \chi_{s_0} H_\eta \big|^2 \mathrm{d}s = -\mathrm{Im} \bigg(\int_D Q_\eta(H_\eta) s_1 \bar H_\eta \chi_{s_0}^2 \mathrm{d}s\bigg) - \eta^{-\frac{\gamma+2}{3}}\mathrm{Im} \bigg( b(0,\eta) \int_D \Phi_\eta s_1 \bar H_\eta \chi_{s_0}^2 \mathrm{d}s\bigg) 
\end{equation}
Let's start with the second term which is simpler, by \eqref{Phi_eta < ...} we have,
\begin{align*}
\eta^{-\frac{\gamma+2}{3}}\bigg|\mathrm{Im} \bigg( b(0,\eta) \int_D \Phi_\eta s_1 \bar H_\eta \chi_{s_0}^2 \mathrm{d}s\bigg) \bigg| &\lesssim \frac{1}{s_0^2} |b(0,\eta)| \| |s|_\eta^{-\gamma} \|_{L^2(D)}\big\| s_1 \chi_{s_0} H_\eta \big\|_{L^2(D)}  \\
&\lesssim \frac{1}{s_0^2} |b(0,\eta)| \bigg(\| |s|^{-\gamma} \|^2_{L^2(D)} + \big\| s_1 \chi_{s_0} H_\eta \big\|^2_{L^2(D)} \bigg) .
\end{align*}
For the first term, we will proceed exactly as for $E^\eta_1$ (first step in the proof of the Proposition \ref{propostion estimation L^2}). By integration by parts,  we write
\begin{align*}
\bigg| \mathrm{Im} \int_D Q_\eta(H_\eta) s_1 \bar H_\eta \chi_{s_0}^2 \mathrm{d}s \bigg| &= \bigg| \mathrm{Im} \int_D \pa_{s_1}\bigg(\frac{H_\eta}{|s|_\eta^{-\gamma}} \bigg)|s|_\eta^{-\gamma}\chi_{s_0}\big[\chi_{s_0}\bar H_\eta + 2 s_1 \bar H_\eta \chi_{s_0}'\big] \mathrm{d}s \bigg|  \\
&\leqslant \bigg\| \pa_{s_1}\bigg(\frac{H_\eta}{|s|_\eta^{-\gamma}} \bigg)|s|_\eta^{-\gamma}\chi_{s_0} \bigg\|_{L^2(D)} \bigg( \| \chi_{s_0} H_\eta\|_{L^2(D)} + 2 \big\| s_1 \chi_{s_0}' H_\eta \big\|_{L^2(D)} \bigg) \\
&\leqslant \frac{1}{2} \bigg\| \pa_{s_1}\bigg(\frac{H_\eta}{|s|_\eta^{-\gamma}} \bigg)|s|_\eta^{-\gamma}\chi_{s_0} \bigg\|^2_{L^2(D)}  + \frac{2}{s_0^2} \| s_1\chi_{s_0} H_\eta\|^2_{L^2(D)}   \\
&+  \frac{s_0}{2} \bigg\| \pa_{s_1}\bigg(\frac{H_\eta}{|s|_\eta^{-\gamma}} \bigg)|s|_\eta^{-\gamma} \chi_{s_0} \bigg\|^2_{L^2(D)} + \frac{1}{s_0}\big\| s_1 \chi_{s_0}' H_\eta \big\|^2_{L^2(D)} .
\end{align*}
Now, since $\chi_{s_0}'\equiv 0$ except on: $D \setminus A^c:=\{\frac{s_0}{2} \leqslant |s_1|\leqslant s_0\} \subset A$ where $|\chi_{s_0}'(s_1)|\lesssim \frac{1}{s_0}$, and since $|H_\eta|\leqslant |K_\eta|+ |s|_\eta^{-\gamma}$ then,
\begin{align*}
\big\| s_1 \chi_{s_0}' H_\eta \big\|^2_{L^2(D)} \lesssim  \| H_\eta \|^2_{L^2(D \setminus A^c)} &\lesssim \| K_\eta \|^2_{L^2(D \setminus A^c)} +\| |s|_\eta^{-\gamma} \|^2_{L^2(D \setminus A^c)} \\
&\lesssim \frac{1}{s_0} \bigg( \big\| |s_1|^{\frac{1}{2}} K_\eta \big\|^2_{L^2(A)} + \big\| |s_1|^{\frac{1}{2}} |s|_\eta^{-\gamma} \big\|^2_{L^2(A)} \bigg) .
\end{align*}
Therefore,
\begin{align}\label{Im int_D Q_eta...}
\bigg| \mathrm{Im} \int_D Q_\eta(H_\eta) s_1 \bar H_\eta \chi_{s_0}^2 \mathrm{d}s \bigg| &\leqslant  s_0 \bigg\| \pa_{s_1}\bigg(\frac{H_\eta}{|s|_\eta^{-\gamma}} \bigg)|s|_\eta^{-\gamma}\chi_{s_0} \bigg\|^2_{L^2(D)} + \frac{2}{s_0^2} \| s_1\chi_{s_0} H_\eta\|^2_{L^2(D)}  \nonumber \\
&+  \frac{C}{s_0^2} \bigg( \big\| |s_1|^{\frac{1}{2}} K_\eta \big\|^2_{L^2(A)} + \big\| |s_1|^{\frac{1}{2}} |s|^{-\gamma} \big\|^2_{L^2(A)} \bigg) .
\end{align}
Let us now deal with the term $\big\| \pa_{s_1}\big(\frac{H_\eta}{|s|_\eta^{-\gamma}} \big)|s|_\eta^{-\gamma}\chi_{s_0} \big\|_{L^2(D)}$. By an integration by parts, we can see that  
$$ \bigg\| \nabla_s\bigg(\frac{H_\eta}{|s|_\eta^{-\gamma}} \bigg)|s|_\eta^{-\gamma}\chi_{s_0} \bigg\|^2_{L^2(D)} = \Re \int_D \bigg[ Q_\eta(H_\eta)\bar H_\eta \chi_{s_0}^2 - 2 \chi_{s_0} \chi_{s_0}' \frac{\bar H_\eta}{|s|_\eta^{-\gamma}}\pa_{s_1}\bigg(\frac{H_\eta}{|s|_\eta^{-\gamma}} \bigg)|s|_\eta^{-2\gamma} \bigg] \mathrm{d}s . $$
Therefore,
\begin{align*}
\bigg\| \pa_{s_1}\bigg(\frac{H_\eta}{|s|_\eta^{-\gamma}} \bigg)|s|_\eta^{-\gamma}\chi_{s_0} \bigg\|^2_{L^2(D)} &\leqslant \bigg|\Re \int_D  Q_\eta(H_\eta)\bar H_\eta \chi_{s_0}^2 \mathrm{d}s\bigg| \\
&+ 2 \|\chi_{s_0}' H_\eta \|_{L^2(D)} \bigg\| \pa_{s_1}\bigg(\frac{H_\eta}{|s|_\eta^{-\gamma}} \bigg)|s|_\eta^{-\gamma}\chi_{s_0} \bigg\|_{L^2(D)}  \\
&\leqslant \bigg|\Re \int_D  Q_\eta(H_\eta)\bar H_\eta \chi_{s_0}^2 \mathrm{d}s\bigg| \\
&+ 2 \|\chi_{s_0}' H_\eta \|^2_{L^2(D)} +\frac{1}{2} \bigg\| \pa_{s_1}\bigg(\frac{H_\eta}{|s|_\eta^{-\gamma}} \bigg)|s|_\eta^{-\gamma}\chi_{s_0} \bigg\|^2_{L^2(D)} .
\end{align*}
Which implies that,
\begin{align*}
\bigg\| \pa_{s_1}\bigg(\frac{H_\eta}{|s|_\eta^{-\gamma}} \bigg)|s|_\eta^{-\gamma}\chi_{s_0} \bigg\|^2_{L^2(D)} &\leqslant 2\bigg|\Re \int_D  Q_\eta(H_\eta)\bar H_\eta \chi_{s_0}^2 \mathrm{d}s\bigg| + 4 \|\chi_{s_0}' H_\eta \|^2_{L^2(D\setminus A^c)} \\
&\lesssim \frac{1}{s_0^2} \bigg(|b(0,\eta)| \| |s|_\eta^{-\gamma} \|_{L^2(D)}\|\chi_{s_0} H_\eta \|_{L^2(D)} + \| H_\eta \|^2_{L^2(D \setminus A^c)} \bigg) \\
&\lesssim \frac{1}{s_0^3} |b(0,\eta)| \bigg(\big\| s_1 \chi_{s_0} H_\eta \big\|^2_{L^2(D)} + \| |s|^{-\gamma} \|^2_{L^2(D)} \bigg)  \\
&+ \frac{1}{s_0^3} \bigg( \big\| |s_1|^{\frac{1}{2}} K_\eta \big\|^2_{L^2(A)} + \big\| |s_1|^{\frac{1}{2}} |s|^{-\gamma} \big\|^2_{L^2(A)} \bigg) .
\end{align*} 
Thus, injecting this last inequality into \eqref{Im int_D Q_eta...} we obtain
\begin{align*}
\bigg| \mathrm{Im} \int_D Q_\eta(H_\eta) s_1 \bar H_\eta \chi_{s_0}^2 \mathrm{d}s \bigg| &\lesssim  \frac{1}{s_0^2}\bigg[ (1+|b(0,\eta)|) \| s_1 \chi_{s_0} H_\eta\|^2_{L^2(D)} + \big\| |s_1|^{\frac{1}{2}} K_\eta \big\|^2_{L^2(A)} \\
&+   \big\| |s_1|^{\frac{1}{2}} |s|^{-\gamma} \big\|^2_{L^2(A)} +  |b(0,\eta)| \| |s|^{-\gamma} \|^2_{L^2(D)} \bigg] ,
\end{align*}
and going back to \eqref{int_D s_1 zeta H_eta}, using the fact that $|b(0,\eta)| \lesssim 1$ by Remark \ref{rmq sur thm d'existence}, we get
\begin{align*}
\| s_1 \chi_{s_0} H_\eta \|^2_{L^2(D)}&\lesssim   \frac{1}{s_0^2}\bigg[ \| s_1 \chi_{s_0} H_\eta\|^2_{L^2(D)} + \big\| |s_1|^{\frac{1}{2}} K_\eta \big\|^2_{L^2(A)} + \big\| |s_1|^{\frac{1}{2}} |s|^{-\gamma} \big\|^2_{L^2(A)} +  \| |s|^{-\gamma} \|^2_{L^2(D)} \bigg] .
\end{align*}
Finally, for $s_0$ large enough, the term $\frac{1}{s_0^2}\| s_1 \chi_{s_0} H_\eta\|^2_{L^2(D)}$ is absorbed and we obtain thanks to the inequality $\| s_1 H_\eta \|^2_{L^2(A^c)} \leqslant \| s_1 \chi_{s_0} H_\eta \|^2_{L^2(D)}$:
\begin{equation}\label{s_1 H_eta sur A^c}
\| s_1 H_\eta \|^2_{L^2(A^c)} \lesssim \frac{1}{s_0^2} \bigg( \big\| |s_1|^{\frac{1}{2}} K_\eta \big\|^2_{L^2(A)} +  \big\| |s_1|^{\frac{1}{2}} |s|^{-\gamma} \big\|^2_{L^2(A)} +  \| |s|^{-\gamma} \|^2_{L^2(D)} \bigg) .
\end{equation}
Now, by injecting  inequality \eqref{s_1 H_eta sur A^c} into \eqref{est. sur A en fonction de A^c}, we get 
$$
\big\||s_1|^{\frac{1}{2}} K_\eta \big\|^2_{L^2(A)} \leqslant  C\bigg(\frac{1}{s_0^3} \big\||s_1|^{\frac{1}{2}} K_\eta \big\|^2_{L^2(A)} + s_0^6 \big\| |s_1|^{\frac{1}{2}} |s|_\eta^{-\gamma} \big\|^2_{L^2(\tilde A)} + s_0^3 \big\| |s|_\eta^{-\gamma} \big\|^2_{L^2(D)} \bigg) .
$$
Where we used the fact that $A^c \subset D$ and $|c_\eta-1| \lesssim 1$ by Remark \ref{remarque c_eta}. Finally, for $s_0$ large enough, the norm $\big\||s_1|^{\frac{1}{2}} K_\eta \big\|^2_{L^2(A)}$ which appears in the right hand side of the previous inequality is absorbed, from where:
\begin{equation}\label{s_1^(1/2)K_eta sur A}
\big\||s_1|^{\frac{1}{2}} K_\eta \big\|^2_{L^2(A)} \lesssim s_0^6 \big\| |s_1|^{\frac{1}{2}} |s|^{-\gamma} \big\|^2_{L^2(\tilde A)} + s_0^3 \big\| |s|^{-\gamma} \big\|^2_{L^2(A^c)} \lesssim 1 ,
\end{equation}
since for $\gamma \in ]\frac{d}{2},\frac{d+1}{2}[$: $|s_1|^{\frac{1}{2}}|s|^{-\gamma} \in L^2(\tilde A)$ and $|s|^{-\gamma} \in L^2(A^c)$.
From the inequality \eqref{s_1 H_eta sur A^c} we deduce that $\big\||s_1|^{\frac{1}{2}} K_\eta \big\|^2_{L^2(A)} \lesssim 1$ implies that $\| s_1 H_\eta \|^2_{L^2(A^c)} \lesssim 1$. Thus we obtain  \eqref{estimation de  H_eta pour gamma<(d+1)/2}. \\

\noindent \textbf{2.} Recall that $\tilde K_\eta:=H_\eta-c_\eta |s|_\eta^{-\gamma}$ satisfies the orthogonality condition \eqref{condition d'orthogonalite} of Hardy-Poincar\'e inequality \eqref{ineg.  Hardy-Poincare} and that $\mathrm{Im} c_\eta = 0$. It follows that, $\mathrm{Im} \tilde K_\eta = \mathrm{Im} K_\eta=\mathrm{Im} H_\eta$, so by \eqref{ineg.  Hardy-Poincare}, we get on the one hand
$$ \bigg\| \frac{\mathrm{Im} H_\eta}{|s|_\eta} \bigg\|^2_{L^2(A)} = \bigg\| \frac{\mathrm{Im} \tilde K_\eta}{|s|_\eta} \bigg\|^2_{L^2(A)} \leqslant \bigg\| \frac{\tilde K_\eta}{|s|_\eta} \bigg\|^2_2  \lesssim_{\gamma,d} \bigg\| \nabla_s\bigg(\frac{\tilde K_\eta}{|s|_\eta^{-\gamma}}\bigg) |s|_\eta^{-\gamma} \bigg\|^2_2 = \bigg\| \nabla_s\bigg(\frac{K_\eta}{|s|_\eta^{-\gamma}}\bigg) |s|_\eta^{-\gamma} \bigg\|^2_2 .$$
On the other hand,
$$\bigg\| \nabla_s\bigg(\frac{K_\eta}{|s|_\eta^{-\gamma}}\bigg) |s|_\eta^{-\gamma} \bigg\|^2_2 = \int_{\RR^d} Q_\eta(K_\eta) \bar K_\eta \mathrm{d}s = \Re \bigg(-\mathrm{i} \int_{\RR^d} s_1 |s|_\eta^{-\gamma} \bar K_\eta \mathrm{d}s \bigg)-\eta^{-\frac{\gamma+2}{3}}b(0,\eta)\int_{\RR^d} \Phi_\eta \bar K_\eta \mathrm{d}s .$$
We have: $b(0,\eta)\int_{\RR^d} \Phi_\eta \bar K_\eta \mathrm{d}s \geqslant 0$. Indeed, by performing the change of variable $s=\eta^{\frac{1}{3}}v$, we obtain:
$$\int_{\RR^d} \Phi_\eta \bar K_\eta \mathrm{d}s = \eta^{\frac{\gamma+d}{3}}\int_{\RR^d} \Phi \bar N_ {0,\eta} \mathrm{d}v = \eta^{\frac{\gamma+d}{3}} \overline{b(0,\eta)} .$$
Now, since $\ds \Re \bigg(-i \int_{\RR^d} s_1 |s|_\eta^{-\gamma} \bar K_\eta \mathrm{d}s \bigg)= \int_{\RR^d} s_1 |s|_\eta^{-\gamma} \mathrm{Im} K_\eta \mathrm{d}s$, then we write:
\begin{align}
\bigg\| \nabla_s\bigg(\frac{K_\eta}{|s|_\eta^{-\gamma}}\bigg) |s|_\eta^{-\gamma} \bigg\|^2_2 &\leqslant \int_{\RR^d} |s_1 |s|_\eta^{-\gamma} \mathrm{Im} K_\eta| \mathrm{d}s \nonumber \\
&= \int_{A} |s_1 |s|_\eta^{-\gamma} \mathrm{Im} K_\eta| \mathrm{d}s + \int_{A^c} |s_1 |s|_\eta^{-\gamma} \mathrm{Im} H_\eta| \mathrm{d}s  \label{grad(K_eta/m_eta) m_eta 1} \\
&\leqslant \big\| s_1 |s|^{1-\gamma} \big\|_{L^2(A)} \bigg\| \frac{\mathrm{Im} \tilde K_\eta}{|s|_\eta} \bigg\|_{L^2(A)} \nonumber \\  
& \quad + s_0^{-\frac{1}{2}} \big\| |s_1|^{\frac{1}{2}} |s|^{-\gamma} \big\|_{L^2(A^c)} \|s_1 H_\eta\|_{L^2(A^c)} . \label{grad(K_eta/m_eta) m_eta}
\end{align}
It remains to estimate the norm $\|s_1 H_\eta\|_{L^2(A^c)}$. Recall that $D:=\{|s_1|\geqslant \frac{s_0}{2}\}$. We start by estimating $\|s_1 \chi_{s_0} H_\eta\|_{L^2(D)}$. We have; as before; the two equalities:
$$
\int_D \big| s_1 \chi_{s_0} H_\eta \big|^2 \mathrm{d}s = -\mathrm{Im} \bigg(\int_D Q_\eta(H_\eta) s_1 \bar H_\eta \chi_{s_0}^2 \mathrm{d}s\bigg) - \eta^{-\frac{\gamma+2}{3}}\mathrm{Im} \bigg( b(0,\eta) \int_D \Phi_\eta s_1 \bar H_\eta \chi_{s_0}^2 \mathrm{d}s\bigg)
$$
and 
$$\bigg| \mathrm{Im} \int_D Q_\eta(H_\eta) s_1 \bar H_\eta \chi_{s_0}^2 \mathrm{d}s \bigg| = \bigg| \mathrm{Im} \int_D \pa_{s_1}\bigg(\frac{H_\eta}{|s|_\eta^{-\gamma}} \bigg)\chi_{s_0}|s|_\eta^{-\gamma}\big[\chi_{s_0}\bar H_\eta + 2 s_1 \bar H_\eta \chi_{s_0}'\big] \mathrm{d}s \bigg|.$$
The term on the right in the first equality is treated in the same way as before and we have:
\begin{align}\label{Im b(0,eta) int_D Phi...}
\eta^{-\frac{\gamma+2}{3}}\bigg|\mathrm{Im} \bigg( b(0,\eta) \int_D \Phi_\eta s_1 \bar H_\eta \chi_{s_0}^2 \mathrm{d}s\bigg) \bigg| &\lesssim s_0^{-\frac{5}{2}} |b(0,\eta)| \big\| |s_1|^{\frac{1}{2}} |s|_\eta^{-\gamma} \big\|_{L^2(D)}\big\| s_1 \chi_{s_0} H_\eta \big\|_{L^2(D)} \nonumber \\
&\lesssim s_0^{-\frac{5}{2}} |b(0,\eta)| \bigg(\big\| |s_1|^{\frac{1}{2}} |s|_\eta^{-\gamma} \big\|^2_{L^2(D)} \nonumber \\
& \hspace{2.5cm} + \big\| s_1 \chi_{s_0} H_\eta \big\|^2_{L^2(D)}  \bigg)
\end{align}
For the left term in the first equality we write: 
$$\bigg| \mathrm{Im} \int_D \pa_{s_1}\bigg(\frac{H_\eta}{|s|_\eta^{-\gamma}} \bigg)\chi_{s_0}|s|_\eta^{-\gamma}\big[\chi_{s_0}\bar H_\eta + 2 s_1 \bar H_\eta \chi_{s_0}'\big] \mathrm{d}s \bigg| \leqslant I^\eta_1 + I^\eta_2.$$
where
$$
I^\eta_1 :=\bigg| \mathrm{Im} \int_D \chi_{s_0}\bar H_\eta \pa_{s_1}\bigg(\frac{H_\eta}{|s|_\eta^{-\gamma}} \bigg)\chi_{s_0}|s|_\eta^{-\gamma} \mathrm{d}s \bigg|\quad \mbox{and}\quad I^\eta_2:= \bigg| \mathrm{Im} \int_D s_1 \chi_{s_0}\bar H_\eta \pa_{s_1}\bigg(\frac{H_\eta}{|s|_\eta^{-\gamma}} \bigg)|s|_\eta^{-\gamma} \chi_{s_0}' \mathrm{d}s \bigg|.$$
Then we write
\begin{align}\label{Im int_D terme 1}
I^\eta_1 &\leqslant \| \chi_{s_0} H_\eta \|_{L^2(D)} \bigg\| \pa_{s_1}\bigg(\frac{H_\eta}{|s|_\eta^{-\gamma}} \bigg)\chi_{s_0}|s|_\eta^{-\gamma} \bigg\|_{L^2(D)} \nonumber \\
&\leqslant \frac{1}{s_0}\| s_1 \chi_{s_0} H_\eta \|_{L^2(D)} \bigg\| \nabla_s\bigg(\frac{H_\eta}{|s|_\eta^{-\gamma}} \bigg) |s|_\eta^{-\gamma} \bigg\|_2  \nonumber \\
&\leqslant \frac{1}{2s_0} \bigg( \| s_1 \chi_{s_0} H_\eta \|^2_{L^2(D)} + \bigg\| \nabla_s \bigg(\frac{H_\eta}{|s|_\eta^{-\gamma}} \bigg) |s|_\eta^{-\gamma} \bigg\|^2_2 \bigg) ,
\end{align}
and
\begin{align}\label{Im int_D terme 2}
I^\eta_2&\leqslant \| s_1 \chi_{s_0} H_\eta \|_{L^2(D)} \| \chi_{s_0}' \|_{L^\infty(D\setminus A^c)} \bigg\| \pa_{s_1}\bigg(\frac{H_\eta}{|s|_\eta^{-\gamma}} \bigg)|s|_\eta^{-\gamma} \bigg\|_{L^2(D)} \nonumber \\
&\lesssim \frac{1}{s_0} \bigg( \| s_1 \chi_{s_0} H_\eta \|^2_{L^2(D)} + \bigg\| \nabla_s\bigg(\frac{H_\eta}{|s|_\eta^{-\gamma}} \bigg) |s|_\eta^{-\gamma} \bigg\|^2_2 \bigg) .
\end{align}
Hence, by  inequalities \eqref{Im b(0,eta) int_D Phi...}, \eqref{Im int_D terme 1} and \eqref{Im int_D terme 2} to estimate $\| s_1 \chi_{s_0} H_\eta \|^2_{L^2(D)}$, and by inequality \eqref{grad(K_eta/m_eta) m_eta} to estimate $\big\| \nabla_s\big(\frac{H_\eta}{|s|_\eta^{-\gamma}} \big) |s|_\eta^{-\gamma} \big\|^2_2$, we get
\begin{align*}
\| s_1 \chi_{s_0} H_\eta \|^2_{L^2(D)} &\lesssim \frac{1}{s_0} \big\| s_1 |s|^{1-\gamma} \big\|_{L^2(A)} \bigg\| \frac{\mathrm{Im} \tilde K_\eta}{|s|_\eta} \bigg\|_{L^2(A)} + s_0^{-\frac{3}{2}} \big\| |s_1|^{\frac{1}{2}} |s|^{-\gamma} \big\|_{L^2(A^c)} \|s_1 H_\eta\|_{L^2(A^c)} \\
&+ \frac{1}{s_0}\| s_1 \chi_{s_0} H_\eta \|^2_{L^2(D)} + s_0^{-\frac{5}{2}} |b(0,\eta)| \bigg(\big\| |s_1|^{\frac{1}{2}} |s|^{-\gamma} \big\|^2_{L^2(D)} + \big\| s_1 \chi_{s_0} H_\eta \big\|^2_{L^2(D)}  \bigg)  \\
&\lesssim \frac{1}{s_0} \big\| s_1 |s|^{1-\gamma} \big\|_{L^2(A)} \bigg\| \frac{\mathrm{Im} \tilde K_\eta}{|s|_\eta} \bigg\|_{L^2(A)} + \frac{1}{s_0}\| s_1 \chi_{s_0} H_\eta \|^2_{L^2(D)} + \frac{1}{s_0^{3/2}}\big\| |s_1|^{\frac{1}{2}} |s|^{-\gamma} \big\|^2_{L^2(D)}.
\end{align*}
Hence, for $s_0$ large enough and since $\| s_1 H_\eta \|^2_{L^2(A^c)} \leqslant \| s_1 \chi_{s_0} H_\eta \|^2_{L^2(D)}$: 
\begin{equation}\label{s_1 H_eta sur A^c pour gamma>(d+1)/2}
\| s_1 H_\eta \|^2_{L^2(A^c)} \lesssim \frac{1}{s_0}  \bigg\| \frac{\mathrm{Im} H_\eta}{|s|_\eta} \bigg\|^2_{L^2(A)} + \frac{1}{s_0} \big\| s_1 |s|^{1-\gamma} \big\|^2_{L^2(A)} + s_0^{-\frac{3}{2}}\big\| |s_1|^{\frac{1}{2}} |s|^{-\gamma} \big\|^2_{L^2(D)} .
\end{equation}
So, going back to \eqref{grad(K_eta/m_eta) m_eta} and using inequality $ab \leqslant C a^2 + \frac{b^2}{4C}$, we get:
\begin{align*}
\bigg\| \frac{\mathrm{Im} H_\eta}{|s|_\eta} \bigg\|^2_{L^2(A)} \leqslant C \bigg\| \nabla_s\bigg(\frac{K_\eta}{|s|_\eta^{-\gamma}} \bigg) |s|_\eta^{-\gamma} \bigg\|^2_2 &\leqslant \frac{1}{4}\bigg\| \frac{\mathrm{Im} H_\eta}{|s|_\eta}  \bigg\|^2_{L^2(A)} + C^2 \big\| s_1 |s|^{1-\gamma} \big\|^2_{L^2(A)} \\ 
&+ \frac{C}{2}\big\| |s_1|^{\frac{1}{2}} |s|^{-\gamma} \big\|^2_{L^2(D)}+ \frac{C}{2s_0} \| s_1 H_\eta \|^2_{L^2(A^c)} 
\end{align*}
Thus, by \eqref{s_1 H_eta sur A^c pour gamma>(d+1)/2} we obtain
$$ \bigg\| \frac{\mathrm{Im} H_\eta}{|s|_\eta} \bigg\|^2_{L^2(A)} \leqslant \bigg(\frac{1}{4}+\frac{C'}{s_0^2}\bigg)\bigg\| \frac{\mathrm{Im} H_\eta}{|s|_\eta} \bigg\|^2_{L^2(A)}+C'\bigg(\big\| s_1 |s|^{1-\gamma} \big\|^2_{L^2(A)}+ \big\| |s_1|^{\frac{1}{2}} |s|^{-\gamma} \big\|^2_{L^2(D)}\bigg) .$$
Finally, for $s_0$large enough
\begin{equation}\label{Im H_eta/|s|_eta sur A}
\bigg\| \frac{\mathrm{Im} H_\eta}{|s|_\eta} \bigg\|^2_{L^2(A)} \lesssim \big\| s_1 |s|^{1-\gamma} \big\|^2_{L^2(A)}+ \big\| |s_1|^{\frac{1}{2}} |s|^{-\gamma} \big\|^2_{L^2(D)} \lesssim 1 .
\end{equation}
By \eqref{s_1 H_eta sur A^c pour gamma>(d+1)/2}, it follows that $\| s_1 H_\eta \|^2_{L^2(A^c)} \lesssim 1$. Note that for $\gamma \in ]\frac{d+1}{2},\frac{d+4}{2}[$ we have: 
$$ s_1|s|_\eta |s|_\eta^{-\gamma} \leqslant s_1 |s|^{1-\gamma} \in L^2(A) \ \mbox{ and } \ |s_1|^{\frac{1}{2}} |s|_\eta^{-\gamma} \leqslant |s_1|^{\frac{1}{2}}|s|^{-\gamma} \in L^2(A^c) .$$ 
Hence  inequality \eqref{estimation de H_eta pour gamma>(d+1)/2} holds.
\epl

\noindent  The third lemma contains some complementary estimates on the rescaled solution. 
\begin{lemma}[Complementary estimates]\label{estimations complementaires}  For all $\eta \in [0,\eta_0]$ and for all $\gamma \in (\frac{d}{2},\frac{d+4}{2})$, the following estimate holds
\begin{equation}\label{estimations de H_eta - |s|_eta^(-gamma)}
\bigg\|  \frac{H_\eta - c_\eta |s|_\eta^{-\gamma}}{|s|_\eta} \bigg\|_{L^2(\RR^d)}  \lesssim_{\gamma,d} \bigg\|  \nabla_s\bigg(\frac{H_\eta}{|s|_\eta^{-\gamma}}\bigg)|s|_\eta^{-\gamma}\bigg\|_{L^2(\RR^d)} \lesssim 1 .
\end{equation}
\end{lemma}

\noindent  The last one gives the formula of the diffusion coefficient.
\begin{lemma}\label{lim kappa_eta}
We have the following limit:
\begin{equation}\label{lim kappa_eta = kappa}
\underset{\eta \rightarrow 0}{\lim}\ \mathrm{i} \eta^{1-\alpha} \int_{\{|v_1| \geqslant R\}} v M_{0,\eta}(v) M(v) \mathrm{d}v = - 2 \int_{0}^\infty \int_{\RR^{d-1}} s_1|s|^{-\gamma} \mathrm{Im} H_0(s) \mathrm{d}s ,
\end{equation}
where $H_0$ is the unique solution to \eqref{eq de H_0} satisfying  conditions \eqref{condition H_0}.
\end{lemma}

\noindent \bpl \ref{estimations complementaires}.  We have by the Hardy-Poincar\'e inequality and  inequality \eqref{grad(K_eta/m_eta) m_eta 1}
\begin{align*}
\Lambda_{\gamma,d} \bigg\| \frac{H_\eta - c_\eta |s|_\eta^{-\gamma}}{|s|_\eta} \bigg\|^2_{L^2(\RR^d)} &\leqslant  \bigg\| \nabla_s\bigg(\frac{H_\eta-|s|_\eta^{-\gamma}}{|s|_\eta^{-\gamma}}\bigg) |s|_\eta^{-\gamma} \bigg\|^2_{L^2(\RR^d)} \\
&\leqslant  \int_{\RR^d} |s_1| |s|_\eta^{-\gamma} |\mathrm{Im} H_\eta| \mathrm{d}s \\
&= \int_{\{|s_1| \leqslant s_0\}} |s_1| |s|_\eta^{-\gamma} |\mathrm{Im} H_\eta| \mathrm{d}s + \int_{\{|s_1|\geqslant s_0\}} |s_1| |s|_\eta^{-\gamma} |\mathrm{Im} H_\eta| \mathrm{d}s .
\end{align*}
\textbf{Case 1: $\gamma \in (\frac{d}{2},\frac{d+1}{2})$.} By Cauchy-Schwarz and inequality \eqref{estimation de  H_eta pour gamma<(d+1)/2} of Lemma \ref{grandes vitesses} we get
$$ \int_{\{|s_1| \leqslant s_0\}} |s_1| s|_\eta^{-\gamma} |\mathrm{Im} H_\eta| \mathrm{d}s \leqslant \big\| |s_1|^{\frac{1}{2}} |s|_\eta^{-\gamma} \big\|_{L^2(\{|s_1| \leqslant s_0\})} \big\| |s_1|^{\frac{1}{2}} \mathrm{Im} H_\eta \big\|_{L^2(\{|s_1| \leqslant s_0\})}  \lesssim 1 $$
and 
$$  \int_{\{|s_1|\geqslant s_0\}} |s_1| |s|_\eta^{-\gamma} |\mathrm{Im} H_\eta| \mathrm{d}s \leqslant \big\| |s|_\eta^{-\gamma} \big\|_{L^2(\{|s_1|\geqslant s_0\})} \big\| s_1 H_\eta \big\|_{L^2(\{|s_1|\geqslant s_0\})} \lesssim 1 .
$$
\textbf{Case 2: $\gamma \in (\frac{d+1}{2},\frac{d+4}{2})$.} Similary, by Cauchy-Schwarz and inequality \eqref{estimation de H_eta pour gamma>(d+1)/2} we get
$$ \int_{\{|s_1| \leqslant s_0\}} |s_1| |s|_\eta^{-\gamma} |\mathrm{Im} H_\eta| \mathrm{d}s \leqslant \big\| s_1 |s|_\eta^{1-\gamma} \big\|_{L^2(\{|s_1| \leqslant s_0\})} \bigg\| \frac{\mathrm{Im} H_\eta}{|s|_\eta} \bigg\|_{L^2(\{|s_1| \leqslant s_0\})}  \lesssim 1$$
and 
$$  \int_{\{|s_1|\geqslant s_0\}} |s_1| |s|_\eta^{-\gamma} |\mathrm{Im} H_\eta| \mathrm{d}s \leqslant \big\| |s|_\eta^{-\gamma} \big\|_{L^2(\{|s_1|\geqslant s_0\})} \big\| s_1 H_\eta \big\|_{L^2(\{|s_1|\geqslant s_0\})} \lesssim 1 .
$$
This completes the proof of the lemma.
\epl

\noindent \bpl \ref{lim kappa_eta}.  First of all, since $\bar{M}_{0,\eta}(-v_1,v')=M_{0,\eta}(v_1,v')$ and $M(-v_1,v')=M(v_1,v')$ for all $v_1\in\RR$ and for all $v' \in \RR^{d-1}$, thus
 $$ \mathrm{i}  \int_{\{|v_1| \geqslant R\}} v_1 M_{0,\eta}(v) M(v) \mathrm{d}v = - 2 \int_{\{v_1 \geqslant R\}} v_1 \mathrm{Im} M_{0,\eta}(v) M(v) \mathrm{d}v . $$
Then, in order to compute the limit
 $$\underset{\eta \rightarrow 0}{\lim} \  \mathrm{i} \eta^{\frac{d+1-2\gamma}{3}} \int_{\{|v_1| \geqslant R\}} v_1 M_{0,\eta}(v) M(v) \mathrm{d}v = - 2  \underset{\eta \rightarrow 0}{\lim}\  \eta^{\frac{d+1-2\gamma}{3}} \int_{\{v_1 \geqslant R\}} v_1 \mathrm{Im} M_{0,\eta}(v) M(v) \mathrm{d}v ,$$
we proceed to a change of variable $v=\eta^{-\frac{1}{3}}s$, which means that we  compute
$$ 
\underset{\eta \rightarrow 0}{\lim} \ \int_{\{|s_1|\geqslant \eta^{\frac{1}{3}}R\}}  s_1 |s|_\eta^{-\gamma} \mathrm{Im} H_\eta(s) \mathrm{d}s .
$$
For that purpose, we use the \textit{weak-strong} convergence in the Hilbert space $L^2(\RR_+\times\RR^{d-1})$. 

\noindent The estimates of Lemma \ref{grandes vitesses} imply that the sequence $\mathsf{H}_\eta$ defined by
$$
\mathsf{H}_\eta(s):=  \left\{\begin{array}{l} s_1^{\frac{1}{2}} \mathrm{Im} H_\eta (s),  \qquad \ \ \qquad \gamma\in (\frac{d}{2},\frac{d+1}{2}],  0< s_1 \leqslant s_0 , \\
 \\
 |s|_\eta^{-1} \mathrm{Im} H_\eta (s) , \ \qquad \quad \ \gamma\in (\frac{d+1}{2},\frac{d+4}{2}), 0 < s_1 \leqslant s_0 , \\
 \\
 s_1 \mathrm{Im} H_\eta (s)  \  \quad \mbox{  for all } \ \ \gamma\in (\frac{d}{2},\frac{d+4}{2}) \mbox{ and } s_1 \geqslant s_0 ,
\end{array}\right.
$$
is bounded in $L^2(\RR^d)$, uniformly with respect to $\eta$, which implies that $\mathsf{H}_\eta$ converges weakly in $L^2(\RR^d)$,  up to a subsequence.  Let's identify this limit that we denote by $\mathsf{H}_0 \in L^2(\RR^d)$.  We have on the one hand, $H_\eta$ converges to $H_0$ in $\mathcal{D}'(\RR^d\setminus\{0\})$. Indeed,  recall that $H_\eta$ satisfies the equation
$$
\bigg[-\Delta_s+\frac{\gamma(\gamma-d+2)}{|s|^2_\eta} +\mathrm{i} s_1 \bigg] H_{\eta}(s)= \eta^{\frac{2}{3}}\frac{\gamma(\gamma+2)}{|s|^4_\eta} H_\eta(s) - \eta^{-\frac{2+\gamma}{3}} b(0,\eta) \Phi_\eta(s).
$$
Let $\varphi \in \mathcal{D}(\RR^d \setminus \{0\})$. Then, by integrating the previous equation against $\varphi$, we obtain:
\begin{align*}
\int_{\RR^d \setminus \{0\}} \bigg[-\Delta_s+ \frac{\gamma(\gamma-d+2)}{|s|^2_\eta} +\mathrm{i} s_1 \bigg]\varphi(s) H_\eta(s)\mathrm{d}s &= \eta^{\frac{2}{3}} \int_{\RR^d \setminus \{0\}} \frac{\gamma(\gamma+2)}{|s|^4_\eta}   \varphi(s)H_{\eta}(s)\mathrm{d}s \\
&- \eta^{-\frac{2+\gamma}{3}} b(0,\eta)  \int_{\RR^d \setminus \{0\}} \Phi_\eta(s) \varphi(s)\mathrm{d}s .
\end{align*}
Thanks to the uniform bound \eqref{estimations de H_eta - |s|_eta^(-gamma)} and since $\Phi_\eta(s) \lesssim \eta^{\frac {2+\gamma}{3}}|s|^{-2-\gamma}$ and $b(0,\eta)\to 0$ then, passing to the limit when $\eta$ goes to $0$ in the last equality, we obtain that $H_\eta$ converges to $H_0$ in $\mathcal{D}'(\RR^d \setminus \{0\})$,  solution to the equation
\begin{equation}\label{equation H_0}
\bigg[-\Delta_s+  \frac{\gamma(\gamma-d+2)}{|s|^2}+ is_1 \bigg]H_0(s)= 0 .
\end{equation}
Moreover, for all $\gamma \in (\frac{d} {2},\frac{d+4}{2})$,  the function $H_\eta$ satisfies the estimate 
$$ \bigg\| \frac{H_\eta - c_\eta |s|_\eta^{-\gamma}}{|s|_\eta}\bigg\|_{L^2(\{|s_1|\leqslant s_0\})} + \|s_1 H_\eta\|_{L^2(\{|s_1|\geqslant s_0\})} \lesssim 1 ,$$
thanks to  inequality \eqref{estimations de H_eta - |s|_eta^(-gamma)} and the first point of Lemma \ref{grandes vitesses} for $\gamma \in (\frac{d}{2},\frac{d+1}{2}]$, and thanks to the second point of Lemma \ref{grandes vitesses} for $\gamma \in (\frac{d+1}{2},\frac{d+4}{2})$. Therefore $H_0$ satisfies the estimate
$$ \bigg\| \frac{H_0 - |s|^{-\gamma}}{|s|}\bigg\|_{L^2(\{|s_1|\leqslant s_0\})}  + \|s_1 H_0\|_{L^2(\{|s_1|\geqslant s_0\})} \lesssim 1 . $$
Now,  $\|s_1 H_0\|_{L^2(\{|s_1|\geqslant s_0\})} \lesssim 1$ implies that $H_0 \in L^2(\{|s_1|\geqslant 1\})$ and $\big\| \frac{H_0 - |s|^{-\gamma}}{|s|}\big\|_{L^2(\{|s_1|\leqslant s_0)} \lesssim 1 $ implies that $H_0(s) \underset{0}{ \sim} |s|^{-\gamma}$, a different behaviour near zero would make the latter norm infinite.
These two conditions imply that $H_0$ is the unique solution of  equation \eqref{equation H_0}.
Thanks to the uniqueness of this limit,  the whole sequence $ \mathsf{H}_\eta$ converges  weakly to 
$$
\mathsf{H}_0(s):=  \left\{\begin{array}{l} s_1^{\frac{1}{2}} \mathrm{Im} H_0 (s),  \ \ \ \quad \qquad \gamma\in (\frac{d}{2},\frac{d+1}{2}],  0< s_1 \leqslant s_0 , \\
 \\
 |s|^{-1} \mathrm{Im} H_0 (s) , \ \ \quad \quad \ \gamma\in (\frac{d+1}{2},\frac{d+4}{2}), 0 < s_1 \leqslant s_0 , \\
 \\
 s_1 \mathrm{Im} H_0 (s)  \  \quad \mbox{  for all } \ \gamma\in (\frac{d}{2},\frac{d+4}{2}) \mbox{ and } s_1 \geqslant s_0 .
\end{array}\right.
$$
Finally, we conclude by passing to the limit in the scalar product $\langle \mathsf{H}_\eta, \mathsf{I}_\eta\rangle$, where $\mathsf{I}_\eta$ definded by 
$$
\mathsf{I}_\eta:=  \left\{\begin{array}{l} s_1^{\frac{1}{2}} |s|_\eta^{-\gamma},  \quad \qquad \gamma\in (\frac{d}{2},\frac{d+1}{2}], \ 0<|s_1| \leqslant s_0 , \\
 \\
 s_1 |s|_\eta^{1-\gamma} , \quad \quad \ \ \gamma\in (\frac{d+1}{2},\frac{d+4}{2}), \ 0< s_1 \leqslant s_0 , \\
 \\
|s|_\eta^{-\gamma} ,   \qquad \qquad  \gamma\in (\frac{d}{2},\frac{d+4}{2}) , \ s_1 \geqslant s_0 ,
\end{array}\right.  
$$
converges strongly in $L^2(\RR_+\times\RR^{d-1})$ to 
$$
\mathsf{I}_0:=  \left\{\begin{array}{l} s_1^{\frac{1}{2}} |s|^{-\gamma},  \quad \qquad \gamma\in (\frac{d}{2},\frac{d+1}{2}], \ 0<s_1 \leqslant s_0 , \\
 \\
 s_1 |s|^{1-\gamma} , \quad \quad \ \ \gamma\in (\frac{d+1}{2},\frac{d+4}{2}), \ 0<s_1 \leqslant s_0 , \\
 \\
|s|^{-\gamma} ,   \qquad \qquad  \gamma\in (\frac{d}{2},\frac{d+4}{2}) , \ s_1 \geqslant s_0 .
\end{array}\right.  
$$
Hence  limit \eqref{lim kappa_eta = kappa} holds true.  
\epl

\noindent \bpp \ref{vp}. By doing an expansion in $\lambda$ for $B$ and by Proposition \ref{contrainte}, we get
$$
B(\lambda,\eta)=\eta^{-\frac{2}{3}}b(\lambda, \eta)=\eta^{-\frac{2}{3}}b(0,\eta)+ \lambda \int_{\RR^d} M_{0,\eta} M \mathrm{d}v+ O(\lambda^2).
$$ 
Then, for $\lambda=\tilde{\lambda}(\eta)$ and since $B(\tilde{\lambda}(\eta),\eta)=0$, we obtain:
$$ 
\tilde{\lambda}(\eta)= - \eta^{-\frac{2}{3}}b(0,\eta)\bigg(\int_{\RR^d} M_{0,\eta} M \mathrm{d}v\bigg)^{-1}+o\big(\eta^{-\alpha}b(0,\eta)\big) ,
$$ 
which implies that 
$$\eta^{-\alpha}\mu(\eta)=\eta^{\frac{2}{3}-\alpha}\tilde{\lambda}(\eta)= -\eta^{-\alpha}b(0,\eta)\bigg(\int_{\RR^d} M_{0,\eta} M \mathrm{d}v\bigg)^{-1} .
$$
By \eqref{lim int  M_eta M} and \eqref{lim kappa_eta = kappa}, 
$$\underset{\eta \rightarrow 0}{\lim} \int_{\RR^d} M_{0,\eta}(v)M(v)\mathrm{d}v = \|M\|^2_2 = C_\beta^{-2} $$ 
and 
$$ \underset{\eta \rightarrow 0}{\lim}\ \eta^{-\alpha}b(0,\eta) = 2C_\beta^2 \int_0^\infty \int_{\RR^{d-1}} s_1|s|^{-\gamma}\mathrm{Im} H_0(s)\mathrm{d}s'\mathrm{d}s_1 $$ 
respectively. Hence,  $\underset{\eta \rightarrow 0}{\lim} \eta^{-\alpha}\mu(\eta) = \kappa$.
For $\eta\in [-\eta_0, 0]$, the symmetry $ \mu(\eta) = \overline \mu (-\eta) $ holds by complex conjugation on the equation.  So it remains to prove the positivity of $\kappa$.  By integrating the equation of $M_{\eta}:=M_{\tilde{\lambda} (\eta),\eta}$ against $\bar M_\eta$ we obtain:
$$ \int_{\RR^d} \bigg| \nabla_v\bigg(\frac{M_\eta}{M}\bigg) \bigg|^2 M^2 \mathrm{d}v + \mathrm{i} \eta \int_{\RR^d} v_1 |M_\eta|^2 \mathrm{d}v = \mu(\eta) \int_{\RR^d} |M_\eta|^2 \mathrm{d}v . $$
Now, taking the real part and using the equality $\ds \mu(\eta)\|M_\eta\|_2^2=\kappa \eta^{\alpha}\big(1+o(\eta^\alpha)\big)$ we get:
\begin{equation}\label{kappa >0}
\int_{\RR^d} \bigg| \nabla_v\bigg(\frac{M_\eta}{M}\bigg) \bigg|^2 M^2 \mathrm{d}v  = \kappa \eta^{\alpha}\big(1+o(\eta^\alpha)\big) .
\end{equation}
Therefore, multiplying this last equality by $\eta^{-\alpha}$ and performing the change of variable $v=\eta^{-\frac{1}{3}}s$ we obtain:
$$ \int_{\RR^d} \bigg| \nabla_s\bigg(\frac{H_\eta}{|s|_\eta^{-\gamma}}\bigg) \bigg|^2 |s|_\eta^{-2\gamma} \mathrm{d}s  = \kappa \big(1+o_\eta(1)\big) . $$
Thus, $ \kappa\geqslant 0$. If $\kappa=0$ then, 
$$ \int_{\RR^d} \bigg| \nabla_s\bigg(\frac{H_0}{|s|^{-\gamma}}\bigg) \bigg|^2 |s|^{-2\gamma} \mathrm{d}s \leqslant \liminf  \int_{\RR^d} \bigg| \nabla_s\bigg(\frac{H_\eta}{|s|_\eta^{-\gamma}}\bigg) \bigg|^2 |s|_\eta^{-2\gamma} \mathrm{d}s = 0 .$$
Therefore, $H_0=|s|^{-\gamma}$. Which leads to a contradiction since $H_0$ is solution to equation \eqref{eq de H_0}.  Hence, the proof of Proposition \ref{vp} is complete.
\epp

\noindent \bpt \ref{main}. The existence of the eigen-solution $(\mu(\eta),M_\eta)$ is given by Proposition \ref{contrainte}. Limit \eqref{M_eta--->M dans H^1} follows from inequality \eqref{estimation de N_lambda,eta dans L^2} for $|\lambda| = |\tilde{\lambda}(\eta)| \lesssim \eta^{\frac{2\gamma-d}{3}} \underset{\eta \to 0}{\longrightarrow} 0$, thanks to \eqref{mu(eta)}, with  limit \eqref{M_lambda,eta-M_0-->0 dans H_0} obtained by Theorem \ref{thm d'existence}. Finally, the second point of Theorem \ref{main} is given by Proposition \ref{vp}.

\section{Derivation of the fractional diffusion equation}
The goal of this section is to prove Theorem \ref{main2}. The proof was taken from Section 3 in \cite{LebPu} and adapted for the dimension $d$. \\

\noindent Let's start by defining the two weighted $L^p$ spaces, $L^p_{F^{1-p}}(\RR^d)$ and $Y^p_F(\RR^{2d})$: 
$$ L^p_{F^{1-p}}(\RR^d) := \left\{f: \mathbb{R}^d\rightarrow \mathbb{R},\int_{\mathbb{R}^{d}} |f|^p \ F^{1-p}\ \mathrm{d}v< \infty\right\} \  \mbox{ and } \  Y^p_F(\RR^{2d}):=L^p\big(\RR^d,L^p_{F^{1-p}}(\RR^d)\big). $$
Recall that our goal is to show that the solution $f^\eps$ of the Fokker-Planck equation \eqref{fp-theta} converges; weakly star in $L^\infty\left([0, T], L^2_{F^{-1}}(\mathbb{R}^{2})\right)$; towards $\rho(t,x) F(v)$ when $\eps$ goes to $0$, where $\rho$ is the solution of the following \emph{fractional diffusion} equation
\begin{equation}
\partial_t\rho +\kappa (-\Delta_x)^{\frac{\beta-d+2}{6}}\rho =0,\quad \rho(0,x)=\int_{\RR^{d}} f_0 \mathrm{d}v .
 \end{equation}
\begin{remark}
Note that we will work with the Fourier transform of $\rho$ and we will prove that 
$ \hat\rho(t,\xi)=\int_{\RR^d} e^{-\mathrm{i}x \cdot \xi} \rho(t,x)\mathrm{d}x$ satisfies
\begin{equation}\label{diff}
\partial_t\hat \rho+\kappa|\xi|^{\frac{\beta-d+2}{3}}\hat\rho=0.
\end{equation}
\end{remark}
\subsection{A priori estimates}
We start by recalling the following compactness lemma.
\begin{lemma}$\cite{LebPu}, \cite{NaPu}$\label{norm}
For initial datum $f_0 \in Y^p_F(\RR^{2d})$ where $p \geqslant 2$ and a positive time $T$.
\begin{enumerate}
\item The solution $f^\varepsilon$ of \eqref{fp-theta} is bounded in $L^\infty\big([0,T];Y^p_F(\RR^{2d})\big)$ uniformly with respect to $\eps$ since it satisfies 
\begin{equation}\label{estimate1.1}\|f^\varepsilon(T)\|^p_{Y^p_F(\RR^{2d})}+\frac{p(p-1)}{\theta(\varepsilon)}\int_0^T\int_{\mathbb{R}^{2d}}\bigg|\nabla_v\bigg(\frac{f^\varepsilon}{F}\bigg)\bigg|^2 \bigg|\frac{f^\varepsilon}{F}\bigg|^{p-2}\ F\ \mathrm{d}v \mathrm{d}x \mathrm{d}t \leqslant \|f_0\|^p_{Y^p_F(\RR^{2d})}.
\end{equation}
\item The density $\rho^\varepsilon(t,x)=\int_{\mathbb{R}^d}f^\varepsilon \ \mathrm{d}v$ is such that
\begin{equation} 
\|\rho^\varepsilon(t)\|^p_{p} \leqslant C_\beta^{-2(p-1)} \|f_0\|^p_{Y^p_F(\RR^{2d})}\quad \mbox{ for all } \quad t\in [0,T].
\end{equation}
\item Up to a subsequence, the density $\rho^\varepsilon$ converges weakly star in $L^\infty \big([0,T]; L^p(\mathbb{R}^d)\big)$ to $\rho$.
\item 
Up to a subsequence, the sequence $f^\varepsilon$ converges weakly star in $L^\infty \big([0,T]; Y^p_F(\RR^{2d})\big)$ to the function $f=\rho(t,x)F(v)$.
\end{enumerate}
\end{lemma}
As a consequence, we have the following estimate:
\begin{corollary}$\cite{LebPu}$ Let $F = C_\beta^2 M^2$ with $M=(1+|v|^2)^{-\frac{\gamma}{2}}$ and $\beta=2\gamma\in (d,d+4)$.
Let $f^\eps$ solution to \eqref{fp-theta} with $\theta(\varepsilon)=\varepsilon^{\frac{2\gamma-d+2}{3}}$.  Assume that $\|f_0/F\|_{\infty} \leqslant C$. Then  $g^\eps=f^\eps F^{-\frac{1}{2}}$ satisfies the following estimate
\begin{equation}\label{estblan}
\int_0^T \int_{\RR^d}\left(\int_{\RR^d} \big|g^\eps-\rho^\eps F^{\frac{1}{2}}\big|^2\mathrm{d}v \right)^{\frac{2\gamma-d+2}{2\gamma-d}} ds \mathrm{d}x \leqslant C \eps^\frac{2\gamma-d+2}{3}.
\end{equation}
\end{corollary}
\bp
Recall the Nash type inequality \cite{CGGR,rw,BaBaCaGu}:  for any $h$ such that $\int h F\mathrm{d}v=0$, we have
\begin{equation}\label{Nash}
\int_{\mathbb R^d} h^2 F\mathrm{d}v\leqslant C\left(\int_{\mathbb R^d} \big|\nabla_v h \big|^2 F\mathrm{d}v\right)^\frac{2\gamma-d}{2\gamma-d+2}(\|h\|_\infty^2)^\frac{2}{2\gamma-d+2} \ .
\end{equation}
Define $h=g^\eps F^{-\frac{1}{2}}-\rho^\eps =\frac{f^\eps}{F}-\rho^\eps$, define $\alpha=\frac{2\gamma-d+2}{3}$.
Observe that  from $\|f\|_{L^p_{F^{1-p}}(\RR^{2d})}=\big\|\frac{f}{F}\big\|_{L^p_F}$ 
and  Lemma \ref{norm}, formula \eqref{estimate1.1}, we have
$$
\|h_0\|_{L^\infty}=\lim_{p\rightarrow\infty } \|h_0\|_{L^p_{F^{1-p}}(\RR^{2d})}\geqslant\lim_{p \rightarrow \infty } \|h\|_{L^p_{F^{1-p}}(\RR^{2d})}
\geqslant\|h\|_{L^\infty} \ .
$$
 Thus by Lemma \ref{norm}, formula \eqref{estimate1.1} taking $p=2$, we get 
$$\begin{array}{rcl}
\dps\int_0^T\int _{\mathbb R^d} \left(\int_{\mathbb R^d} |g^\eps-\rho^\eps F^{\frac{1}{2}}|^2\mathrm{d}v\right)^{\frac{2\gamma-d+2}{2\gamma-d}} \mathrm{d}s\mathrm{d}x
&=& \dps\int_0^T\int _{\mathbb R^d} \left(\int_{\mathbb R^d} h^2 F\mathrm{d}v\right)^{\frac{2\gamma-d+2}{2\gamma-d}}\mathrm{d}s\mathrm{d}x \\ 
&\leqslant& 
\dps C\int_0^T\int _{\mathbb R^d} \left(\int_{\mathbb R^d}\big|\nabla_v h \big|^2 F\mathrm{d}v\right)(\|h\|_\infty^2)^\frac{2}{2\gamma-d} \mathrm{d}s \mathrm{d}x\\
&\dps\leqslant &C\dps \int_0^T\int _{\mathbb R^d} \bigg(\int_{\RR^d} \bigg|\nabla_v \bigg(\frac{f^\eps}{F}\bigg)\bigg|^2F \mathrm{d}v\bigg) \mathrm{d}s \mathrm{d}x
\leqslant C \eps^\alpha.
\end{array}$$
\ep
\subsection{Weak limit and proof of Theorem \ref{main2}}
By solving equation (\ref{rescaled}), we write
$$\hat g^\eps(t,\xi,v)=e^{-t\theta(\eps)\mathcal{L}_\eta}\hat g(0,\xi,v) ,
$$
which gives going back to the rescaled space variable $x$
$$
 g^\eps(t,x,v)=\frac{1}{(2\pi)^d} \int_{\RR^d} e^{\mathrm{i} x \cdot \xi}\hat g^\eps(t,\xi,v)\mathrm{d}\xi  .
$$
Our purpose is to pass to the limit when $\eps\rightarrow 0$.\\
Recall that $f^\eps(t,x,v) \geqslant 0$ and $\int f^\eps(t,x,v)\mathrm{d}x\mathrm{d}v=\int f_0(x,v)\mathrm{d}x\mathrm{d}v$ for all $t\geqslant 0$.\\
Let $\hat \rho^\eps (t,\xi)= \int_{\RR^d} e^{-\mathrm{i}x \cdot \xi}\rho^\eps (t,x)\mathrm{d}x \ $ be the Fourier transform in $x$ of 
$\rho^\eps = \int_{\RR^d} f^\eps \mathrm{d}v= \int_{\RR^d} g^\eps F^{\frac{1}{2}} \mathrm{d}v$.

\begin{proposition}\label{lem:equilibrium}
For  all $\xi\in \RR^d$,  $\hat\rho^\varepsilon(\cdot,\xi)$ converges  to $\hat \rho(\cdot,\xi)$, 
unique solution  to the ode
\begin{equation}\label{fracdifeq}
\partial_t\hat\rho+\kappa|\xi|^\alpha \hat \rho=0, \quad \hat \rho_0 = \int_{\RR^d}\hat  f_0 \mathrm{d}v \ .
\end{equation}
\end{proposition}

\noindent \bp Let $\xi \in \RR^d $ and let
$M_\eta$ be the unique solution in $L^2(\RR^d,\CC)$ of $\ \mathcal{L}_\eta(M_\eta)=\mu(\eta)M_\eta \ $ given in Theorem \ref{main}. One has
$$
\begin{array}{rcl}
\dps \frac{\mathrm{d}}{\mathrm{d}t}\int_{\RR^d}  \hat g^\eps(t,\xi,v) M_\eta \mathrm{d}v&=&\dps  \int_{\RR^d} \partial_t\hat g^\eps  M_\eta \mathrm{d}v
= -\eps^{-\alpha}\int_{\RR^d} \mathcal{L}_\eps(\hat g^\eps) M_\eta \mathrm{d}v\\
&=&
\dps -\eps^{-\alpha}\int_{\RR^d} \hat g^\eps    \mathcal{L}_\eps(M_\eta) \mathrm{d}v = -\eps^{-\alpha}\mu(\eta )\int_{\RR^d} \hat g^\eps (t,\xi,v) M_\eta \mathrm{d}v  .  
\end{array}
$$
Therefore one has, with $F^\eps (t,x)=C_\beta \int_{\RR^d}   g^\eps(t,x,v) M_\eta \mathrm{d}v$, 
\begin{equation}\label{gl9}
 \hat F^\eps (t,\xi) =e^{-t\eps^{-\alpha}\mu(\eps |\xi|)}  \hat F^\eps (0,\xi) , \quad \forall t \geqslant 0.
\end{equation} 
By Theorem \ref{main}, we have
$ \eps^{-\alpha}\mu(\eps |\xi|)\rightarrow \kappa  |\xi|^\alpha$. 
Moreover, the following limit holds true:
\begin{equation}\label{gl11}
\forall \xi \in \RR^d, \quad \hat F^\eps (0,\xi) = C_\beta \int_{\RR^d} \hat g^\eps(0,\xi,v)  M_\eta \mathrm{d}v  \rightarrow  \hat \rho_0(\xi) \ . 
\end{equation}
The verification of \eqref{gl11} is easy. One has $ \hat g^\eps(0,v,\xi)=  \hat f_0(v,\xi) F^{-\frac{1}{2}}(v) = \frac{\hat f_0(v,\xi)}{C_\beta M(v)}$ and $M_\eta \to M$ in $L^2(\RR^d)$ thanks to \eqref{M_eta--->M dans H^1}. Thus, \eqref{gl11} holds true by Cauchy-Schwarz inequality by writing: 
$$ \bigg| C_\beta \int_{\RR^d} \hat g^\eps(0,\xi,v)  M_\eta \mathrm{d}v - \hat \rho_0(\xi) \bigg| \leqslant C_\beta \bigg(\int_{\RR^d} \frac{f_0^2}{F} \mathrm{d}v \bigg)^{\frac{1}{2}}\bigg(\int |M_\eta-M|^2\mathrm{d}v \bigg)^{\frac{1}{2}} .$$
It remains to verify 
\begin{equation}\label{gl10}
\forall \xi \in {\RR^d}, \quad C_\beta\int_{\RR^d} \hat g^\eps(t,\xi,v)  M_\eta \mathrm{d}v  \longrightarrow  \hat \rho(t,\xi)  
\quad \text{in} \ \mathcal D'\big(]0,\infty[\times\RR^{d}\big).
\end{equation}
By \eqref{gl9} and \eqref{gl11}, for all $\xi\in \RR^d $ and $t\geqslant 0$, one has
$\underset{\eps\rightarrow 0}{\lim}\ \hat F^\eps (t,\xi)=e^{-t\kappa |\xi|^\alpha}\hat\rho_0(\xi)$, thus \eqref{gl10}
will be consequence of the weaker
\begin{equation}\label{gl10bis}
\quad C_\beta\int_{\RR^d}  g^\eps(t,x,v)  M_\eta \mathrm{d}v  \rightarrow   \rho(t,x)  
\quad \text{in} \ \mathcal D'\big(]0,\infty[\times\RR^d\big) \ .
\end{equation}

\noindent Let us now verify \eqref{gl10bis}. For that purpose, we write
$$
C_\beta \int_{\RR^d}  g^\eps  M_\eta \mathrm{d}v  -\rho= 
C_\beta  \int_{\RR^d} ( g^\eps- \rho^\eps F^{\frac{1}{2}})M_\eta \mathrm{d}v +  \rho^\eps \int_{\RR^d} (C_\beta M_\eta-F^{\frac{1}{2}})F^{\frac{1}{2}} \mathrm{d}v 
+\rho^\eps -\rho \ .
$$
By using (\ref{estblan}) and the first point of Theorem \ref{main},  limit \eqref{M_eta--->M dans H^1}, we pass to the limit.
The proof of Proposition \ref{lem:equilibrium} is complete.
\ep

\noindent \textbf{\bpt \ref{main2}. }From the two last items in Lemma \ref{norm}, we have just to prove that 
for any given $\xi$, the Fourier transform $\hat\rho(t,\xi)$ of the weak limit
$\rho(t,y)$, is solution  of  equation \eqref{diff}, which is precisely Proposition \ref{lem:equilibrium}.  \ept

\ni \textbf{Acknowledgment.} The authors would like to thank Gilles Lebeau for the fruitful discussions, as well as Pierre Rapha\"el for drawing their attention to Herbert Koch's article.

{\footnotesize
\bibliography{biblio-FP}
}
\bibliographystyle{abbrv}

\end{document}